\documentclass[preprint,3p,onecolumn,times]{elsarticle}

\usepackage{natbib}
\usepackage{graphicx}
\usepackage{subfigure}
\usepackage{url}
\usepackage{amssymb}
\usepackage{amsmath}
\usepackage{amsfonts}
\usepackage{wrapfig}
\usepackage{afterpage}
\usepackage{float}
\usepackage{algorithm}
\usepackage{colortbl}
\usepackage{multirow}
\usepackage[dvipsnames*,svgnames]{xcolor}
\usepackage{lscape}
\usepackage{threeparttable}
\definecolors{Fuchsia}
\usepackage{threeparttable}
\usepackage[noend]{algpseudocode}
\usepackage[toc,page]{appendix}
\usepackage[bookmarks=true]{hyperref}
\graphicspath{{figures/}}

\newtheorem{thm}{Theorem}
\newtheorem{cor}{Corollary}

\def\eb{\begin{equation}}
\def\ee{\end{equation}}

\newcommand{\bsig}{\mbox{\boldmath$\sigma$}}
\newcommand{\bveps}{\mbox{\boldmath$\varepsilon$}}

\newcommand{\btau}{\mbox{\boldmath$\tau$}}

\newcommand{\bpsi}{\mbox{\boldmath$\psi$}}
\newcommand{\bPsi}{\mbox{\boldmath$\Psi$}}
\newcommand{\bPhi}{\mbox{\boldmath$\Phi$}}

 \def\bB{\mathbf{B}}

\def\bD{\mathbf{D}}

\def\bf{\mathbf{f}} \def\bF{\mathbf{F}}

\def\bI{\mathbf{I}}

\def\bk{\mathbf{k}} \def\bK{\mathbf{K}}

\def\bM{\mathbf{M}}

\def\bP{\mathbf{P}}
 \def\bR{\mathbf{R}}

\def\bu{\mathbf{u}} \def\bU{\mathbf{U}}
 \def\bN{\mathbf{N}}

\def\bx{\mathbf{x}} \def\bX{\mathbf{X}}
\def\bv{\mathbf{v}}

\def\calX{{\cal X}}

\def\bF{\mathbf{F}}

\def\bB{\mathbf{B}}

\def\bff{\mathbf{f}}
\def\bg{\mathbf{g}}

\def\bD{\mathbf{D}}

\def\bK{\mathbf{K}}

\def\bI{\mathbf{I}}
\def\bP{\mathbf{P}}
\def\bN{\mathbf{N}}

\def\bq{\mathbf{q}}
\def\bQ{\mathbf{Q}}
\def\bu{\mathbf{u}}
\def\u{\mathbf{u}}
\def\n{\mathbf{n}}
\def\bU{\mathbf{U}}
\def\bv{\mathbf{v}}

\def\bx{\mathbf{x}}
\def\bX{\mathbf{X}}
\def\x{\mathbf{x}}

\def\bt{\mathbf{t}}

\def\rdV{\mathrm{dV}}
\def\rd{\mathrm{d}}

\def\calM{\mathcal{M}}

\def\calV{\mathcal{X}}

\def\RR{\mathbb{R}}

\begin{document}
\begin{frontmatter}
\title{Analysis of Heterogeneous Structures of Non-separated Scales on Curved Bridge Nodes}
\author[cadaddress]{Ming Li\corref{mycorrespondingauthor}}
\cortext[mycorrespondingauthor]{Corresponding author: liming@cad.zju.edu.cn}
\address[cadaddress]{State Key Laboratory of CAD$\&$CG, Zhejiang University, Hangzhou, China}
\author[cadaddress]{Jingqiao Hu}


\begin{abstract}
  Numerically predicting the performance of heterogenous structures without scale separation represents a significant challenge to meet the critical requirements on computational scalability and efficiency -- adopting a mesh fine enough to fully account for the small-scale heterogeneities leads to prohibitive computational costs while simply ignoring these fine heterogeneities tends to drastically over-stiffen the structure's rigidity.

  This study proposes an approach to construct new material-aware shape (basis) functions per element on a coarse discretization of the structure with respect to each \emph{curved bridge nodes (CBNs)} defined along the elements' boundaries. Instead of formulating their derivation by solving a nonlinear optimization problem, the shape functions are constructed by building a map from the CBNs to the interior nodes and are ultimately presented in an explicit matrix form as a product of a B\'ezier interpolation transformation and a boundary-interior transformation.
  The CBN shape function accomodates more flexibility in closely capturing the coarse element's heterogeneity, overcomes the important and challenging issues of inter-element stiffness and displacement discontinuity across interface between coarse elements, and improves the analysis accuracy by orders of magnitude; they also meet the basic geometric properties of shape functions that avoid aphysical analysis results. Extensive numerical examples, including a 3D industrial example of billions of degrees of freedom, are also tested to demonstrate the approach's performance in comparison with results obtained from classical approaches.
\end{abstract}
\end{frontmatter}

\section{Introduction}
Heterogeneous structures comprise varied material properties at different locations within their interior and are found in different types of natural objects such as human bones or organs~\cite{khademhosseini2016decade}, or engineered alloys, polymers, reinforced composites~\cite{fratzl2007nature}. The numerical prediction of the physical performance of such heterogenous structures is a perpetual and fundamental issue in engineering design~\cite{Karel2017A}; however, it remains a significant challenge to develop elaborate numerical methods to meet the critical requirements on computational scalability and efficiency~\cite{1997olson,2013key,alexandersen2015topology,yvonnet2019scales,YUZH2020,raschi2021high}.
Classical finite element (FE) methods only capture properly the structures' behavior (only elasticity is studied here) if one adopts a mesh fine enough to account for the small-scale heterogeneities~\cite{2007Heterogeneous}, leading to prohibitive computational costs particulary for structures of highly complex geometries and material distributions. Simply ignoring these fine heterogeneities however tends to result in an important issue of \emph{inter-element stiffness}~\cite{nesme2009preserving}, which renders the structure deformation dramatically more rigid than in reality.

A possible strategy to address the issue is via parallel computation based on domain decomposition methods (DDM)~\cite{1991FarhatA,2006Domain,Spillane2016} or to significantly reduce the scope of the problem by using a coarse grid via geometric multigrid~\cite{2000GMG} or algebraic multigrid~\cite{K1983Algebraic,K2001A}. The efficacy of these methods are challenged by a loss of accuracy for structures containing large heterogeneities or high contrast of materials, particularly when the subdomain interface intersects the heterogeneities. We will not go further into the topic, and refer interested readers~\cite{2005Domain,2011goddeke,le2020coarse}.

Multiscale methods are being increasingly applied to predict the behavior of heterogenous structures. The analysis in such cases is usually achieved via two levels of FE simulations---macroscale and microscale---that use the analysis results on each microstructure in parallel to aid the prediction of the overall performance of the structure in the macroscale and vice versa. Numerical homogenization is a typical mean-field multiscale analysis approach that replaces each microstructure with an effective elasticity tensor using the calculation results from the microstructure analysis via the asymptotic approach~\cite{pinho2009asymptotic,chung2001asymptotic,andreassen2014determine} and the energy-based approach~\cite{sigmund1994materials,xia2015design}. However, the method is limited in its use of linear models. The multi-level FE method (FE$^2$) is another important multiscale approach that typically conducts FE analysis iteratively by transiting between fields (stress and strain) in the macroscale and microscale until convergence~\cite{1998Prediction,feyel2003multilevel,xia2014concurrent}. The FE$^2$ approach is able to more accurately capture microscopic heterogenous information, although at more expensive computational costs.
Both approaches of numerical homogenization and FE$^2$ are usually built on the assumption of scale separation, that is, the length scale of the microstructure is much lower than that of the structural length scale. The assumption is, however, no longer valid for the purpose of analysis of heterogeneous structures without scale separation, as studied here. Researchers have developed various approaches to address this issue, including the high-order computational homogenization~\cite{2002Multi,kouznetsova2004multi,yvonnet2019computational}, fiber-based homogenization~\cite{tognevi2016multi} or direct FE$^2$~\cite{tan2020direct}. A comprehensive literature review on FE$^2$ is referred to in~\cite{2014FE2Book}, and on multiscale in ~\cite{Karel2017A}.

Substructuring is also studied for the analysis of heterogeneous structures~\cite{wu2019topology,liu2020data}. It treats all the structures as a set of substructures connected by boundary nodes between the coarse elements, called \emph{super-elements}. Based on a local FE formulation of each coarse element, a matrix condensation strategy ultimately produces a linear equation about the super-elements, whose solution consequently yields the global solution in fine mesh. Substructuring is able to produce a high-accuracy solution, but faces two main challenges that prohibit its industrial applications. First, the local analysis problem per coarse mesh element involves solution computations to a very large number of linear equation systems and is costly. More importantly, in contrast to the fine-scale analysis problem, it produces a dense global stiffness matrix with more non-zero elements; see also the example in Fig.~\ref{fig:sparsity}. In addition, the substructuring approach is only applicable to linear problems. 

Constructing tailored material-aware shape (or "basis") functions has shown great promise for analysis of heterogeneous structure of non-separated scales; it is also called multiscale FE method~\cite{hou1997multiscale,efendiev2011multiscale}. These approaches substitute classical FE shape functions per coarse mesh element with newly constructed complex ones obtained from fine scale calculations. These approaches meet with two main challenges: closely capturing the coarse element's heterogeneity and maintaining the global solution continuity in the fine mesh. Most previous studies focus on the first challenge, and articulate the shape function construction as a spectral expansion~~\cite{hou1997multiscale,efendiev2011multiscale} or constrained nonlinear optimization problem~\cite{chen2018numerical,le2020coarse}.
%
%
More recently, Le at al developed a novel CMCM (Coarse Mesh Condensation Multiscale Method) approach for a better solution approximation via using second-order strain fields~\cite{le2020coarse}. These previous approaches partially overcome the inter-element stiffness caused by the usage of linear shape functions in conventional FE methods. However, the produced shape functions generally do not meet the basic property of partition of unity (PU), and may result in deformations of aphysical behaviours. To resolve the issue, a set of discontinuous and matrix-valued shape functions were derived by Chen~\cite{chen2018numerical}, where the the basic geometric properties of shape functions are imposed as constraints in an optimization problem. These studies however have not (fully) addressed the issue of the global solution continuity. Further discussion on the previous approaches is presented in Section~\ref{sec:other_methods}.

In this study, an approach is proposed for the analysis of heterogeneous structure of non-separated scales on a new concept of \emph{curved bridge nodes (CBNs)}, induced from a subset of the boundary nodes. In notable contrast to previous approaches solely working on corner nodes, the CBN analysis approach accommodates more DOFs in analysis by constructing a cubic B\'ezier curve along the interfaces between the coarse elements. Ultimately, it generates in an \emph{explicit} form a set of new CBN shape functions, which are applicable to both linear and nonlinear elasticity analysis problems.
The novel CBN shape functions further overcome the challenging issues of inter-element stiffness and ensures the global solution continuity in the fine mesh scale. It also meets the basic geometric properties of shape functions that avoid aphysical analysis results. Their analysis accuracy and efficiency are tested using various numerical examples, including a 3D industrial example of billions of DOFs, in comparison with classical approaches.

\begin{table}[htb]
  \centering
  \caption{Important notations in this paper}
  \label{tab:notations}
  \begin{tabular}{ccl}
  \hline
  $\Omega^{\alpha}$   & : & Coarse element                                                                  \\
  $\omega^{\alpha}_e$ & : & Fine element, simplified as $\omega_e$                                          \\
  $\calM^H$           & : & Coarse mesh, set of discrete coarse elements in the whole domain                             \\
  $\calM^{h}$         & : & Global fine mesh ,set of discrete fine elements in the whole domain                              \\
  $\calM^{\alpha,h}$  & : & Local fine mesh ,set of discrete fine elements in $\Omega^{\alpha}$                              \\
  $M$                 & : & Number of coarse elements of $\calM^H$                                          \\
  $m$                 & : & Number of fine elements of $\calM^{\alpha,h}$                                            \\
  $\calV_c$           & : & Corner nodes of local fine mesh $\calM^{\alpha,h}$                 \\
  $\calV_b$           & : & Boundary nodes of local fine mesh $\calM^{\alpha,h}$                                      \\
  $\calV_i$           & : & Interior nodes of local fine mesh $\calM^{\alpha,h}$                                      \\
  $\calV_r$           & : & Bridge nodes as subset of $\calV_b$ of local fine mesh $\calM^{\alpha,h}$               \\
  $\bQ$               & : & Vector of displacements of all CBNs of $\calM^H$                  \\
  $\bQ^{\alpha}$      & : & Vector of displacements of all CBNs of $\Omega^{\alpha}$                 \\
  $\bq_s$             & : & Vector of displacements of nodes in $\calV_s$, where $s$ could be $b,i$ \\
  $\bq$               & : & Vector of displacements of nodes in $\calM^{\alpha,h}$, including $\bq_i$ and $\bq_b$ \\
  $\bK^{\alpha}$      & : & Stiffness matrix of a coarse element $\Omega^{\alpha}$           \\
  $\bk^{\alpha}$      & : & Stiffness matrix of a local fine mesh $\calM^{\alpha,h}$                    \\
  $\bPsi$             & : & The B\'ezier interpolation matrix relating $\bQ^{\alpha}$ to $\bq_b$            \\
  $\bM^{\alpha}$      & : & The boundary-interior transformation matrix relating $\bq_b$ to $\bq_i$                    \\
  $\bN_e(\bx)$        & : & The basic bilinear shape function on point $\bx$ of a fine element $\omega_e$     \\
  $\bN^{h}(\bx)$      & : & Assembly of all $\bN_e(\bx)$ in local fine mesh $\calM^{\alpha,h}$                              \\
  $\bN^{\alpha}(\bx)$ & : & The CBN shape functions of a coarse element $\Omega^{\alpha}$     \\
  $\bI_{d,1}$         & : & An all-one vector with size of $d \times 1$ \\
  $\bI_d$             & : & An identity matrix with size of $d \times d$ \\ \hline
  \end{tabular}
\end{table}

\section{Problem statement and approach overview}
\begin{figure}[bt]
  \centering
  \subfigure[Structure and its coarse and fine elements]{\includegraphics[width=0.46\textwidth]{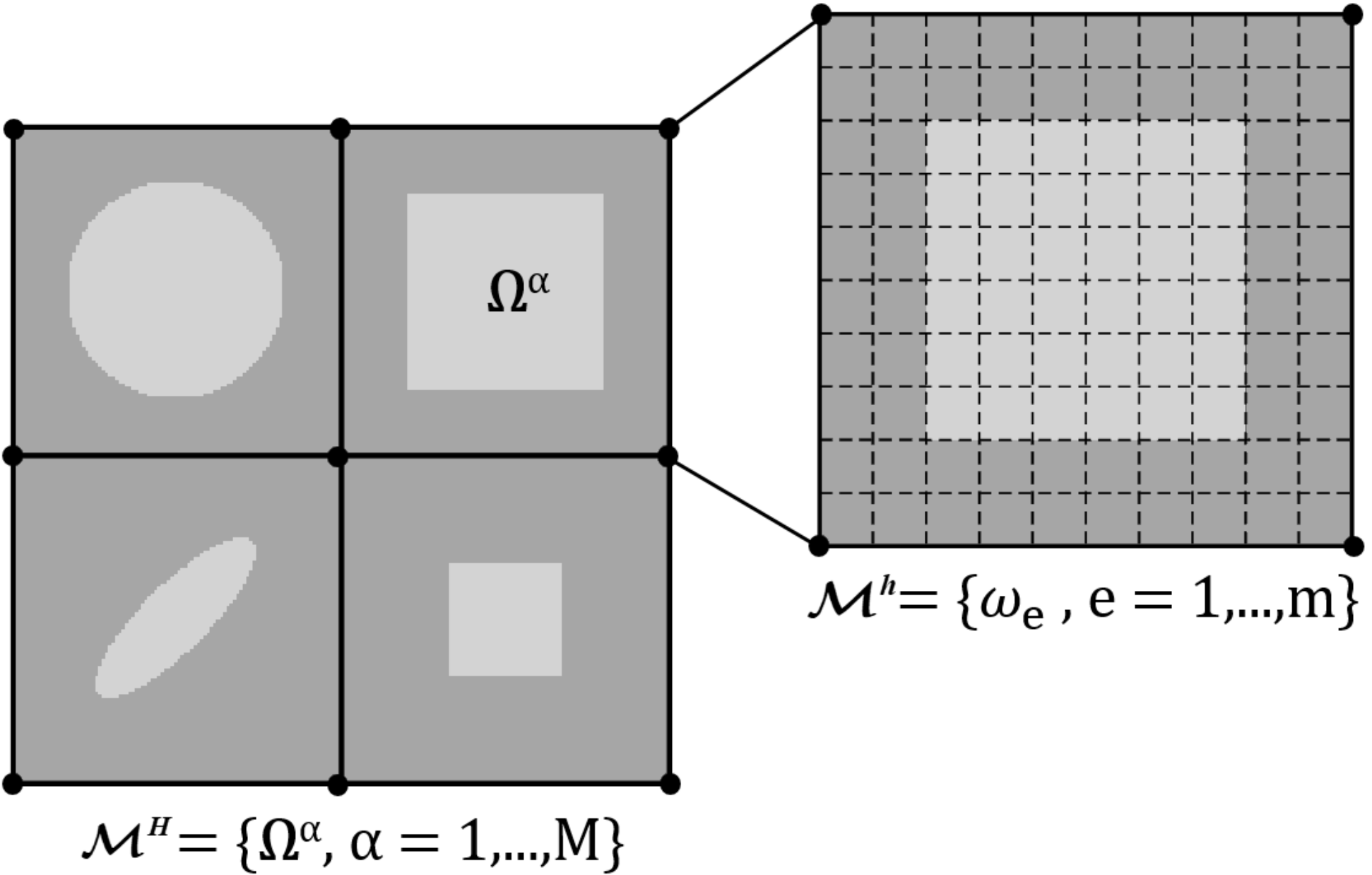}}\quad\quad\quad
  \subfigure[Various nodes and curved bridge nodes (CBNs)]{\includegraphics[width=0.24\textwidth]{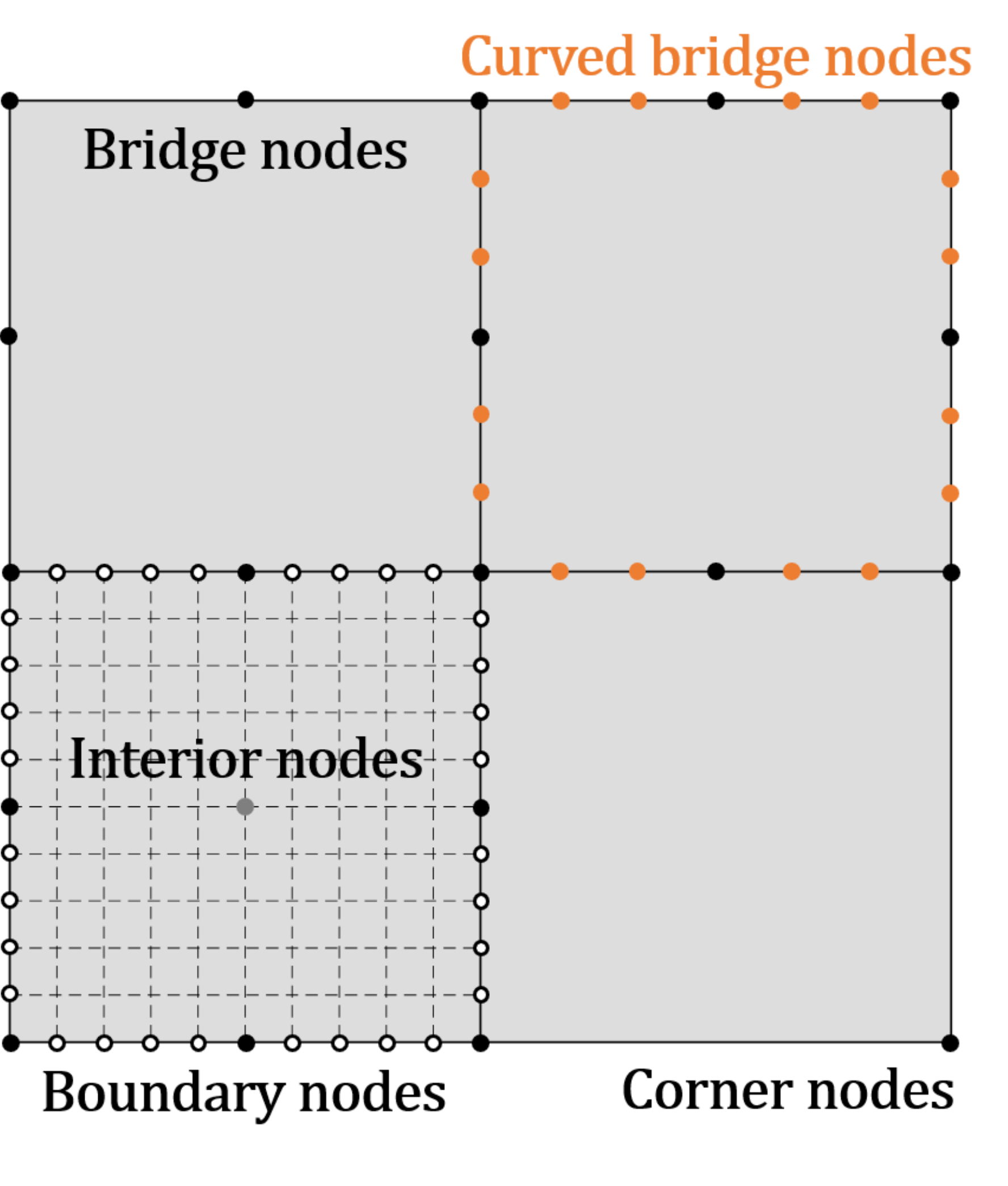}}
  \caption{(a). Structure $\Omega$, coarse mesh $\calM^H=\{\Omega^{\alpha}\}$, local fine mesh $\calM^{\alpha,h}=\{\omega^{\alpha}_e,\ e=1,2,\ldots, m\}$. (b). Bridge nodes are subset of boundary nodes; curved bridge nodes (CBNs) introduce additional nodes between adjacent bridge nodes, and are taken as analysis DOFs in our CBN heterogenous structure analysis approach. }
  \label{fig-domain}
\end{figure}

We mainly describe the approach for analysis of 2D linear elastic body. Its extensions to 3D case and to nonlinear models are explained later in Section~\ref{sec:extension}.

\subsection{Linear elasticity analysis of heterogeneous structures}\label{sec-elas}
As illustrated in Fig.~\ref{fig-domain}(a), let $\Omega \subset \RR^d$ for dimension $d=2,3$ be a heterogeneous solid structure under study, which may have different elasticity tensors $\bD(\bx)$ at different locations $\bx\in \Omega$. The linear elasticity analysis of $\Omega$ is described by a displacement vector for each point $\bx$ as $\bu(\bx)=(u(\bx),v(\bx))^T$. The \emph{strain vector} $\bveps(\bx,\bu)$ is defined as a linear approximation to the Green's strain and is represented in vector form as,
\eb\label{eq-strain}
\bveps(\bx,\bu)=(\bveps_{11} ,\bveps_{22},\sqrt{2}\bveps_{12})^T,
\ee
and the \emph{stress vector} $\bsig(\bu)$ is defined via Hooke's law,
\eb\label{eq-stress}
\bsig(\bx)=(\bsig_{11},\bsig_{22},\sqrt{2}\bsig_{12})^T=\bD(\bx)~\bveps(\bu).
\ee

The linear elasticity analysis of $\Omega$ aims to find the displacement $\bu$ satisfying
\eb\label{eq-elasticity}
\left\{
\begin{array}{lll}
&-\mbox{div} \bsig(\bu(\bx)) =\bg, & \mbox{in}\ \Omega,  \\
& \bsig(\bu(\bx))\cdot\n=\btau, & \mbox{on}\ \Gamma_N, \\
& \bu(\bx)=\bu_0, &\mbox{on}\ \Gamma_D,
\end{array}
\right.
\ee
where $\Gamma_D$ is a fixed boundary of a prescribed displacement $\bu_0$, $\Gamma_N$ the loading boundary of an external loading $\btau$, and $\bg$ is the body force.

The differential form in Eq.~\eqref{eq-elasticity} can also be stated in a weak form to induce its FE analysis: find the displacement $\bu\in H^1(\Omega)$ satisfying
\begin{equation}\label{eq-weak}
a(\u,\bv)=l(\bv),\quad \forall~ \bv\in H_0^1(\Omega),
\end{equation}
where
\eb
a(\u,\bv)=\int_{\Omega} \bsig(\bu)\cdot \bveps(\bv)~\rdV=\int_{\Omega} \bveps(\bu)^T~\bD(\bx)~\bveps(\bv)~\rdV,
\ee
and
\eb
l(\bv)=\int_{\Omega} \bg\cdot \bv~\rdV+\int_{\Gamma_N} \btau\cdot \bv~dS,
\ee
where $H^{1}(\Omega)$ and $H_{0}^{1}(\Omega)$ are the usual Sobolev vector spaces.

\subsection{Preliminary: bridge nodes and shape functions}
The FE analysis of the linear elasticity problem in Eq.~\eqref{eq-elasticity} or~\eqref{eq-weak} is usually conducted on a discretized domain of $\Omega$. Two different meshes are involved in this study. The \emph{coarse mesh} $\calM^H=\{\Omega^{\alpha},\ \alpha=1,2,\ldots, M\}$ contains a set of disjoint discrete \emph{heterogenous} \emph{coarse elements} $\Omega^{\alpha}$ of large size. Each coarse element $\Omega^{\alpha}$ further consists of a \emph{local fine mesh} $\calM^{\alpha,h}= \{\omega^{\alpha}_e,\ e=1,2,\ldots, m\}$ made of disjoint \emph{homogeneous} \emph{fine elements} $\omega^{\alpha}_e$ of much smaller size, which together induces the \emph{global fine mesh} $\calM^{h}= \{\omega^{\alpha}_e,\ e=1,2,\ldots, m,\ \alpha=1,2,\ldots, M\}$.

Each FE is formed by a set of nodes. Given a local fine mesh $\calM^{\alpha,h}$ for a coarse element $\Omega^{\alpha}$, we classify their nodes into three different subsets depending on their locations: sets of \emph{corner nodes, boundary nodes or interior nodes}, respectively denoted as $\calV_c,\ \calV_b,\ \calV_i$ if they are on the corners, boundaries, or interiors of $\Omega^{\alpha}$. See also Fig.~\ref{fig-domain}(b).

We also introduce the concept of \emph{bridge node} set as a subset of the boundary node set, denoted $\calV_r$ and defined below
\[
\calV_c \subseteq \calV_r \subseteq \calV_b.
\]

Given a segment determined by a pair of adjacent bridge nodes, a set of equally spaced nodes are inserted, which together with those in $\calV_r$ form the set of \emph{curved bridge nodes (CBNs)}. These nodes are used in the downstream task in high-order interpolation curve construction. See also Fig.~\ref{fig-domain}(b).

Shape functions play a role of basis functions in FE analysis, whose linear combination describes a deformation of structure under study. We first explain the definition on the local fine mesh $\calM^{\alpha,h}$. Given a master FE $\omega_e$ with four corner nodes within the range $(x,y)\in [-1,1]\times[-1,1]$ and numbered from $1$ to $4$, the scalar \emph{bilinear shape function} $N_i(\bx)$ is defined below:
\eb\label{eq-linear-sf}
N_i(\bx):\Omega\rightarrow \RR,\quad  N_i(\bx)=\frac{1}{4}(1+x_ix)(1+y_iy) \quad, i=1,2,3,4.
\ee

Accordingly, the solution $\bu(\bx)$ to problem~\eqref{eq-elasticity}, \eqref{eq-weak} is interpreted as a linear combination of the shape functions, or in a matrix form,
\eb\label{eq:feelmentint}
\bu_e(\bx) ={\bN_e(\bx)}~{\bq_e},\quad \bx\in\omega_e,
\ee
where $\bq_e$ is the displacement vector of dimension $8\times 1$, and $\bN_e(\bx)$ is the matrix form of $N_i(\bx)$ of dimension $2\times 8$,
\eb\label{eq:sfm}
\bN_e(\bx)=
\left[
\begin{array}{cccccccc}
N_1(\bx) & 0 & N_2(\bx) & 0 & N_3(\bx) & 0 & N_4(\bx) & 0 \\
0 & N_1(\bx) & 0 & N_2(\bx) & 0 & N_3(\bx) & 0 & N_4(\bx)
\end{array}
\right].
\ee

Consequently, the \emph{partition of unity (PU)} and \emph{Lagrange property} are satisfied for $N_i(\bx),\ i=1,2,3,4$, that is,
\eb
\bN_e(\bx) ~ \bI_{8,1} = \bI_{2,1}, \quad\mbox{and}\quad N_i(\bx_j)=\delta_{ij},
\ee
where $\bI_{d,1}$ is an all-one vector with size of $d\times 1$, and $\delta_{ij}$ is the Kronecker delta function.

The property of \emph{node interpolation} also comes directly from the Lagrange property,
\eb
{\bN_e(\bx_j)}~{\bq_e}=\bq_{e,j},
\ee
where $\bq_{e,j}$ is the displacement of $\bq_e$ on a node $\bx_j$.

The above bilinear shape functions can be defined on the global fine mesh $\calM^h$ or the coarse mesh $\calM^H$. On $\calM^h$, the shape functions are defined on a homogeneous fine element $\omega_e^{\alpha}$, and effectively approximate the target solution. However, its numerous fine elements result in a problem that is too computationally expensive. In contrast, on coarse mesh $\calM^H$, the shape functions are defined on a heterogeneous coarse element $\Omega^{\alpha}$, potentially losing a high amount of  solution accuracy. We resolved the issues by constructing a set of material-aware shape functions for the CBNs, as discussed in the following section.

\subsection{Approach overview}
\begin{figure}[tb]
  \centering
  \subfigure[Bilinear shape functions on $\calM^h$]{\includegraphics[width=0.18\textwidth]{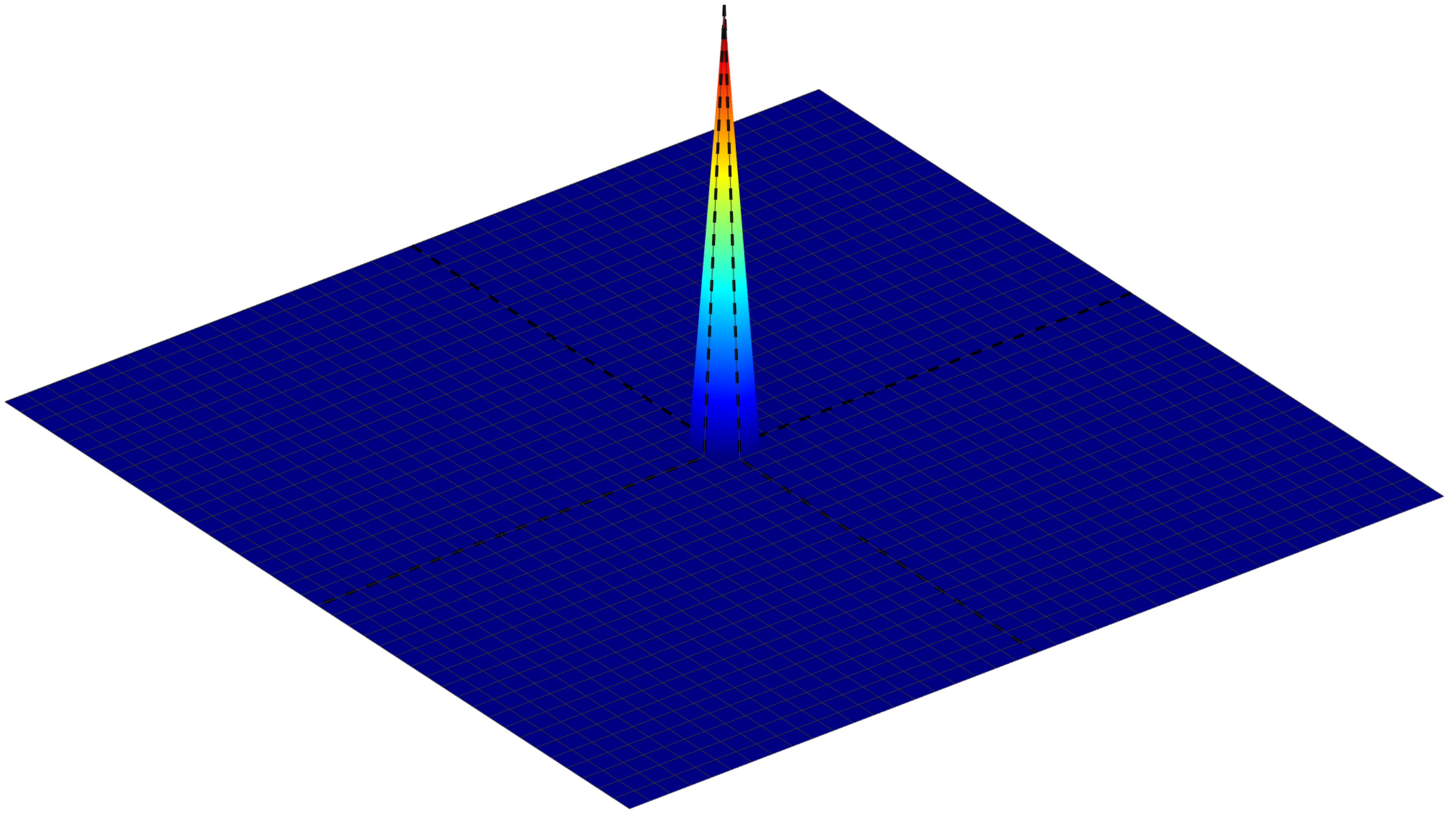}}\quad
  \subfigure[Bilinear shape functions on $\calM^H$]{\includegraphics[width=0.18\textwidth]{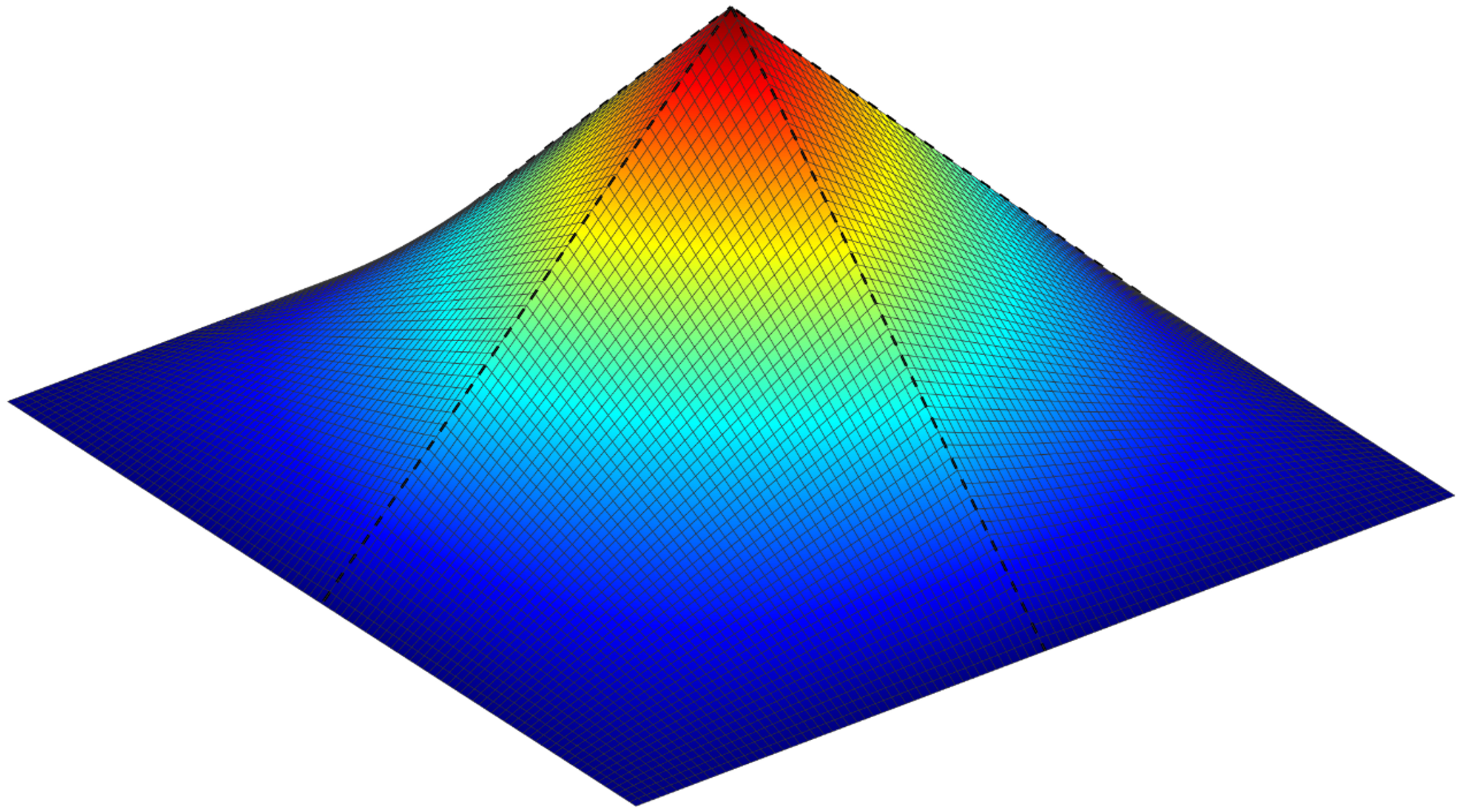}}\quad
  \subfigure[Shape functions without continuity consideration]{\includegraphics[width=0.18\textwidth]{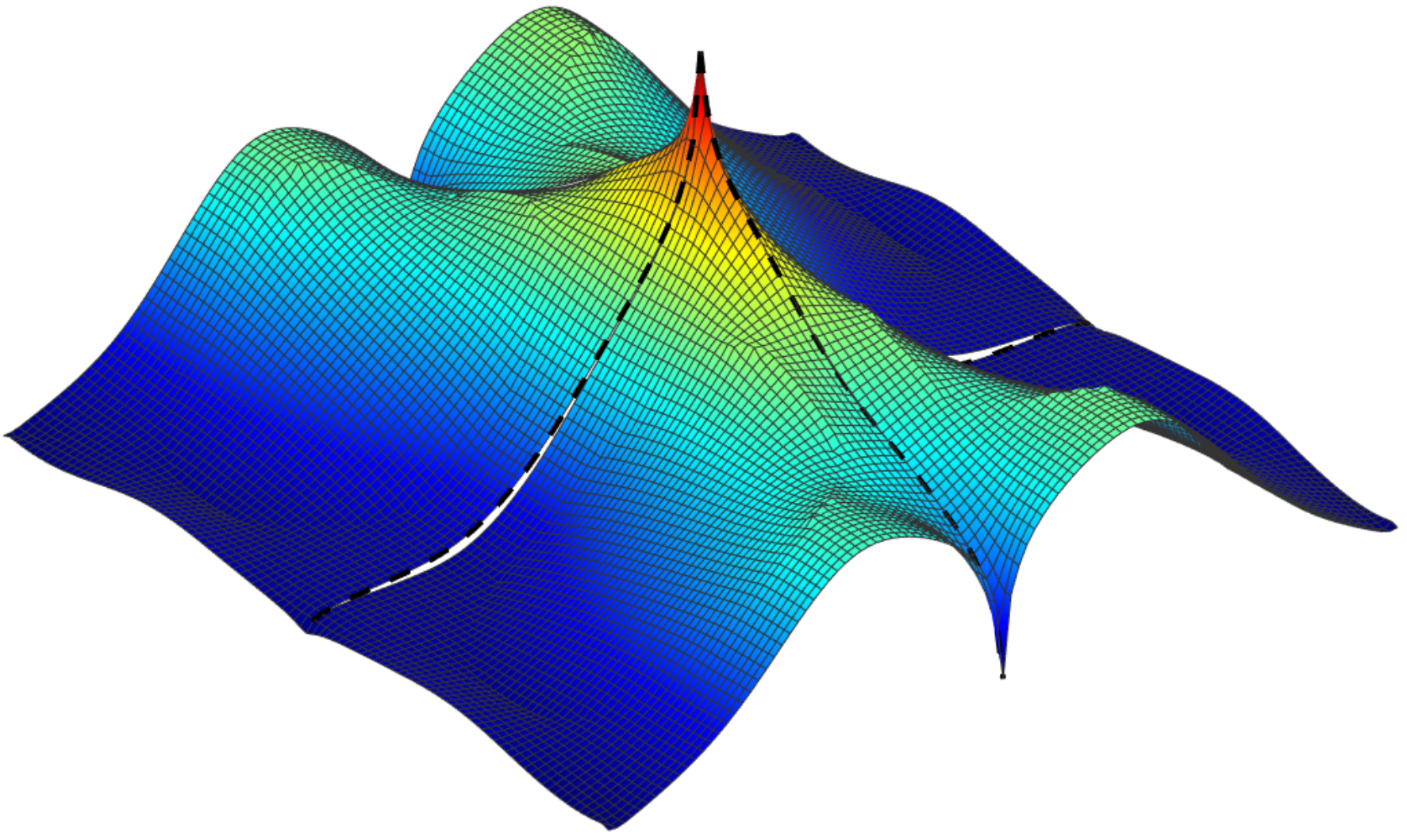}}\quad
  \subfigure[Shape functions with linear interpolation on interface]{\includegraphics[width=0.18\textwidth]{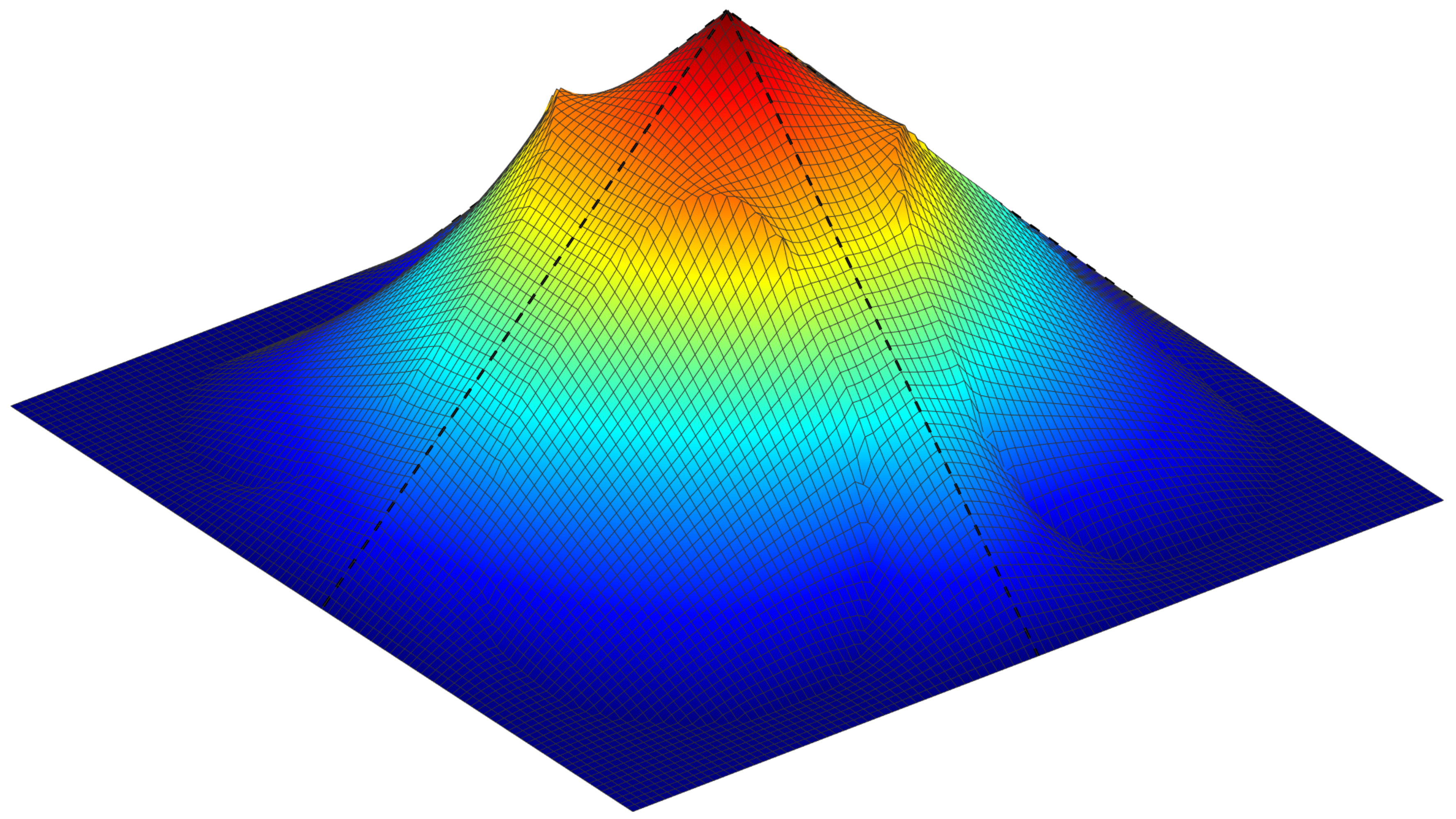}}\quad
  \subfigure[Shape functions on CBNs with curved interpolation]{\includegraphics[width=0.18\textwidth]{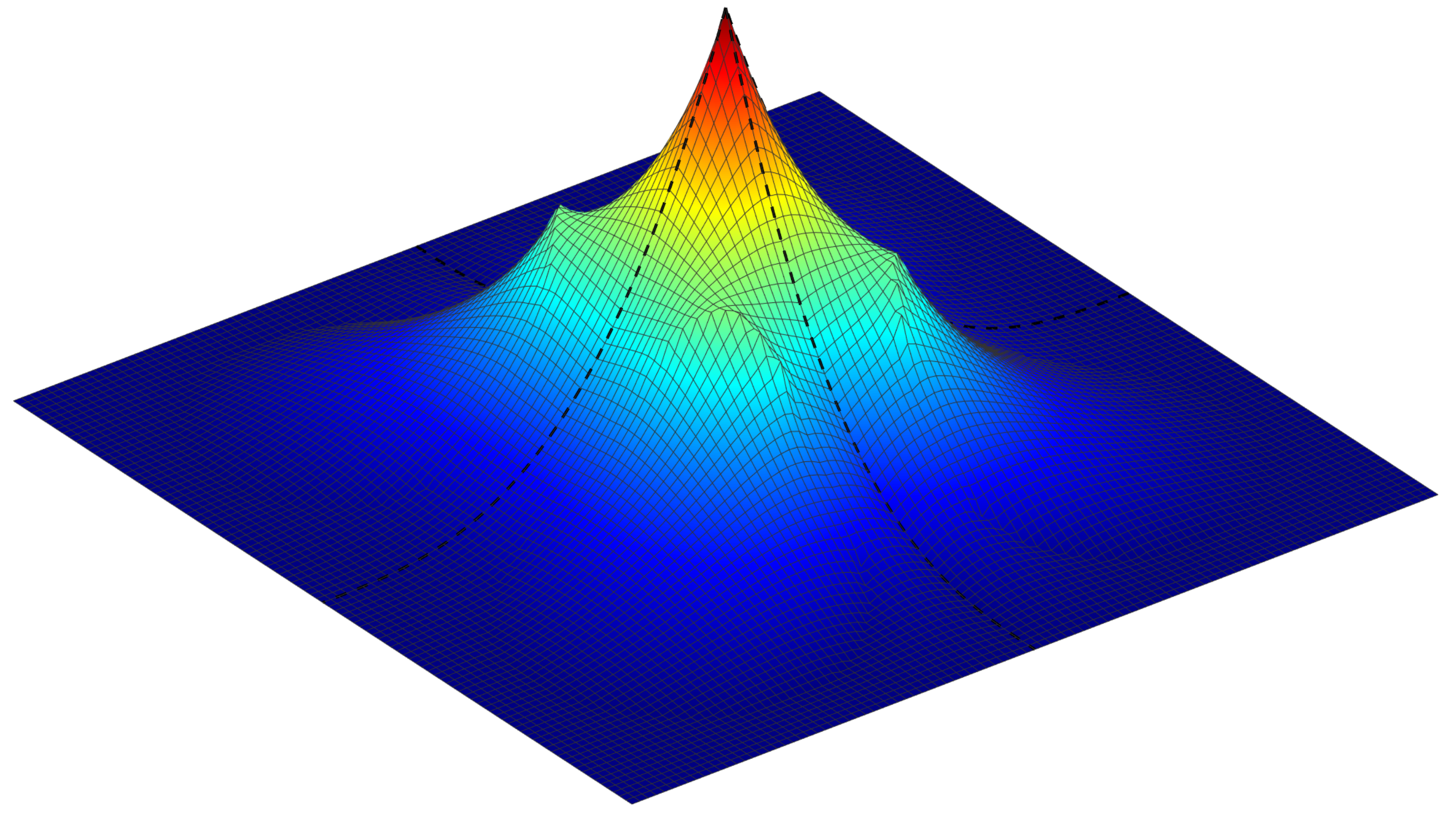}}
  \caption{Shape functions under different construction strategies.}
  \label{fig:surf_diff_NH}
  \end{figure}
  \begin{figure}
  \centering
  \subfigure[On $\calM^h$]{\includegraphics[width=0.16\textwidth]{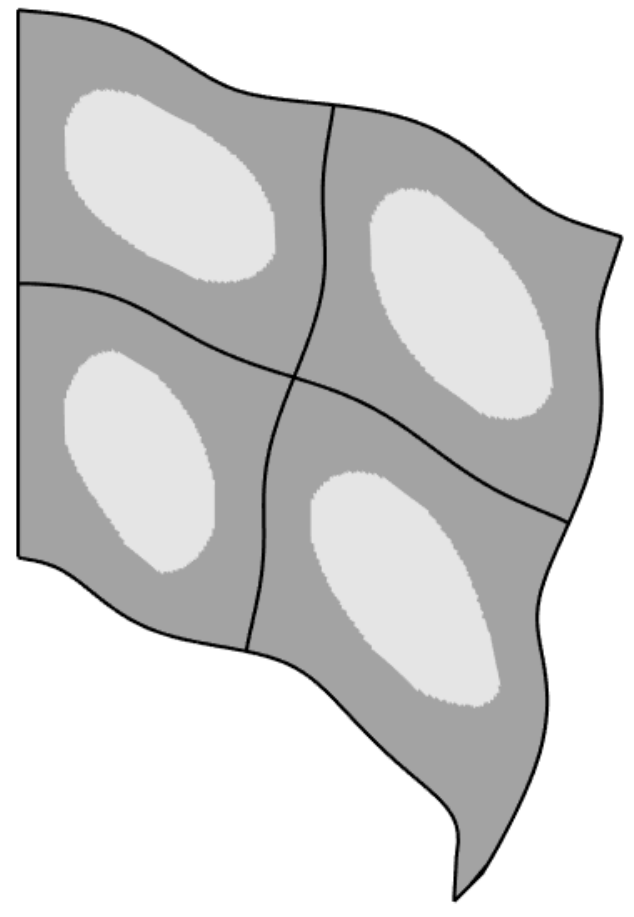}} \quad
  \subfigure[On $\calM^H$]{\includegraphics[width=0.16\textwidth]{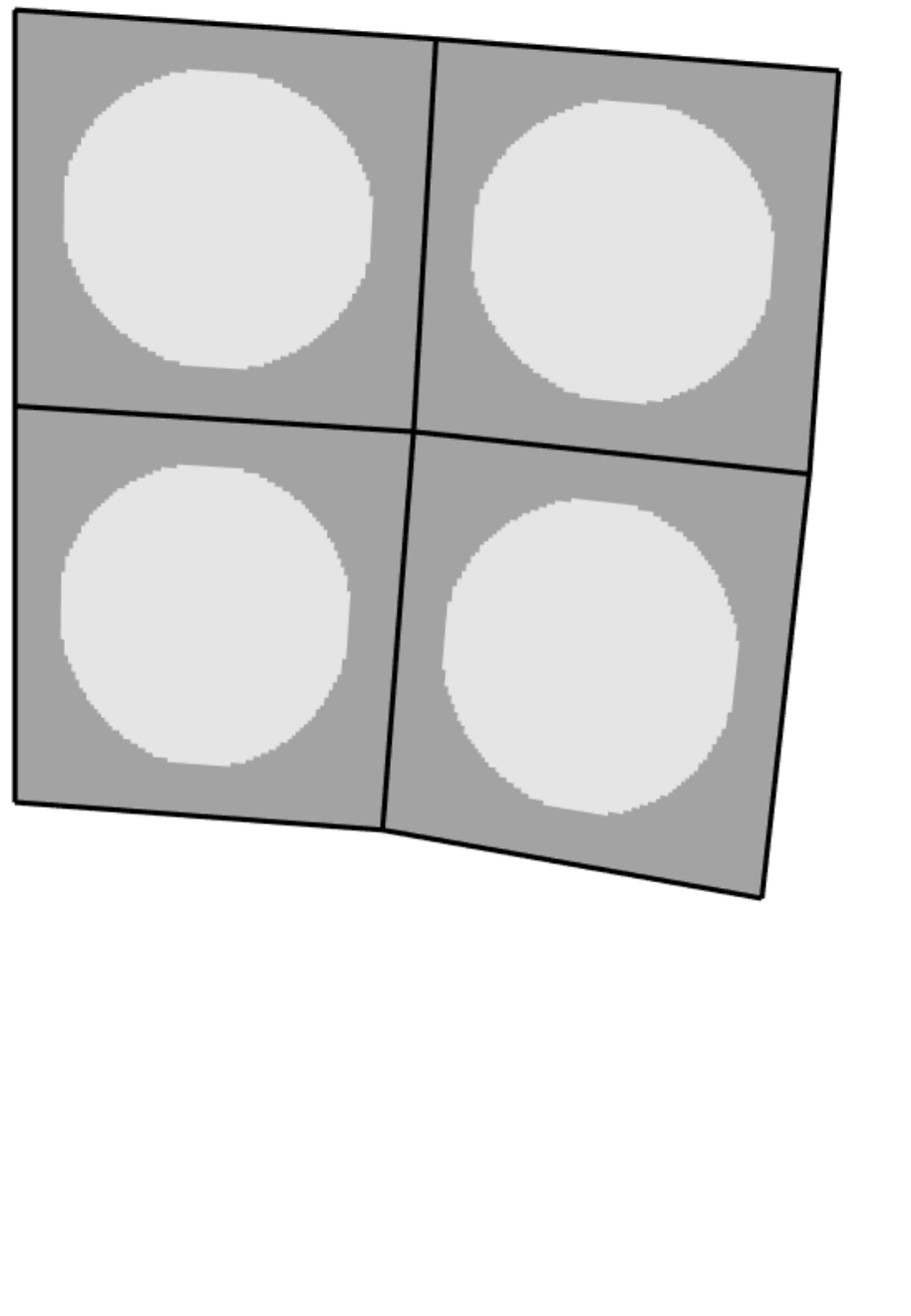}} \quad
  \subfigure[Without continuity consideration]{\includegraphics[width=0.16\textwidth]{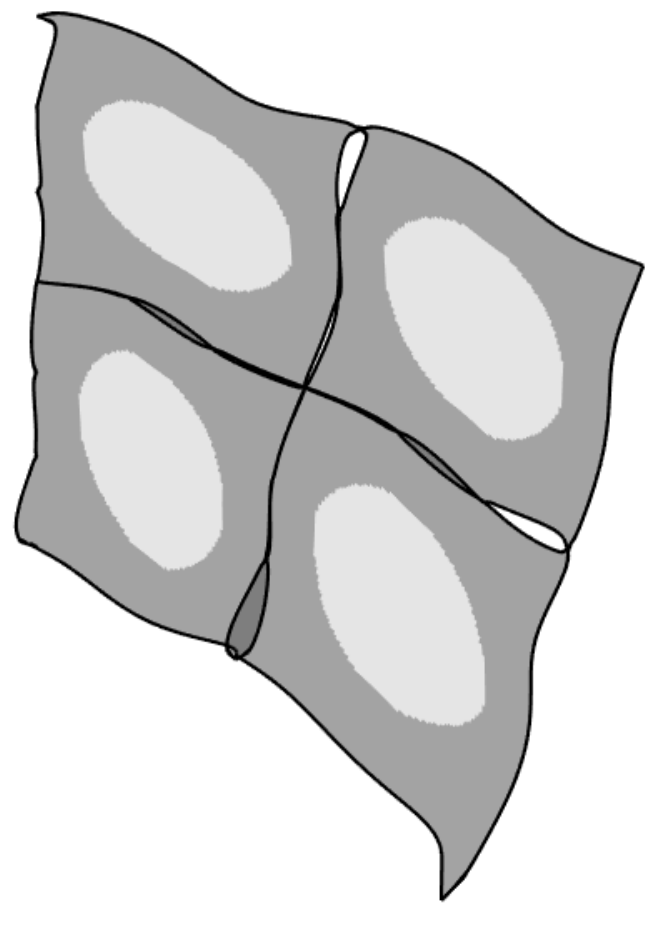}} \quad
  \subfigure[With linear interpolation on interface]{\includegraphics[width=0.16\textwidth]{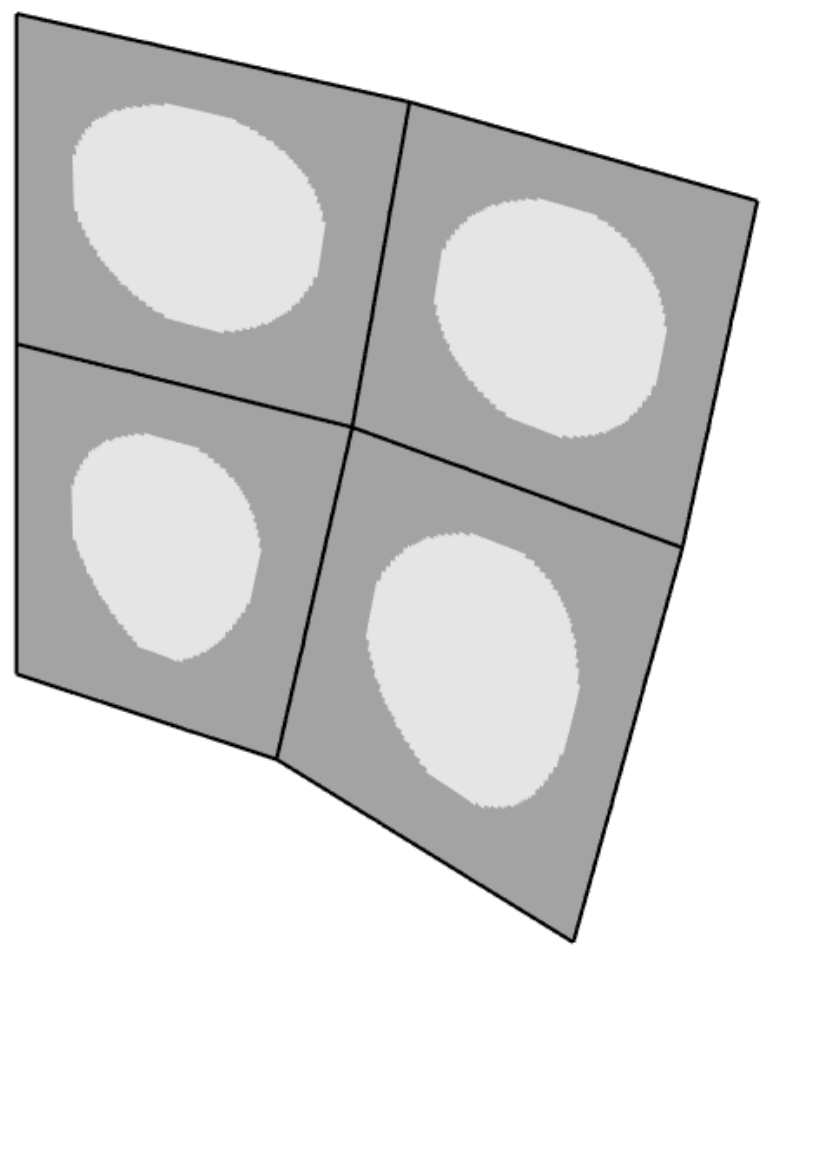}} \quad
  \subfigure[On CBNs with curved interpolation]{\includegraphics[width=0.16\textwidth]{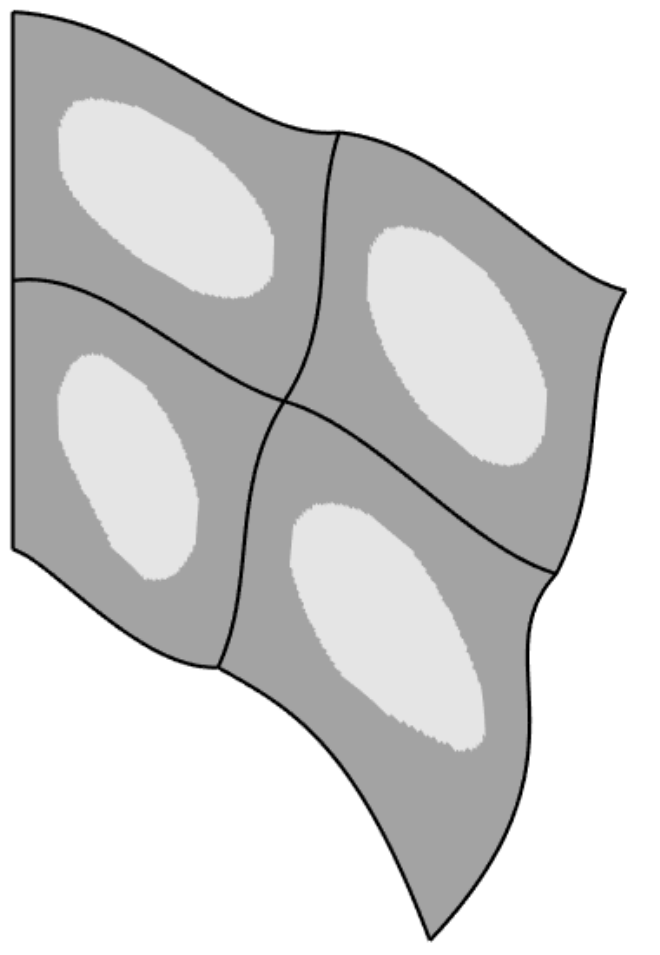}}
  \caption{Analysis results using the different shape functions defined in Fig.~\ref{fig:surf_diff_NH}.}
  \label{fig:continuity}
\end{figure}

Following classical Galerkin FE method, shape functions play a role of bases to produce the overall displacement with respect to a vector of discrete nodal values. Instead of choosing the corner nodes, the CBNs are introduced here and set as the coarse nodes for more analysis DOFs and more flexibilities of shape descriptions.

In this study, a set of material-aware CBN shape functions $\bN^{\alpha}(\bx)$ is to be constructed for each coarse element $\Omega^{\alpha}$. Let $\bQ$ be the vector of discrete displacements on CBNs in $\calM^H$ to be determined, and $\bQ^{\alpha}$ be its component on $\Omega^{\alpha}$. Accordingly, the displacement to Eq.~\eqref{eq-elasticity} takes the following form,
\eb
\bu(\x)=\sum_{\alpha=1}^M\bN^{\alpha}(\bx)~\bQ^{\alpha}.
\ee

In order to improve the ability for describing the heterogeneity of fine mesh, we construct the CBN shape function $\bN^{\alpha}(\bx)$ as a piecewise-bilinear function (in 2D) defined over the local fine mesh $\calM^{\alpha,h}(\bx)$,
\eb\label{eq-gshapef}
\bN^{\alpha}(\bx)= \bN^h(\bx)~\tilde{\bPhi},
\ee
where $\tilde{\bPhi}$ is a matrix of DOFs to be determined to closely capture the coarse element's heterogeneity, and $\bN^{h}(\bx)$ is a matrix of the fine-mesh shape functions $\bN_e(\bx)$ in $\calM^{\alpha,h}$,
\eb\label{eq:Nh}
\bN^{h}(\bx) = {\sum_{e=1}^m} \bN_e(\bx), \quad \bx \in \Omega^{\alpha},
\ee
where $\sum$ denotes the \emph{assembly sum} in numerical FE assembly process that conducts the summation on the same location; given a specific point $\bx_0$, the value of $\bN^{h}(\bx_0)$ can then be directly evaluated.


In fact, the construction of effective shape functions is confronted by at least two known challenges: the inter-element stiffness and the displacement discontinuity across the coarse element interface. The inter-element stiffness issue originates from the usage of linear interpolation in the reconstruction of the global fine-mesh displacement from the discrete coarse nodal displacement~\cite{nesme2009preserving}. The interface discontinuity issue is mainly due to the fact that the shape functions are usually locally constructed without considering the adjacency of the coarse elements, and thus they may have different values along the common interface~\cite{chen2018numerical}.
In practice, different values of $\tilde{\bPhi}$ determine different shape functions, as indicated in Fig.~\ref{fig:surf_diff_NH}(c),(d),(e). These consequently result in very different analysis results as shown in Fig.~\ref{fig:continuity}(c),(d),(e); their counterparts from the fine mesh $\calM^h$ or the coarse mesh $\calM^H$ are respectively shown in Figs.~\ref{fig:surf_diff_NH}(a),(b) and Fig.~\ref{fig:continuity}(a),(b). Specifically, the inter-element stiffness is observed in Fig.~\ref{fig:continuity}(d) owing to simple linear interpolation, and the construction without continuity consideration leads to the interface overlap and discontinuity in Fig.~\ref{fig:continuity}(c).

In an effort to further address the above-mentioned challenges, this study aims to develop a new class of material-aware shape functions, known as CBN shape functions. Given a master coarse element, this is achieved by treating the shape functions construction as a process to build up a map from the coarse DOFs (or displacements on CBNs) to the local fine displacements per coarse element, instead of formulating it as a constrained nonlinear optimization problem.

Firstly, it constructs cubic B\'ezier interpolation curves from the CBNs along the coarse element's boundaries, which not only ensures the continuity of the global displacement in fine mesh but also improves its accuracy by allowing for more deformation flexibilities. Secondly, it maps the boundary nodal displacements to those on the interior nodes, which builds on its intrinsic physical properties and closely capture the coarse element's heterogeneities. Ultimately, the shape functions are derived in an explicit matrix form as a product of two matrix transformations: the B\'ezier interpolation transformation and the boundary-interior transformation, and preserve the basic geometric properties of shape functions.

In summary, the derived shape functions present the following properties.
\begin{enumerate}
\item Defined with respect to CBNs with flexible analysis DOFs, further overcoming inter-element stiffness while maintaining the global displacement smoothness.

\item Expressed in an explicit matrix form as a product of two transformations: B\'ezier interpolation and the boundary-interior transformations.

\item Preserving the basic geometric properties of shape functions: node interpolation, translation and rotation invariants, to avoid aphysical analysis behaviours.

\item Applicability to both linear and nonlinear analysis problems.

\end{enumerate}

\section{Construction of CBN shape functions}
Considering Fig.~\ref{fig-domain}, let $\Omega^{\alpha}$ be a master coarse element, $\calM^{\alpha,h}$ its local fine mesh, and $\calX_b$ the given bridge node set. The CBN shape functions $\bN^{\alpha}(\bx)$ are achieved as a product of two matrix transformations: the B\'ezier interpolation transformation and the boundary-interior transformation, as detailed below.

\subsection{Construction of B\'ezier interpolation transformation}\label{sec:Bezier}
\begin{figure}[t]
  \centering
  \subfigure[Linear interpolation]{\includegraphics[width=0.33\textwidth]{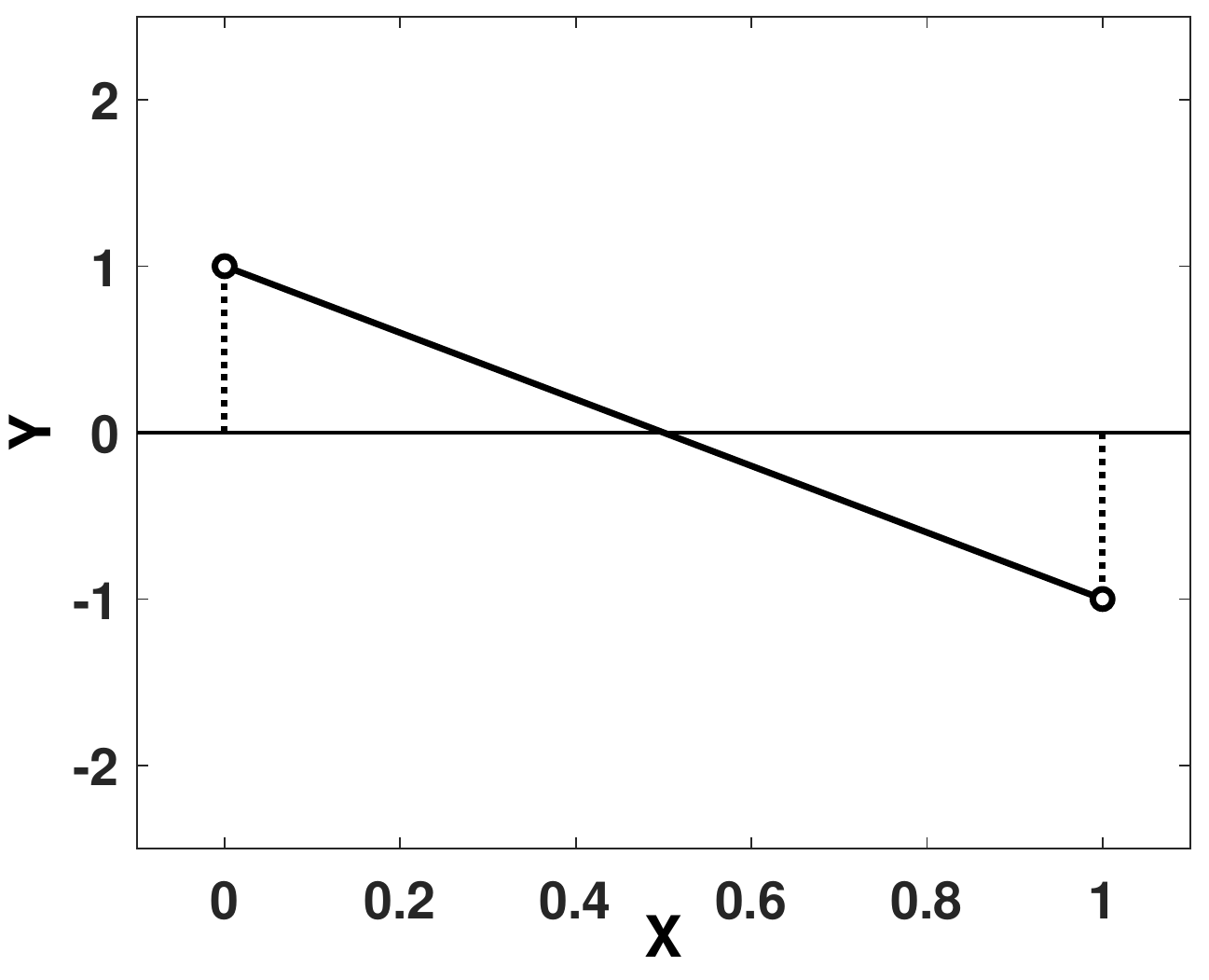}}
  \subfigure[Linear interpolation with additional nodes]{\includegraphics[width=0.33\textwidth]{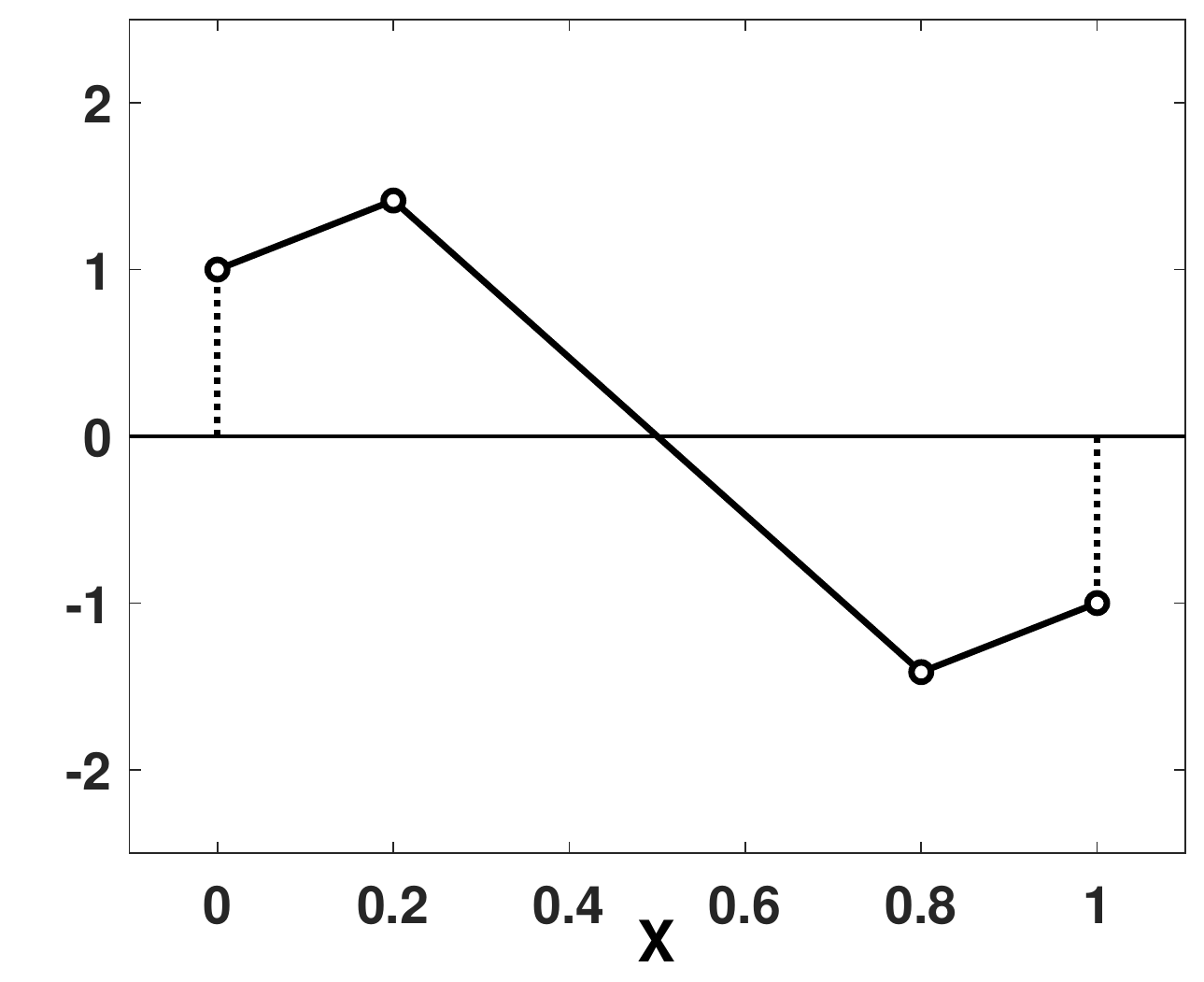}}
  \subfigure[B\'ezier interpolation with additional nodes]{\includegraphics[width=0.33\textwidth]{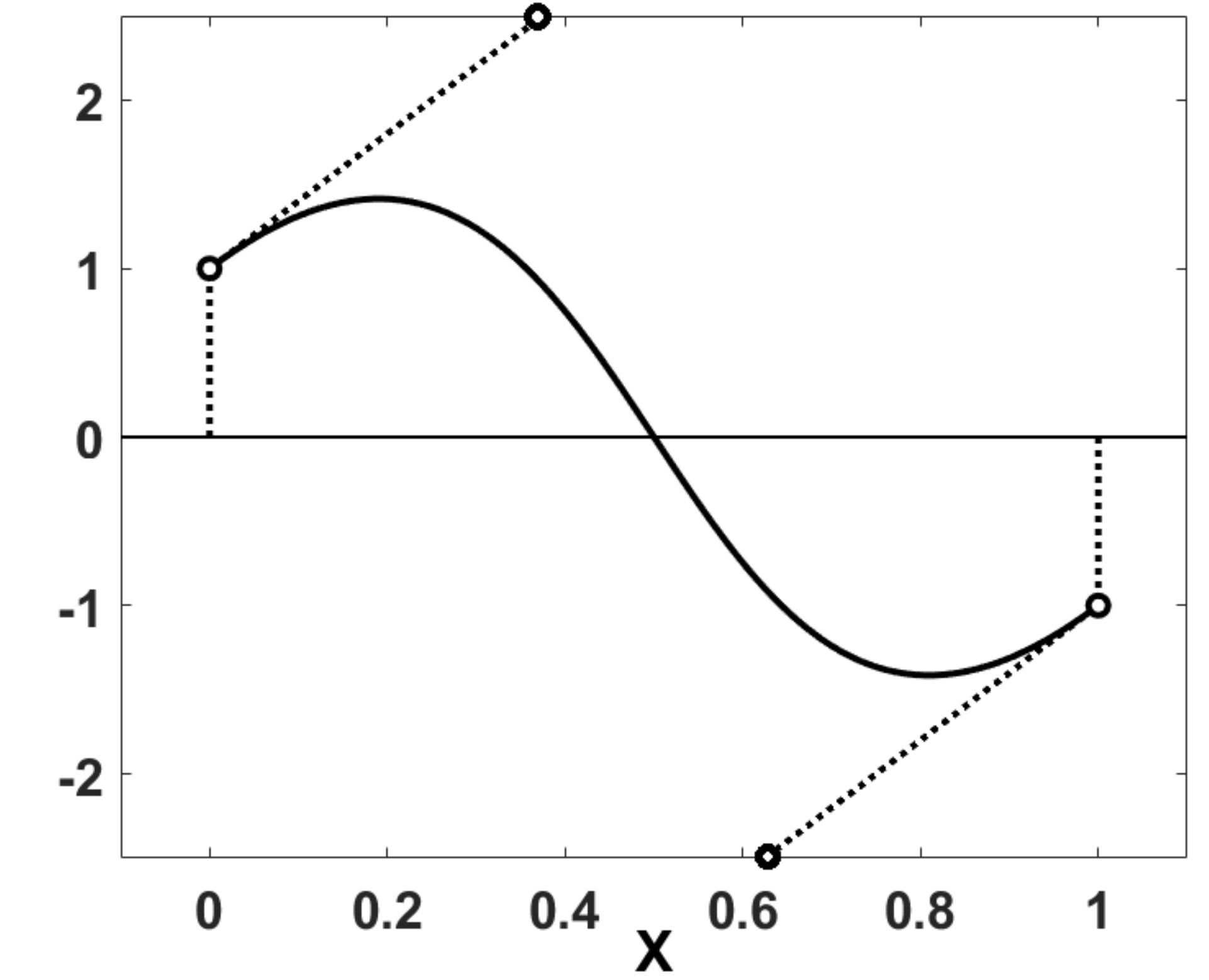}}
  \caption{Different interpolation strategies produce different boundary curves, where the higher order cubic B\'ezier curve has more flexibility.}
  \label{fig:strategy_interpolation}
\end{figure}

The B\'ezier interpolation transformation constructs a map from the discrete displacements on the CBNs $\bQ^{\alpha}$ to those on the bridge nodes in $\calX_b$. The higher-order interpolation of B\'ezier curve not only allows for more control flexibility but also has the favorable geometric properties of PU and translation/rotation invariants. Fig.~\ref{fig:strategy_interpolation} illustrates the difference between different interpolation strategies.

A cubic B\'ezier curve $\bP(t)$ is taken here in the following form,
\eb\label{eq-cbezier}
\bP(t)=\sum_{i=0}^3\psi_{i}(t)~\bP_i=\sum_{i=0}^3  C^i_nt^i(1-t)^{3-i}~\bP_i,
\ee
where $\psi_{i}(t)$ is a cubic Bernstein base, $\bP_i \in \RR^2$ are the control points, and $C_n^i$'s are the binomial coefficients.

Note $\bP_i$ are the values of $P(t)$ at node $t=i/3$ for $i=0,1,2,3$, or,
\eb\label{eq:bzsmp}
\bP(0)=\bP_0,\quad \bP(1/3)=\bP_1,\quad \bP(2/3)=\bP_2,\quad \bP(1)=\bP_3.
\ee
In addition, $\bP(t)$ has the following nice properties of the translation invariant and rotation invariant for a constant angular velocity $\hat{\theta}$,
\eb
\sum_{i=0}^3\psi_i(t)=1, \quad \hat{\theta}\times \bP(t)=\sum_{i=0}^3\psi_{i}(t)~\hat{\theta}\times \bP_i,
\ee
which are used in delivering such properties of our CBN shape functions.

%
%

Now we consider the approach to construct the B\'ezier transformation matrix.  Following the idea of FE displacement expression in~\eqref{eq:feelmentint}, we rewrite $\bP(t)$ in the following matrix form,
\eb\label{eq-bcmatrix}
\underset{(2\times 1)}{\bP(t)\vphantom{\bpsi}}=\underset{(2\times 8)}{\bpsi(t)}~\underset{(8\times 1)}{\bP\vphantom{\bpsi}},
\ee
for
\eb\label{eq-bp}
\bP=[\bP_0,\bP_1,\bP_2,\bP_3]^T,
\ee
and
\eb\label{eq:basis}
    \begin{aligned}
      \bpsi(t) = \bpsi_0(t)\otimes \bI_2=\begin{bmatrix}
        1 & t & t^{2} & t^{3}
      \end{bmatrix} \begin{bmatrix}
        1  & 0  & 0  & 0 \\
        -3 & 3  & 0  & 0 \\
        3  & -6 & 3  & 0 \\
        -1 & 3  & -3 & 1
      \end{bmatrix} \otimes \bI_2
    \end{aligned}
\ee
\eb
=\left[
\begin{array}{cccccccc}
\psi_0(t) & 0 & \psi_1(t) & 0 & \psi_2(t) & 0 & \psi_3(t) & 0  \\
0 & \psi_0(t) & 0 & \psi_1(t) & 0 & \psi_2(t) & 0 & \psi_3(t)
\end{array}
\right],
\ee
where the Kronecker product with identity matrix $\bI_2$ is to match the dimensions of the column vector $\bP$.

Suppose $E$ is a \emph{bridge segment} bounded by a pair of adjacent bridge nodes in $\calX_b$, and $t(\cdot)$ is a function that re-parameterizes $E$ into a parametric curve in a range of $[0,1]$. By inserting two additional equally spaced nodes along $E$, or at $t=1/3,2/3$, in total, we have the associated four CBNs (see Fig.~\ref{fig-domain}(b)).

Let $\bQ_E$ be the vector of the x, y displacements at the four CBNs. Taking $\bQ_E$ as control points in the cubic B\'ezeir curve function~\eqref{eq-bcmatrix}, we have the interpolation displacement function $\bu_E(\bx)$ along segment $E$,
\eb\label{eq-bezier-x}
\underset{(2\times 1)}{\bu_E(\bx)} = \underset{(2\times 8)}{\bpsi_E(t(\bx))}~\underset{(8\times 1)}{\bQ_E}.
\ee

From Eq.~\eqref{eq:bzsmp}, it can be noted that evaluating $\bu_E(\bx)$ at the four CBNs in $E$ gives the four control points $\bQ_E$. The relation is now to be derived on the full boundary nodes in $\calX_b$.

Consider a specific boundary node $\bx_0$ in $\calV_b$ located at a bridge segment $E$. Evaluating $\bpsi_E(t(\bx))$ at $\bx_0$ gives its interpolated displacement $\bpsi_E(t(\bx_0))$. Following a similar FE assembly process, we have
\eb
\underset{(2\times 6r)}{\bpsi(\bx_0)}= \sum_E \bpsi_E(t(\bx_0)) , \quad \bx_0 \in E.
\ee

Iterating $\bx_0$ for all the bridge nodes in $\calX_b$, we consequently have
\eb\label{eq-qb}
\underset{(2b\times 1)}{\bq_b}=\underset{(2b\times 6r)}{\bPsi\vphantom{\bQ^{\alpha}}}~\underset{(6r\times 1)}{\bQ^{\alpha}},
\ee
where the \emph{B\'ezier interpolation matrix} is
\eb
\label{eq:bezier-psi}
\underset{(2b\times 6r)}{\bPsi}=\left[\bpsi(\bx_0), \ \bx_0\in\calV_b\right],
\ee
and we list their corresponding base $\bpsi(\bx_0)$ row by row, and $\bQ^{\alpha}$ is the vector of all displacements on CBNs of $\Omega^{\alpha}$.

The dimension $6r$ of $\bQ^{\alpha}$ can be seen from the fact: we have $r$ bridge segments from $t$ bridge nodes, which together have $3r$ CBNs, and thus of a dimension $6r$ considering its x-, y- components.

\subsection{Construction of boundary-interior transformation} \label{sec:NH-construction}
Given the vector $\bq_b$ of the boundary displacements in Eq.~\eqref{eq-qb}, we construct a boundary-interior transformation to map it to the interior displacements. Let $\bk^{\alpha}$ be the associated stiffness matrix on the local fine mesh $\calM^{\alpha,h}$. Reordering all DOFs to partition them into internal and boundary entries indexed by $i$ and $b$, denoted as vectors $\bq_i, \bq_b$. Their relation is determined by the following FE equilibrium equation,
\begin{align}\label{eq:micros}
  \begin{bmatrix}
    \bk_{b} & \bk_{bi} \\
    \bk_{ib} & \bk_{i}
  \end{bmatrix}\begin{bmatrix}
    \bq_b \\
    \bq_i
  \end{bmatrix} = \begin{bmatrix}
    \bf_b \\
    0
  \end{bmatrix},
\end{align}
where $\bk_{b},\ \bk_{i},\ \bk_{bi},\ \bk_{ib}$ are the associated stiffness sub-matrices of $\bk^{\alpha}$, and $\bf_b$ is vector of exposed forces on the boundary nodes.

We have from the second-row
\eb
  \label{eq:micro-interior}
  \underset{(2i\times 1)}{\bq_i} = \underset{(2i\times 2b)}{\bM^{\alpha}\vphantom{\bq_x}}~\underset{(2b\times 1)}{\bq_b},\quad \mbox{for}\quad   \underset{(2i\times 2b)}{\bM^{\alpha}\vphantom{\bq_x^{-1}}} = \underset{(2i\times 2i)}{-\bk_{i}^{-1}}~\underset{(2i\times 2b)}{\bk_{ib}\vphantom{\bq_x^{-1}}}.
\ee
Accordingly, we have the vector of the displacements on $\calM^{\alpha,h}$,
\eb\label{eq:bdintmatrix}
\underset{(2i+2b)\times 1}{\bq}= [\bq_b,\bq_i]^T= \underset{(2i+2b)\times 2b} {\tilde{\bM}^{\alpha}\vphantom{\bq}}~\bq_b,
\ee
where $\tilde{\bM}^{\alpha}$ is the desired material-aware \emph{boundary-interior transformation matrix}
\eb\label{eq:bitrans}
\underset{(2i+2b)\times 2b}{\tilde{\bM}^{\alpha}\vphantom{\bI_{2b}}}=[\underset{2b\times 2b}{\bI_{2b}},\underset{2i\times 2b}{\bM^{\alpha}\vphantom{\bI_{2b}}}]^T.
\ee

\subsection{Shape functions in a matrix form} \label{sec:NH}
Substituting Eqs.~\eqref{eq-qb} into Eq.~\eqref{eq:bdintmatrix} further gives
\eb \label{eq:fine-q}
\underset{(2i+2b)\times 1}{\bq}= \underset{(2i+2b)\times 2b} {\tilde{\bM}^{\alpha}\vphantom{\bq}}~\underset{(2b\times 6r)}{\bPsi
\vphantom{\bq}}~\underset{(6r\times 1)}{\bQ^{\alpha}}.
\ee
In combination with the shape functions $\bN^h(\bx)$ to the local fine mesh $\calM^{\alpha,h}$, as defined in Eq.~\eqref{eq:Nh},  the displacement on any point $\bx\in\Omega^{\alpha}$ is interpolated
\eb\label{eq-3maps}
\underset{2\times 1}{\bu^{\alpha}(\bx)\vphantom{_x}} =
\underset{2\times (2i+2b)} {\bN^{h}(\bx)}~\underset{(2i+2b)\times 1}{\bq}=
\underset{2\times (2i+2b)} {\bN^{h}(\bx)} ~ \underset{(2i+2b)\times 2b} {\tilde{\bM}^{\alpha}\vphantom{_x}} ~ \underset{(2b\times 6r)}{\bPsi\vphantom{_x}}~ \underset{(6r\times 1)}{\bQ^{\alpha}},\quad \bx\in\Omega^{\alpha}.
\ee
The above equation maps the discrete nodal displacement $\bQ^{\alpha}$ to the continuous interpolated displacement $\bu^{\alpha}(\bx)$. It accordingly gives our CBN shape functions in a matrix form,
\begin{align}
  \label{eq:NH}
\underset{2\times 6r}{\bN^{\alpha}(\bx)} =  \underset{2\times (2i+2b)} {\bN^{h}(\bx)} ~ \underset{(2i+2b)\times 2b} {\tilde{\bM}^{\alpha}\vphantom{_x}} ~ \underset{(2b\times 6r)}{\bPsi\vphantom{_x}},
\end{align}
determined from the product of the B\'ezier interpolation matrix $\bPsi$ and the boundary-interior transformation matrix $\tilde{\bM}^{\alpha}$.

The shape functions $\bN^{\alpha}(\bx)$ on the coarse element $\calM^{\alpha}$ are the assembly of the shape function to each CBN, represented as a $2\times2$ matrix possibly of all non-zero entry values. In significant contrast to the conventional scalar bilinear shape functions where displacement interpolates each coordinate independently, the matrix-valued shape function tightly couples different dimensions and handles anisotropy naturally; the phenomenon was first observed and studied by~\cite{chen2018numerical}.

\begin{figure}[tb]
  \centering
  \includegraphics[width=1\textwidth]{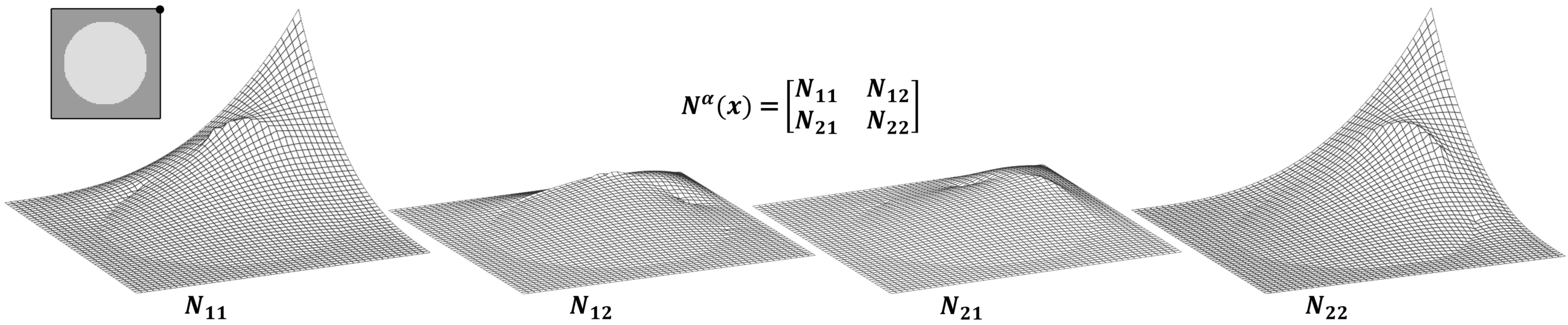}
  \caption{Surfaces of the four components of our matrix-valued CBN shape function $\bN^{\alpha}(\bx)$ for the corner node in solid; the corner nodes are taken as bridge nodes.}
  \label{fig:NH-components}
  \includegraphics[width=0.95\textwidth]{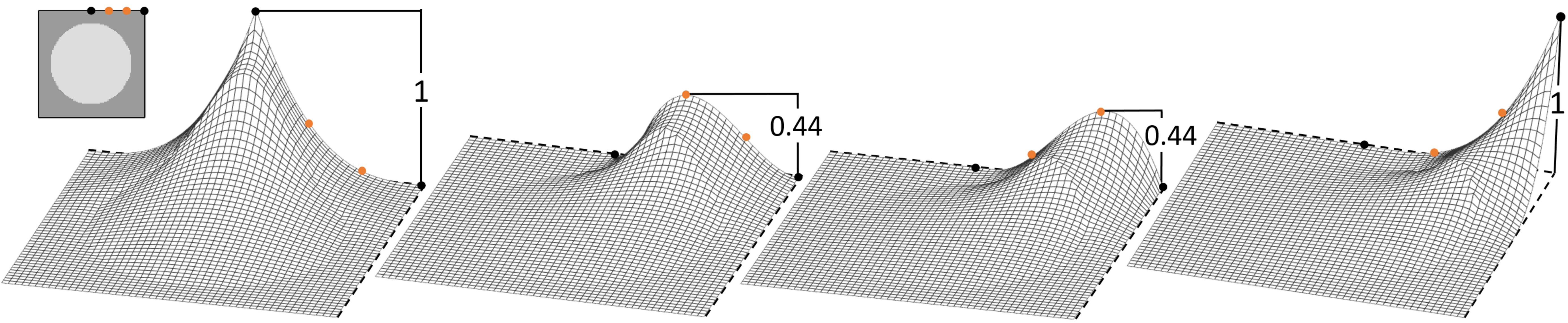}
\caption{Shape functions (in its first component) for the four different CBNs; the example has $8$ bridge nodes.}
\label{fig:sfelement}
\end{figure}

Consider the surface functions in Fig.~\ref{fig:NH-components}, where the bridge nodes are chosen as the corner nodes. The example has $12$ CBNs, and due to its symmetry, we just consider the shape function to a corner and plot surfaces of its four components. We note from the plots that $\bN_{11}(\bx)$ and $\bN_{22}(\bx)$ have much larger height values than those of $\bN_{12}(\bx)$ and $\bN_{21}(\bx)$. This is consistent with our assumption that $\bN_{11}(\bx)$ and $\bN_{22}(\bx)$ play key roles while $\bN_{12}(\bx)$ and $\bN_{21}(\bx)$ regulate the interpolation by coupling different axes. In addition, the complex surfaces in Fig.~\ref{fig:NH-components} within a coarse domain differ significantly from the bilinear shape functions and are expected to better expose the heterogeneity of the coarse element.

\subsection{Geometric properties of CBN shape functions}
The derived shape functions $\bN^{\alpha}(\bx)$ also satisfy the basic geometric properties required by FE shape functions to avoid aphysical behavior, as explained below.

\paragraph{Node interpolation}
Let $\bQ^{\alpha}_j$ be the displacement on a coarse node $\bx_j$. From the CBN shape function expression in Eq.~\eqref{eq:NH} and the fine nodal displacement expression $\bq$ in Eq.~\eqref{eq:fine-q}, we have
\begin{align}
{\bN^{\alpha}(\bx_j)} ~ {\bQ^{\alpha}\vphantom{\bN^{\alpha}(\bx_j)}} = {\bN^{h}(\bx_j)} ~ {\tilde{\bM}^{\alpha}\vphantom{\bN^{\alpha}(\bx_j)}} ~ {\bPsi\vphantom{\bN^{\alpha}(\bx_j)}}~{\bQ^{\alpha}\vphantom{\bN^{\alpha}(\bx_j)}}
= {\bN^{h}(\bx_j)} ~ {\bq\vphantom{\bN^{\alpha}(\bx_j)}}.
\end{align}
As the local fine mesh shape functions $\bN^{h}(\bx)$ naturally satisfy the node interpolation property,
\eb
{\bN^{h}(\bx_j)} ~ {\bq\vphantom{\bN^{\alpha}(\bx_j)}} = {\bQ^{\alpha}_j},
\ee
so does ${\bN^{\alpha}(\bx)}$, or
\eb \label{eq:itp}
{\bN^{\alpha}(\bx_j)\vphantom{\bU^{\alpha}_j}} ~ {\bQ^{\alpha}\vphantom{\bU^{\alpha}_j}} = {\bQ^{\alpha}_j}.
\ee

\paragraph{Partition of unity (PU)} The property is to ensure the property of translation invariant. Note that PU is satisfied for the fine node shape function $\bN^{h}(\bx)$ and the B\'ezier interpolation matrix $\bPsi$, or
\begin{align}
  \bN^{h}(\bx) ~ \bI_{2i+2b,1} &= \bI_{2,1}, \\
  \bPsi ~ \bI_{6r,1} &= \bI_{2b,1}.
\end{align}

Substituting $\bq_b$ and $\bq_i$ with an all-one vector in Eq.~\eqref{eq:micro-interior} gives
\begin{align}
  {\bM}^{\alpha} ~ \bI_{2b,1} = \bI_{2i,1},
\end{align}
and thus,
\begin{align}
  \tilde{\bM}^{\alpha} ~ \bI_{2b,1} = [\bI_{2b}, \bM^{\alpha}] ~ \bI_{2b,1} = \bI_{2i+2b,1}.
\end{align}
Consequently, we have the PU property for $\bN^{\alpha}(\bx)$,
\eb
\begin{aligned}
{\bN^{\alpha}(\bx)} ~ \bI_{6r,1}  = {\bN^{h}(\bx)} ~ {\tilde{\bM}^{\alpha}\vphantom{_x}} ~ {\bPsi\vphantom{_x}} ~ \bI_{6r,1} =
{\bN^{h}(\bx)} ~ {\tilde{\bM}^{\alpha}\vphantom{x}} ~ \bI_{2b,1} = {\bN^{h}(\bx)} ~ \bI_{2i+2b,1} = \bI_{2,1}.
\end{aligned}
\ee

\paragraph{Rotation invariant} Given a constant angular velocity $\hat{\theta}$, we need to show that
\eb
\hat{\theta}\times \bx = \bN^{\alpha}(\bx) ~(\hat{\theta}\times \bx).
\ee
The results can be similarly proved as the case of the translation invariant by noticing that the B\'ezier curve is rotation invariant and that Eq.~\eqref{eq:micro-interior} is also satisfied under a rotation transformation.
\begin{thm}\label{thm-gp} The CBN shape functions $\bN^{\alpha}(\bx)$ in Eq.~\eqref{eq:NH} has the following basic geometric properties,
\begin{align}
  &\bN^{\alpha}(\bx_j) ~ \bQ^{\alpha} = \bQ^{\alpha}_j,~\mbox{for a coarse node}~\bx_j\in \Omega^{\alpha}\\
  &\bN^{\alpha}(\bx) ~ \bI_{6r,1} = \bI_{2,1}, ~ \forall \bx \in \Omega^{\alpha}\\
  &\hat{\theta}\times \bx = \bN^{\alpha}(\bx) ~(\hat{\theta}\times \bx),~ \forall \hat{\theta},~\bx\in\Omega^{\alpha},
\end{align}
where $\bQ^{\alpha}_j$ is the nodal displacement on CBN nodes, $\hat{\theta}$ is a constant angular velocity.
\end{thm}

The above properties are important to produce physically reasonably analysis results using an FE analysis framework. Otherwise, it may, for example, produce a lower global stiffness if PU cannot be satisfied~\cite{nesme2009preserving}. These basic properties are imposed as constraints in an optimization problem in previous study on constructing the material-aware shape functions~\cite{chen2018numerical}. They are naturally satisfied for our CBN shape functions. We summarize the results below.

\begin{thm}\label{thm-sfd}
Given a coarse element $\Omega^{\alpha}$, its local fine mesh $\calM^{\alpha,h}$, and the fine-mesh shape function $\bN^{h}(\bx)$ in~\eqref{eq:Nh}, we have the CBN shape functions
\eb\label{eq-thm-sf}
{\bN^{\alpha}(\bx)} =  {\bN^{h}(\bx)}~ {\tilde{\bM}^{\alpha}}~\bPsi,
\ee
where the boundary-interior transformation matrix $\tilde{\bM}^{\alpha}$ is defined in Eq.~\eqref{eq:bitrans} and the B\'ezier interpolation matrix $\bPsi$ is defined in Eq.~\eqref{eq:bezier-psi}. In addition, the basic geometric properties stated in Theorem~\ref{thm-gp} are satisfied.
\end{thm}

\subsection{Numerical aspects}
\subsubsection{Computation reduction} \label{sec:reduce}
The main computational costs to derive the shape functions $\bN^{\alpha}(\bx)$ in Eq.~\eqref{eq-thm-sf} mainly involve the computation of $\bM^{\alpha}$ in Eq.~\eqref{eq:micro-interior}, or the product of $\bk_i^{-1}$ with $\bk_{ib}$. It is formulated as a solution to the following linear equation systems,
\eb\label{eq-solM}
  \underset{2i\times 2i}{\bk_{i}\vphantom{\bM^{\alpha}_x}}~\underset{2i\times 2b}{\bM^{\alpha}\vphantom{_x}} = -\underset{2i\times 2b}{\bk_{ib}}.
\ee
The column number $2b$ of the right terms can be a very large number, and computing solutions to such a large number of equations would be costly even if using a pre-computed LU decomposition. Specifically, it is usually unaffordable bearing in mind that such equation systems have to be computed for all the different coarse elements $\Omega^{\alpha}$ of different stiffness matrices $\bk_{i}$.

Our special introduction of CBNs and the associated B\'ezier interpolation transformation $\bPsi$ provides an alternative to reduce the computational costs. Multiplying both sides in Eq.~\eqref{eq-solM} by $\bPsi$ yields
\begin{equation}
  \underset{2i\times 2i}{\bk_{i}} ~ \underset{2i\times 2b}{\bM^{\alpha}\vphantom{_x}} ~ \underset{2b\times 6r}{\bPsi\vphantom{_x}}  = -\underset{2i\times 2b}{\bk_{ib}} ~ \underset{2b\times 6r}{\bPsi\vphantom{_x}}.
\end{equation}
Instead of computing $\bM^{\alpha}$, we directly compute the product $\bM^{\alpha}~\bPsi$ as the matrix $\bPhi$ defined below,
\eb \label{eq:division}
  \underset{2i\times 2i}{\bk_{i}} ~ \underset{2i\times 6r}{\bPhi\vphantom{_x}}  = \underset{2i\times 2r}{\bk_{\phi}\vphantom{_x}},\quad \mbox{for} \quad
  \underset{2i\times 2r}{\bk_{\phi}\vphantom{_x}} = -\underset{2i\times 2b}{\bk_{ib}} ~ \underset{2b\times 6r}{\bPsi\vphantom{_x}}.
\ee
Now the number of linear equation systems to be solved is greatly reduced from $2b$ to a much smaller number of $6r$.

According to Eq.~\eqref{eq:bitrans}, by letting
\eb\label{eq-psiphi}
\tilde{\bPhi}=[\bPsi,\bPhi]^T,
\ee
the CBN shape function now takes the following form,
\eb
{\bN^{\alpha}(\bx)} =  {\bN^{h}(\bx)}~\tilde{\bPhi}.
\ee

\begin{cor} The CBN shape functions $\bN^{\alpha}(\bx)$ in Theorem~\ref{thm-sfd} can be numerically derived as follows
\eb\label{eq:sfbphi}
{\bN^{\alpha}(\bx)} =  {\bN^{h}(\bx)}~\tilde{\bPhi},
\ee
where $\tilde{\bPhi}$ is defined in Eq.~\eqref{eq-psiphi}.
\end{cor}

\subsubsection{Heterogeneous structure analysis on CBNs}
Once the CBN shape functions $\bN^{\alpha}(\bx)$ have been constructed for each coarse element $\Omega^{\alpha}$, computing the displacement to the linear elasticity problem in~\eqref{eq-elasticity}, on a heterogenous structure $\Omega$, can then be achieved following a traditional FE analysis framework. The overall algorithm is described in Algorithm~\ref{alg:frame}.

Let $\bQ^{\alpha}$ be the vector of discrete nodal displacements on $\calM^H$ to be determined. Then the continuous displacement $\bu^{\alpha}(\bx)$ is interpolated using the CBN shape functions $\bN^{\alpha}(\bx)$,
\eb
\bu^{\alpha}(\bx)= \bN^{\alpha}(\bx)~\bQ^{\alpha}.
\ee
Accordingly, we have the strain,
\eb\label{eq-disc-strain}
{\bveps^{\alpha}(\bx)\vphantom{_x}}={\bB^{\alpha}(\bx)\vphantom{_x}}~{\bQ^{\alpha}},
\ee
where $\bB^{\alpha}(\bx)$ is the derivative of $\bN^{\alpha}(\bx)$ with respect to $\bx$,
\eb\label{eq-bexp}
\bB^{\alpha}(\bx) = \frac{\partial \bN^{\alpha}(\bx)}{\partial \bx} = \frac{\partial \bN^{h}(\bx)}{\partial \bx} ~ \tilde{\bPhi},
\ee
according to Eq.~\eqref{eq:sfbphi}.

Substituting Eq.~\eqref{eq-disc-strain} into the weak formulation Eq.~\eqref{eq-weak}, the coarse nodal displacement $\bQ$ is then
computed as,
\eb\label{eq-kuf}
\bK~\bQ=\bff, \quad \text{for} \quad \bK=\sum_{\alpha = 1}^M \bK^{\alpha},
\ee
where
\begin{eqnarray}\label{eq:kalpha}
\bK^{\alpha} &= & \int_{\Omega^{\alpha}} (\bB^{\alpha}(\bx))^T~\bD(\bx)~\bB^{\alpha}(\bx)~\rd\Omega^{\alpha}\nonumber  \\
&=&\sum_{e=1}^m\int_{\omega_e} (\bB^{\alpha}(\bx))^T~\bD_e(\bx)~\bB^{\alpha}(\bx)~\rd\omega_e\nonumber \\
&=& \sum_{e=1}^m  \tilde{\bPhi}^T(\int_{\omega_e}(\frac{\partial \bN^{h}(\bx)}{\partial \bx})^T~\bD_e(\bx)~\frac{\partial \bN^{h}(\bx)}{\partial \bx}~\rd\omega_e)~ \tilde{\bPhi}\nonumber\\
&=& \sum_{e=1}^m\bK^{\alpha}_e,
\end{eqnarray}
for
\eb \label{eq:k-alpha-e}
\bK^{\alpha}_e =\tilde{\bPhi}^T(\int_{\omega_e}(\frac{\partial \bN^{h}(\bx)}{\partial \bx})^T~\bD_e(\bx)~\frac{\partial \bN^{h}(\bx)}{\partial \bx}~\rd\omega_e)~ \tilde{\bPhi},
\ee
by noticing Eq.~\eqref{eq-bexp}.

\paragraph{Numerical computation for $\bK^{\alpha}$} The computation of $\bK^{\alpha}$ involves an integration computation on a coarse element $\Omega^{\alpha}$. For a homogeneous $\Omega^{\alpha}$, it is to be achieved using Gauss integration at $2^d$ Gauss points, $d=2,3$. The integration here, however, works on a heterogeneous coarse element $\Omega^{\alpha}$ and a piecewise shape function $\bN^{\alpha}(\bx)$.
To achieve computation accuracy, it thus has to be conducted on each fine element $\omega_e$ and assembled together following Eq.~\eqref{eq:kalpha}; we further take $2^d$ Gauss points, $d=2,3$, for each fine element. Note also that $\bN^{h}(\bx)$ is locally supported, and the numerical integration in Eq.~\eqref{eq:k-alpha-e} only involves $\bN_e(\bx)$ instead of $\bN^{h}(\bx)$ for a specific fine element $\omega_e$.


\begin{cor} Let $\Omega$ be a heterogenous solid structure, $\calM^{H},\calM^{\alpha,h}$ are its coarse mesh and fine mesh, respectively. The solution $\bQ$ to the linear elasticity analysis problem in Eq.~\eqref{eq-elasticity} can be computed as
\eb\label{eq:CKQ}
\bK~\bQ=\bff,
\ee
for
\eb
\bK=\sum_{\alpha = 1}^M \bK^{\alpha}= \sum_{\alpha = 1}^M \sum_{e = 1}^m  \bK^{\alpha}_e,
\ee
and $\bK^{\alpha}_e$ given in~\eqref{eq:k-alpha-e}.
\end{cor}

\begin{algorithm}
  \caption{Heterogeneous structure analysis using CBN}
  \textbf{Input}: a heterogeneous structure $\Omega$, its coarse mesh $\calM^H=\{\Omega^{\alpha},\ \alpha=1,2,\ldots, M\}$, fine mesh $\calM^h= \{\omega^{\alpha}_e,\ e=1,2,\ldots, m,\ \alpha = 1,2,\ldots,M \}$. \\
  \textbf{Output}: the CBN shape functions $\bN^{\alpha}(\bx)$, and the approximated displacement $\bQ$ to Eq.~\eqref{eq-elasticity}
  \begin{algorithmic}[1]
    \State \textbf{Prepare} B\'ezier interpolation matrix $\bPsi$ in Eq.~\eqref{eq:bezier-psi}
    \State \textbf{Prepare} boundary-interior transformation matrix $\tilde{\bM}^{\alpha}$ for each coarse element $\Omega^{\alpha}$ in Eq.~\eqref{eq:bitrans}
    \State \textbf{Construct} the transformation matrix $\tilde{\bPhi}$ in Eq.~\eqref{eq:division}
    \State \textbf{Construct} the shape functions $\bN^{\alpha}(\bx)$ for each $\Omega^{\alpha}$ in Eq.~\eqref{eq:sfbphi}
    \State \textbf{Compute} the elemental stiffness matrix $\bK^{\alpha}$ for each $\Omega^{\alpha}$ in Eq.~\eqref{eq:kalpha}
    \State \textbf{Assemble} the global stiffness matrix $\bK$
    \State \textbf{Compute} the displacement solution $\bQ$ by Eq.~\eqref{eq:CKQ}
  \end{algorithmic}
  \label{alg:frame}
\end{algorithm}

\section{Extension to 3D cases and nonlinear analysis}\label{sec:extension}

\subsection{Extension to 3D cases}\label{sec:3d}
\subsubsection{3D B\'ezier interpolation matrix $\bPsi$}
The 3D case takes a bicubic B\'ezier surface in the following form:
\eb\label{eq-cbezier-3d}
\bP(\bt)=\sum_{i=0}^3\sum_{j=0}^3\psi_{i}(t_u)~ \psi_{j}(t_v) ~ \bP_{ij}=\sum_{i=0}^3\sum_{j=0}^3  C^i_3t_u^i(1-t_u)^{3-i}~ C^j_3t_v^j(1-t_v)^{3-j}~\bP_{ij},
\ee
where $\bt=(t_u,t_v)$, $\bP_{ij} \in \RR^3$ are the control points, $\psi_{i}(t_u)$ or $\psi_{j}(t_v)$ is a cubic Bernstein basis, and $C_3^i$ is the binomial coefficient.

We can similarly write $\bP(\bt)$ in a matrix form as
\eb\label{eq-bcmatrix-3d}
\underset{(3\times 1)}{\bP(\bt)\vphantom{\bpsi_F}}=\underset{(3\times 48)}{\bpsi_F(t_u, t_v)}~\underset{(48\times 1)}{\bP\vphantom{\bpsi_F}}
\ee
for
\eb
\bP=[\bP_0,\bP_1,...,\bP_{15}]^T, \quad
\bpsi_F={\bpsi_1(t_u)}~{\bpsi_2(t_v)},
\ee
and
\eb
\underset{(3\times 12)}{\bpsi_1(t_u)} = \underset{(1\times 4)}{\bpsi_0(t_u)} \otimes \underset{(3\times 3)}{\bI_3\vphantom{\bpsi_0}},\quad \underset{(12\times 48)}{\bpsi_2(t_v)} = \underset{(1\times 4)}{\bpsi_0(t_v)} \otimes \underset{(12\times 12)}{\bI_{12}\vphantom{\bpsi_0}},
\ee
where $\bpsi_0(t)$ is referred in Eq.~\eqref{eq:basis}

\begin{figure}[tb]
  \centering
  \subfigure[Corner nodes as bridge nodes, $\calV_c = \calV_r$]{\includegraphics[width=0.35\textwidth]{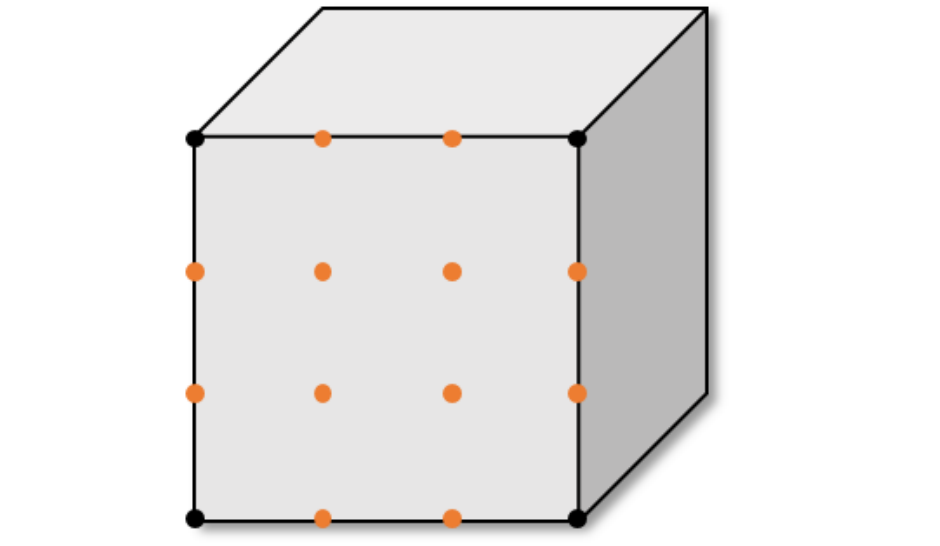}}
  \subfigure[$3\times 3$ bridge nodes, $\calV_c \subset \calV_r$]{\includegraphics[width=0.35\textwidth]{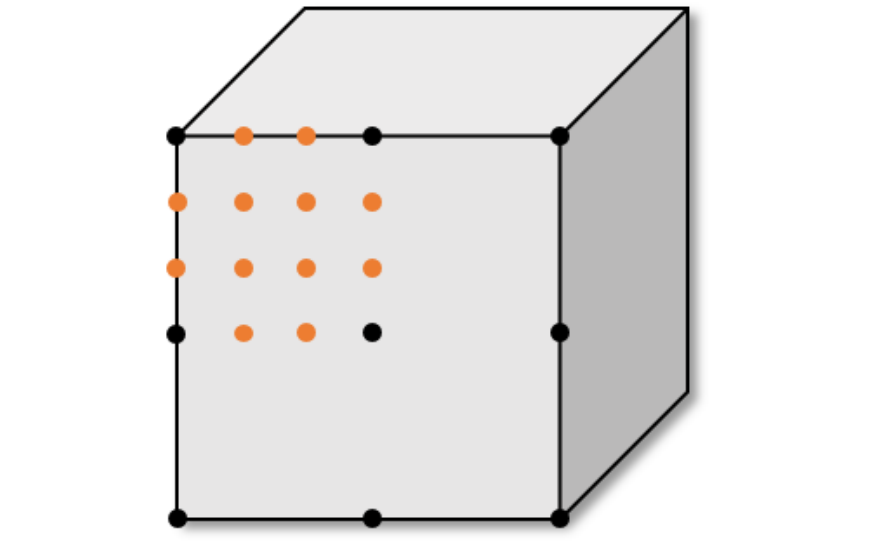}}
  \caption{Strategies to set CBN in 3D cases for shape functions construction.}
  \label{fig-bezier-3d}
\end{figure}

Given a bridge face $F$ determined by a set of bridge nodes, we have similarly as the 2D case the displacement on any point $\bx$ on face $F$,
\eb
\bu_{F}(\bx) = \bu(t_u(\bx), t_v(\bx))={\bpsi_F(t_u(\bx), t_v(\bx))}~{\bQ_F},
\ee
where $\bQ_F$ is the vector of the displacements on the associated CBNs as indicated in Fig.~\ref{fig-bezier-3d}.

Assembling $\bu_{F}(\bx)$ over all faces $F$ and evaluating it on all the boundary nodes in $\calV_b$ , we have the 3D B\'ezier interpolation matrix as
\eb
\label{eq:bezier-psi-3d}
\underset{(3b\times q)}{\bPsi}=\left[\sum_F\bpsi_F(\bt(\bx_0)), \ \bx_0\in \calV_b\right],
\ee
where $q$ is the number of DOFs on CBNs given as
\eb \label{eq:3d_dofs_num}
q = 3 ~ (54 x^2 - 108x + 56)\quad \mbox{for} \quad x = 1+\frac{1}{6}\sqrt{6r-12}.
\ee

\subsubsection{Boundary-interior transformation matrix $\tilde{\bM}^{\alpha}$}
In the 3D case, the shape function matrix $\bN_e(\bx)$ has a dimension of $3\times 24$ in the following form:
\eb
\bN_e(\bx)=
\left[
\begin{array}{cccccccc}
\bN_1(\bx) & \bN_2(\bx)  & \bN_3(\bx) & \bN_4(\bx) & \bN_5(\bx) & \bN_6(\bx) & \bN_6(\bx) & \bN_8(\bx)
\end{array}
\right],
\ee
where each submatrix $\bN_i(\bx)$ is
\eb
{\bN_i(\bx)}=
\left[
\begin{array}{cccccccc}
N_i(x) & 0  & 0  \\
0 & N_i(y) & 0 \\
0 & 0 & N_i(z)
\end{array}
\right],
\ee
and $N_i(\bx)$ is the \emph{trilinear shape function} defined on each of the eight corner nodes,
\eb\label{eq-linear-sf}
N_i(\bx)=\frac{1}{8}(1+x_ix)(1+y_iy)(1+z_iz) \quad, 1\leq i\leq 8,\quad (x,y,z)\in [-1,1]\times[-1,1]\times[-1,1].
\ee
Following a similar approach in 2D case, the boundary-interior transformation matrix $\tilde{\bM}^{\alpha}$ can then be derived.

\subsection{Extension to nonlinear analysis}\label{sec:nonlinear}
The approach works directly for nonlinear analysis. We just need to replace the bilinear or trilinear shape functions for each coarse element $\Omega^{\alpha}$
with our CBN shape functions $\bN^{\alpha}(\bx)$ in the deformation gradient $\bF(\bx)$, as shown below:
\eb
\bF(\bx) = \frac{\partial \bu(\bx)}{\partial \bX} + \bI_d =  \frac{\partial \bN^{\alpha}(\bx)}{\partial \bX}\bQ^{\alpha} + \bI_d,
\ee
where $\bx$ is in the deformed shape, $\bX$ is in the reference shape, and $\bI_d$ is a $d \times d$ identity matrix for $d=2,3$.

Conducting nonlinear elasticity analysis using the deformation gradients $\bF(\bx)$ follows a classical FE analysis process. More details can be found in~\cite{bower2009applied}.

\paragraph{Accuracy improvement in nonlinear cases} To improve the elastic behaviors of the nonlinear model, the following corotational formulation~\cite{chen2018numerical} of the displacement is taken into consideration,
\eb
\label{eq:local-frame}
\bu(\bx) = \bR^{\alpha} (\bX + \bN^{\alpha}(\bx) (\bR^{\alpha} \bx - \bX)) - \bX,
\ee
where $\bR^{\alpha}$ is the local frame to each shape function~\cite{chen2018numerical}.

\section{Discussions}\label{sec:other_methods}
We further discuss the relation and difference between the proposed CBN analysis approach (denoted Our-CBN) and previous classical approaches for heterogeneous structure analysis:  Homogenization, FE$^2$ using Voigt-Taylor model~\cite{tan2020direct}, the second-order CMCM~\cite{le2020coarse} (CMCM for short) and Substructuring~\cite{wu2019topology}.

\begin{figure}[tb]
  \centering
  \subfigure[Our CBN approach]{\includegraphics[width=0.24\textwidth]{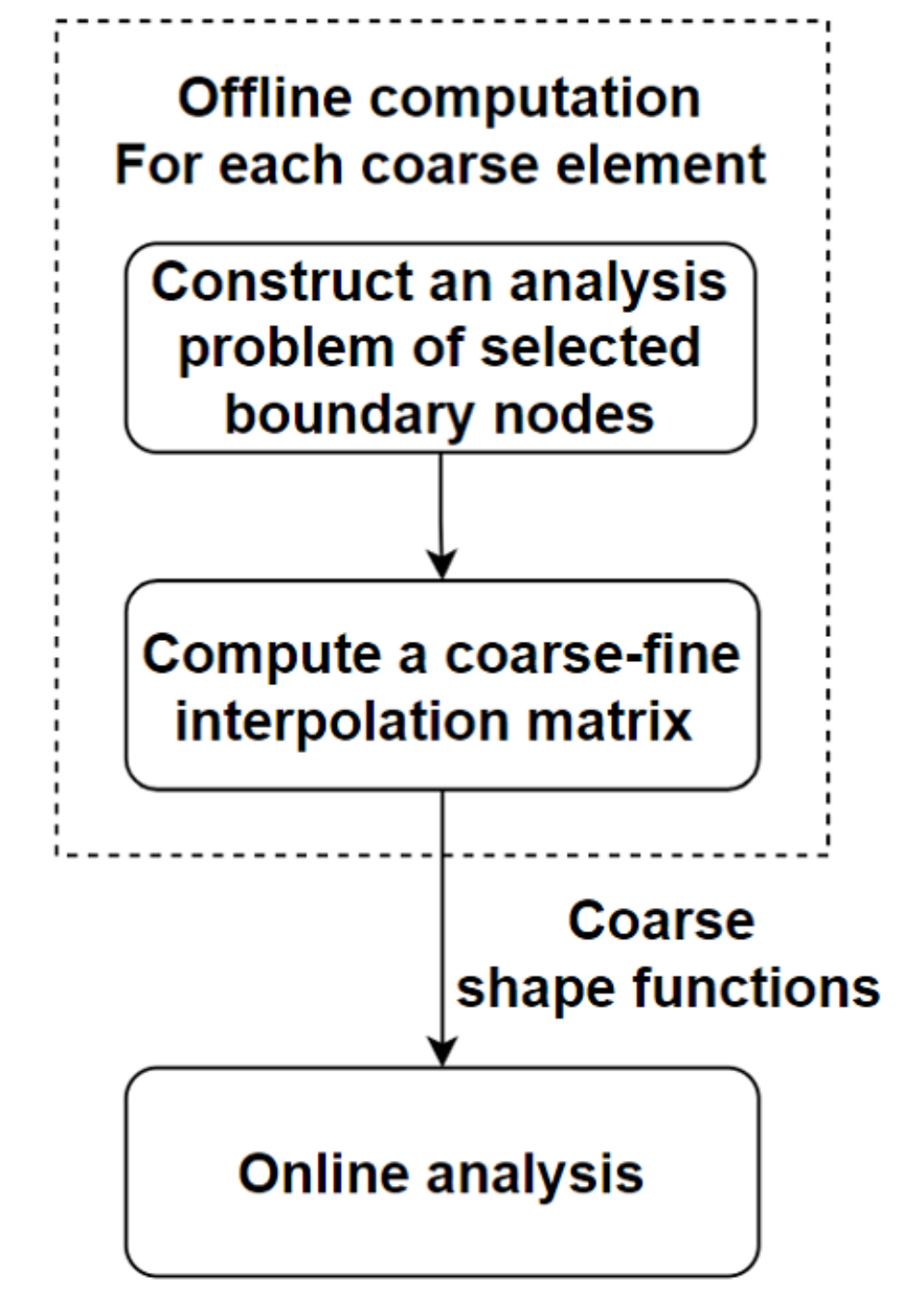}}
  \subfigure[Homogenization]{\includegraphics[width=0.24\textwidth]{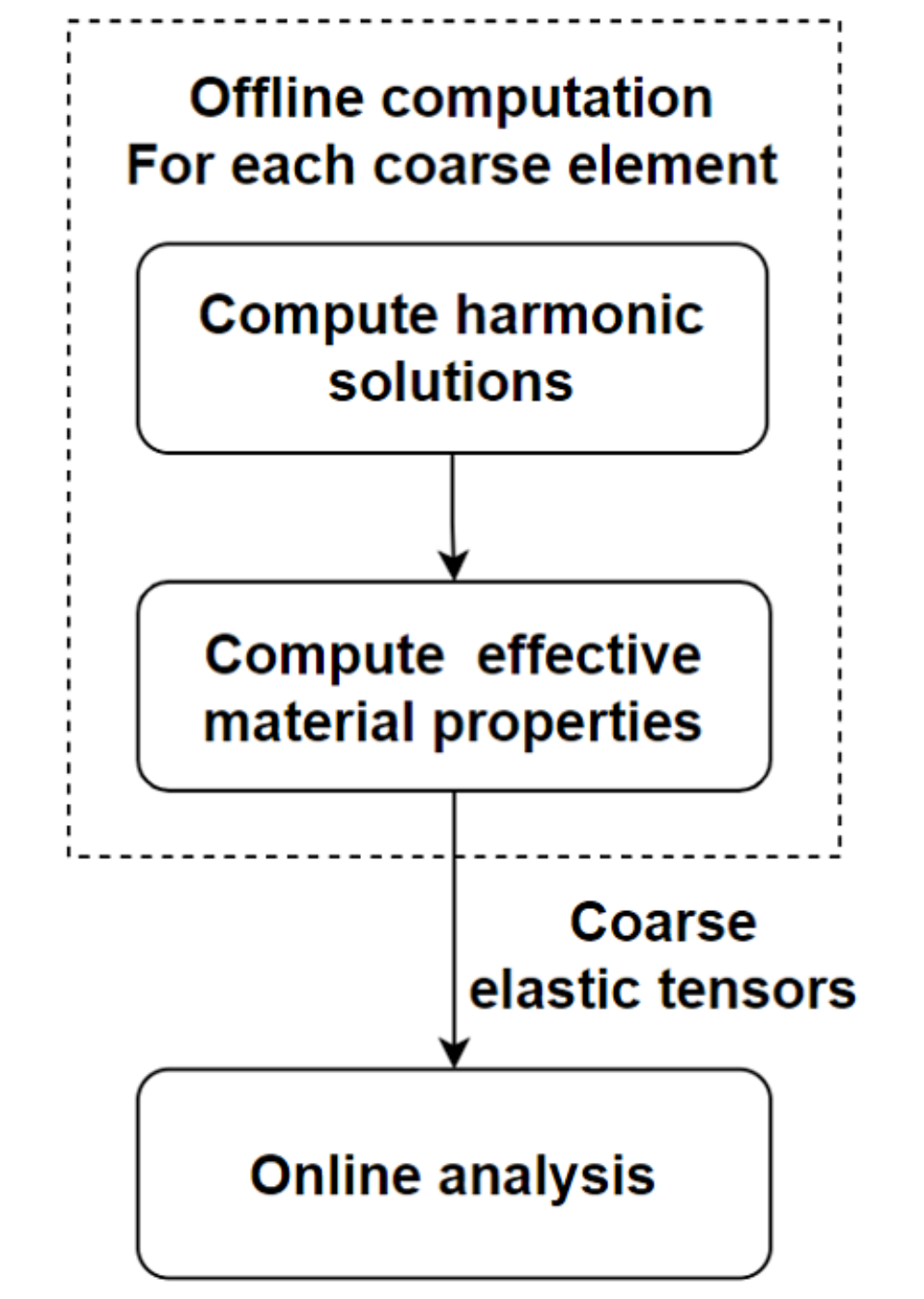}}
  \subfigure[CMCM]{\includegraphics[width=0.24\textwidth]{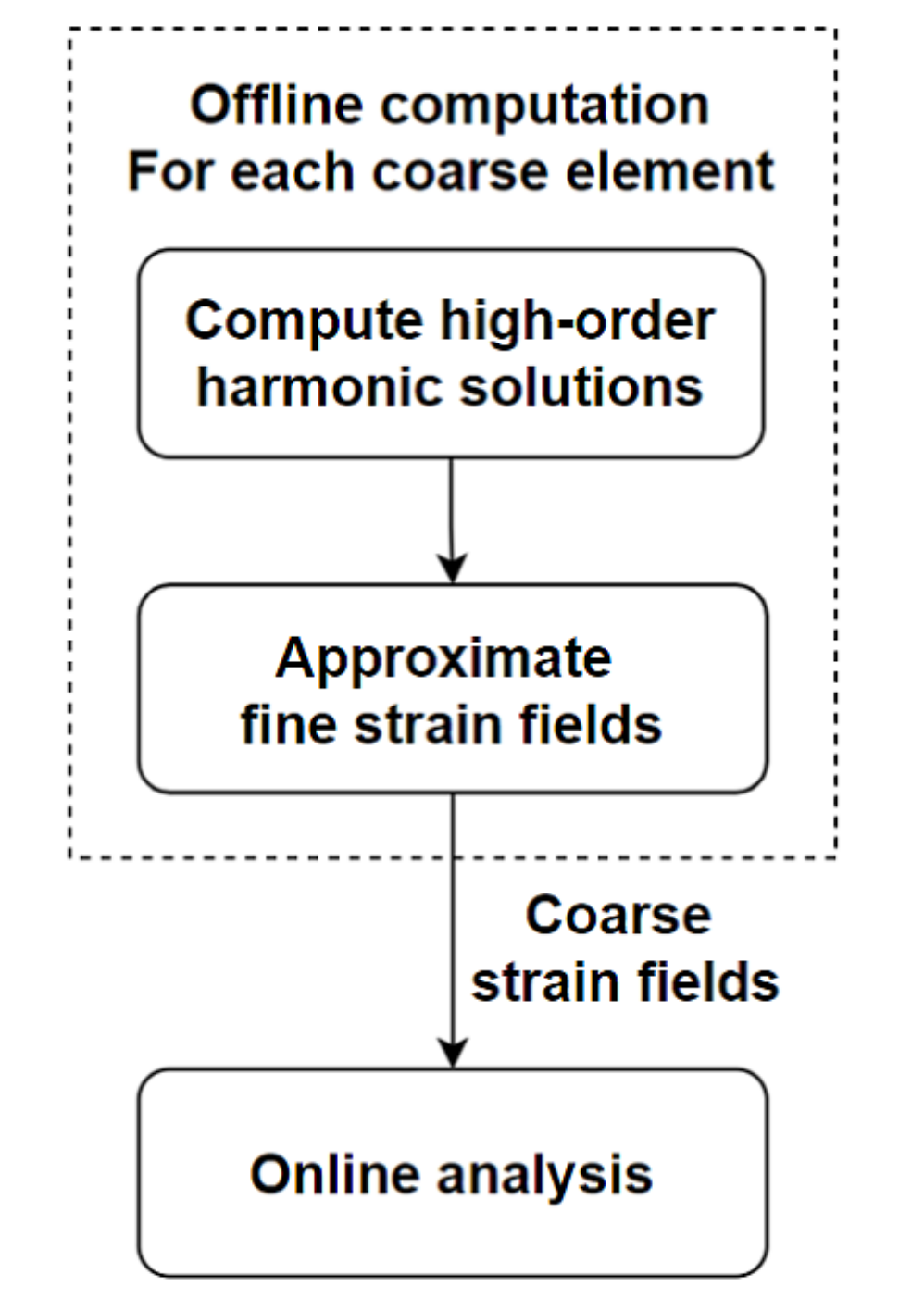}}
  \subfigure[Substructuring]{\includegraphics[width=0.24\textwidth]{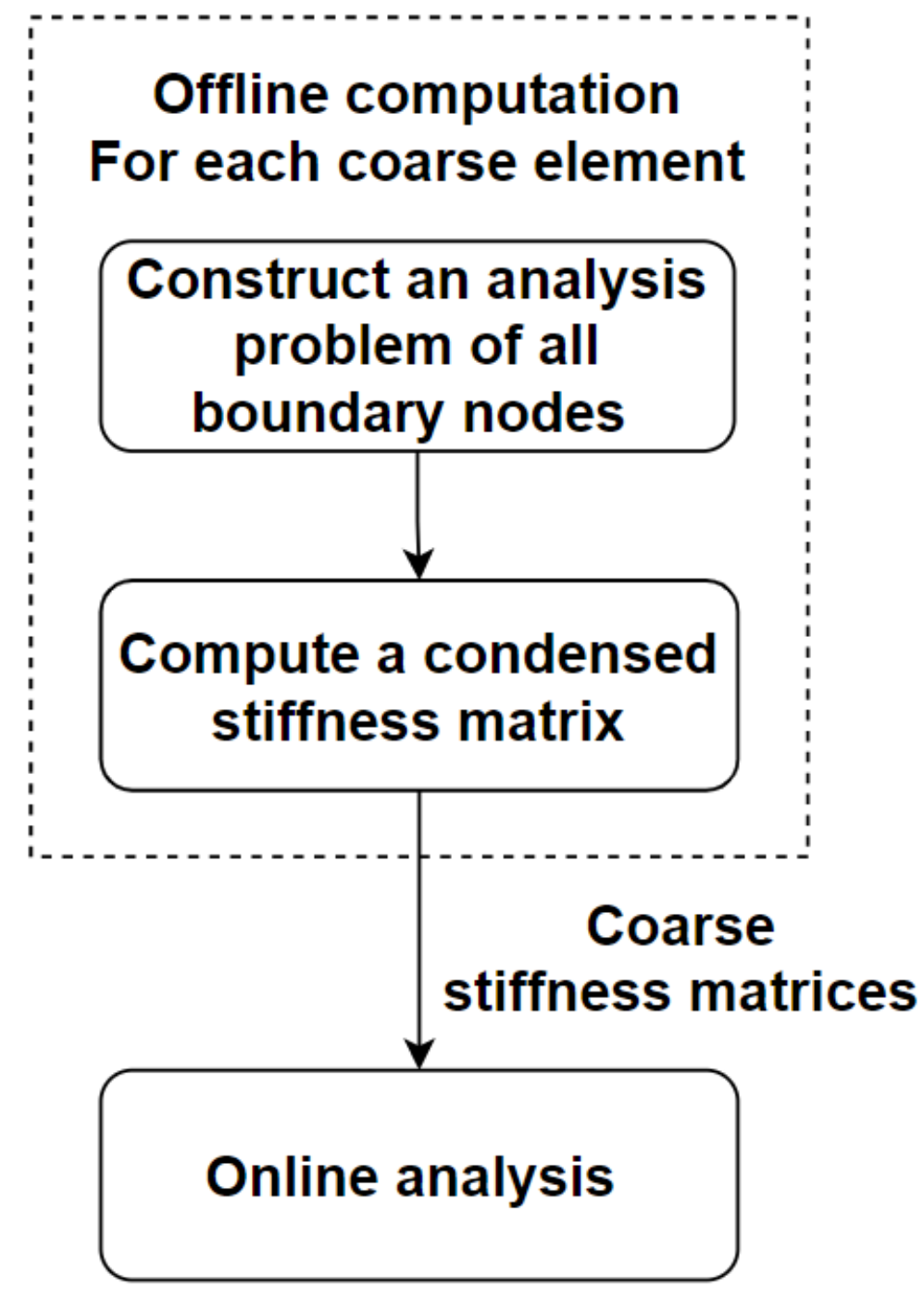}}
  \subfigure[FE$^2$]{\includegraphics[width=0.35\textwidth]{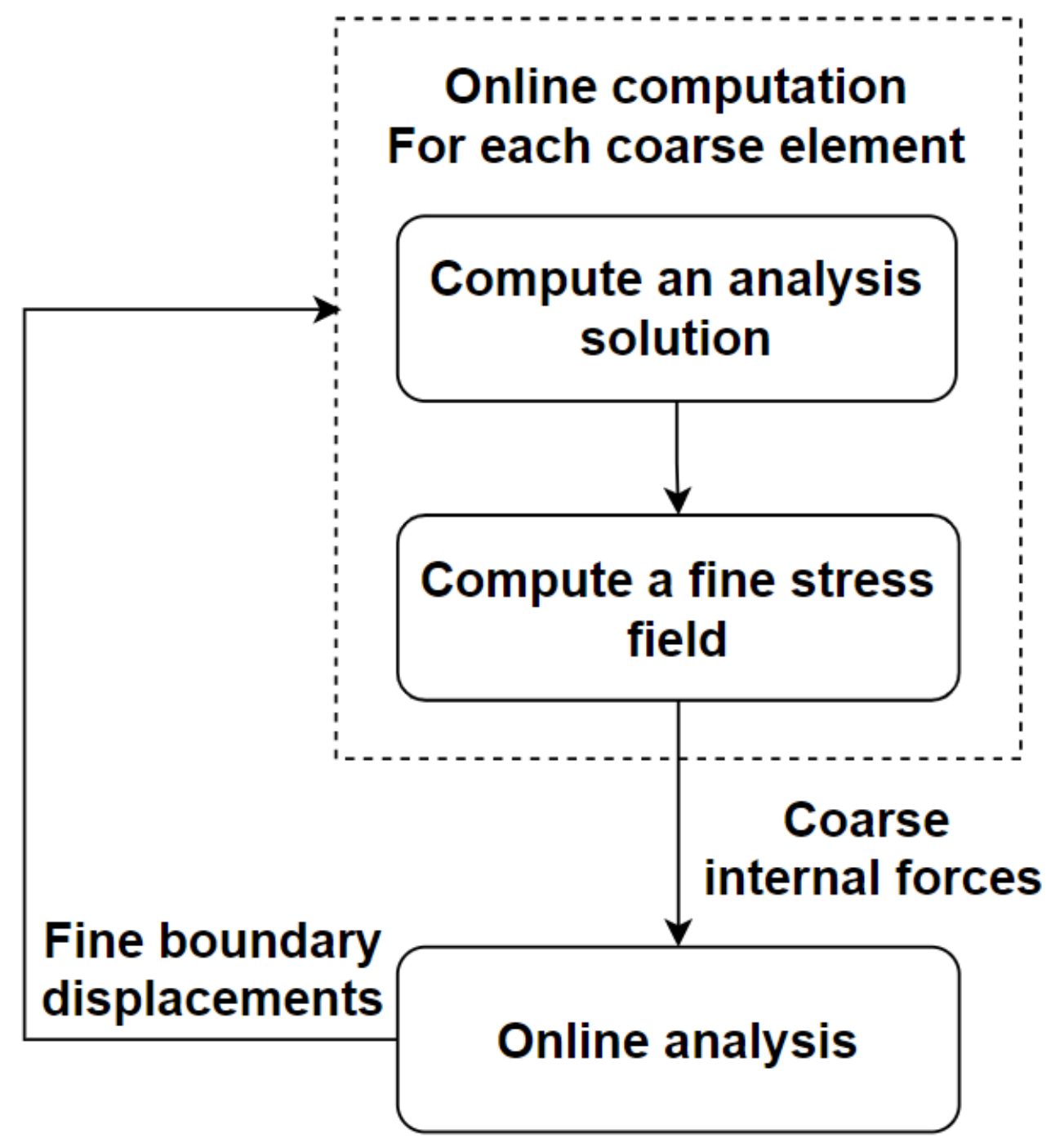}}
  \caption{Flowcharts of different approaches for heterogeneous structure analysis, including Our-CBN, homogenization, CMCM, substructuring and FE$^2$.}
  \label{fig:flowchart}
\end{figure}


\subsection{Technical differences}
Most approaches, including homogenization, CMCM, substructuring and our CBN, generally conduct the analysis following two main procedures: an offline process to compute local fine mesh displacement for each coarse element, and an online process to conduct the global analysis on the coarse mesh. FE$^2$ is different, where both local analysis and global analysis are iteratively conducted online. See also Fig.~\ref{fig:flowchart} for an illustration of the flowcharts of these approaches, with their differences explained below.

\paragraph{Scale separation assumption} Homogenization assumes scale separation as a precondition and may lose much of its analysis accuracy if the assumption is broken. Other approaches, including FE$^2$ of Voigt-Taylor model, substructuring, CMCM and Our-CBN, are not based on the assumption.

\paragraph{Local computations} All these local computations involve solution computations to linear equation systems of the same left-hand stiffness matrix. The right hand is different in two aspects: the number of columns and the entry values. The column number is determined by the involved coarse analysis DOFs. The entry values are calculated from the imposed boundary conditions: testing BCs in homogenization and CMCM, linear interpolation of coarse corner displacements in FE$^2$, or a submatrix of the local fine-mesh stiffness matrix in substructuring and Our-CBN.

\paragraph{Fine--coarse transmission} Homogenization or FE$^2$ transmits specific physical quantities for global coarse-mesh analysis, such as an effective material elasticity tensor and internal force vector. FE$^2$ attempts to further improve the analysis accuracy via iterative computations between the local fine mesh and coarse mesh, and thus may encounter a convergence issue. By contrast, substructuring, CMCM or Our-CBN constructs an explicit physical field, specifically strain fields or the derivatives of the shape functions, local stiffness matrix of the super-elements, and CBN shape functions. They are used to generate the local stiffness matrix to the coarse element.

\paragraph{Global solution reconstruction} Given the coarse mesh solution, reconstructing the local response in the local fine mesh is important for various industrial applications. Homogenization only computes the coarse displacement and is not directly applicable to recover the fine-mesh displacement. Other approaches studied here are all able to reconstruct the fine-mesh strain or stress field directly. In terms of global displacement smoothness, CMCM is not continuous along the coarse boundary as its coarse strain fields are constructed locally without adjacency consideration. All the other three approaches, substructuring, FE$^2$ and Our-CBN, are able to generate a globally smooth displacement, although they each have significantly different computational costs, as discussed later in this paper.

\begin{table}[]
  \centering
  \caption{Computational complexity comparisons of different heterogeneous structure analysis approaches: Homogenization, FE$^2$, CMCM, Substructuring and Our-CBN. }
  \label{tab:complexity}
  \begin{threeparttable}
  \begin{tabular}{c|c|c}
    \hline
    \rowcolor[HTML]{C0C0C0}
    {\color[HTML]{333333} Methods} & \begin{tabular}[c]{@{}c@{}}Number of linear equation systems\\for each coarse element (2D or 3D) \end{tabular} & \begin{tabular}[c]{@{}c@{}}DOFs of each coarse element\\in online analysis (2D or 3D) \end{tabular} \\ \hline
    \rowcolor[HTML]{FFFFFF}
    Homogenization                 & 3 ~\mbox{or}~6              & 8~\mbox{or}~ 24   \\
    \rowcolor[HTML]{EFEFEF}
    FE$^2$                         & 1 for one iteration        & 8~\mbox{or}~24                     \\
    \rowcolor[HTML]{FFFFFF}
    CMCM                           & 5 ~\mbox{or}~ 15                     & 8~\mbox{or}~24                     \\
    \rowcolor[HTML]{EFEFEF}
    Substructuring                 & $2b$ ~\mbox{or}~ $3b$                & $2b$~\mbox{or}~$3b$                     \\
    \rowcolor[HTML]{FFFFFF}
    Our-CBN                            & $6r$ ~\mbox{or}~ q              & 6r~\mbox{or}~q                   \\ \hline
  \end{tabular}
  \begin{tablenotes}
    \footnotesize
    \item[*] $b$ is the number of boundary nodes in $\calV_b$, $r$ is the number of coarse bridge nodes in $\calV_r$, and $q$ is in Eq.~\eqref{eq:3d_dofs_num}.
  \end{tablenotes}
\end{threeparttable}
\end{table}

\subsection{Complexity analysis}\label{sec:complexity}
The complexity mainly depends on two aspects: local displacement computation to each coarse element and the global displacement computation on the coarse mesh, as further analyzed below.

\subsubsection{Complexity analysis of local displacement computations}
For all the approaches mentioned above, the local analysis problem involves all the efforts to compute displacements $\bar{\bq}$ to a set of linear equation systems as
\eb\label{eq-com}
\bk_{i}~\bar{\bq} = \bar{\bf},
\ee
where $\bk_{i}$ is the submatrix defined in Eq.~\eqref{eq:micros}, and $\bar{\bf}$ represents a set of column vectors.

Let $\kappa$ be the number of the vectors, which determines the number of equation systems to be computed and consequently the computational complexity; see also Table~\ref{tab:complexity}.
For homogenization and CMCM, $\kappa$ depends on the number of testing boundaries: 3 in 2D or 6 in 3D in homogenization, and 5 in 2D and 15 in 3D in CMCM (as CMCM imposes high-order boundary conditions). With reference to substructuring, we have $\kappa=2b$ in 2D or $3b$ in 3D, which typically can be as high as tens of thousands.
For Our-CBN, we have $\kappa=6r$ in 2D or $q$ (in Eq.~\eqref{eq:3d_dofs_num}) in 3D, which is typically in hundreds. In FE$^2$, $\kappa$ is equal to the number of iterations; it only has a single column vector in each iteration.

We also point out that these different linear equation systems all have the same left-hand matrix, and thus the displacements to different right-hand column vectors can be efficiently computed by performing KL-decomposition in advance. However, this strategy does not apply for the FE$^2$ directly.

\begin{figure}[tbh]
  \centering
  \subfigure[Homogenization, FE$^2$ and CMCM]{\includegraphics[width=0.3\textwidth]{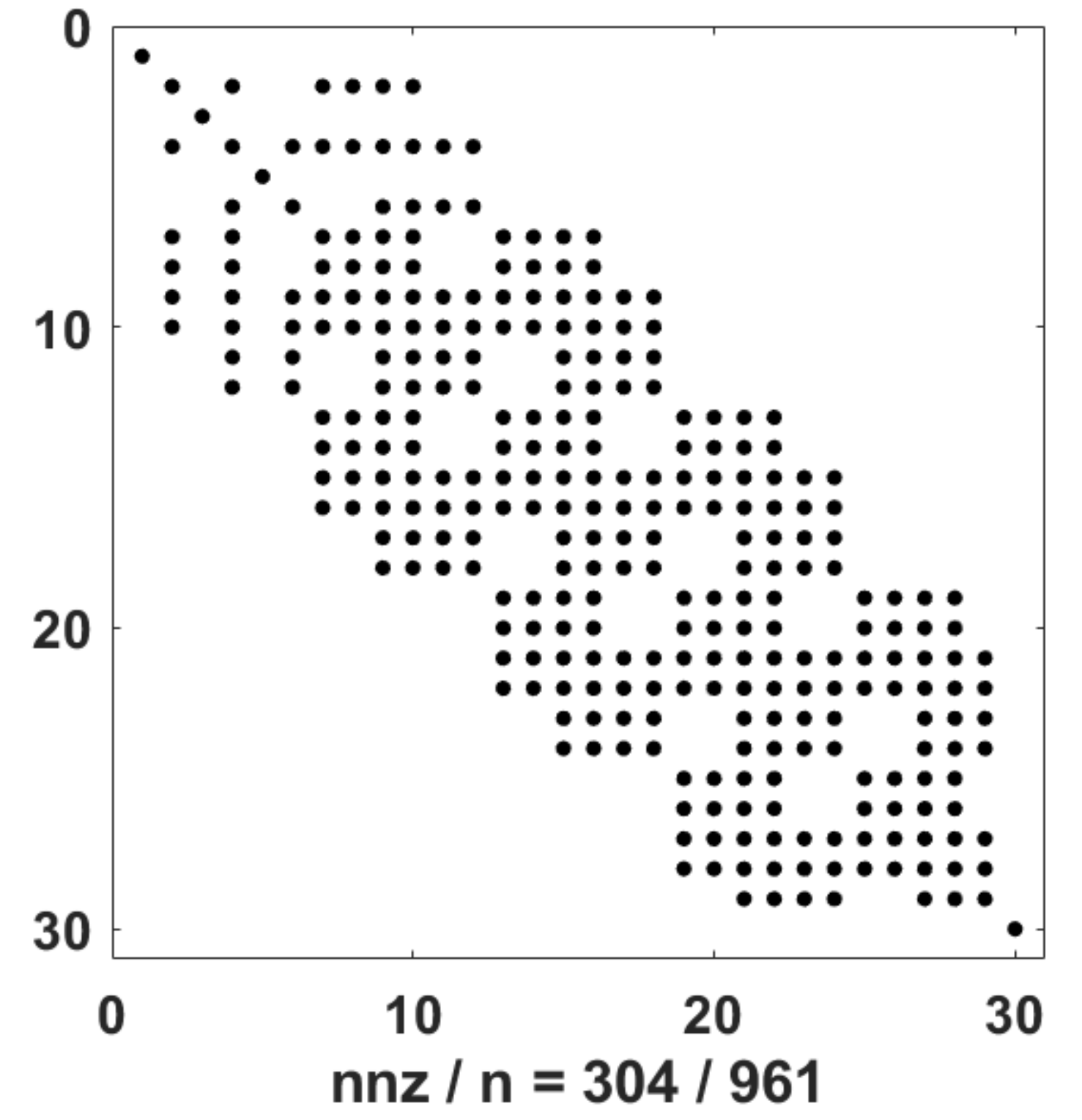}}\quad
  \subfigure[Our-CBN]{\includegraphics[width=0.3\textwidth]{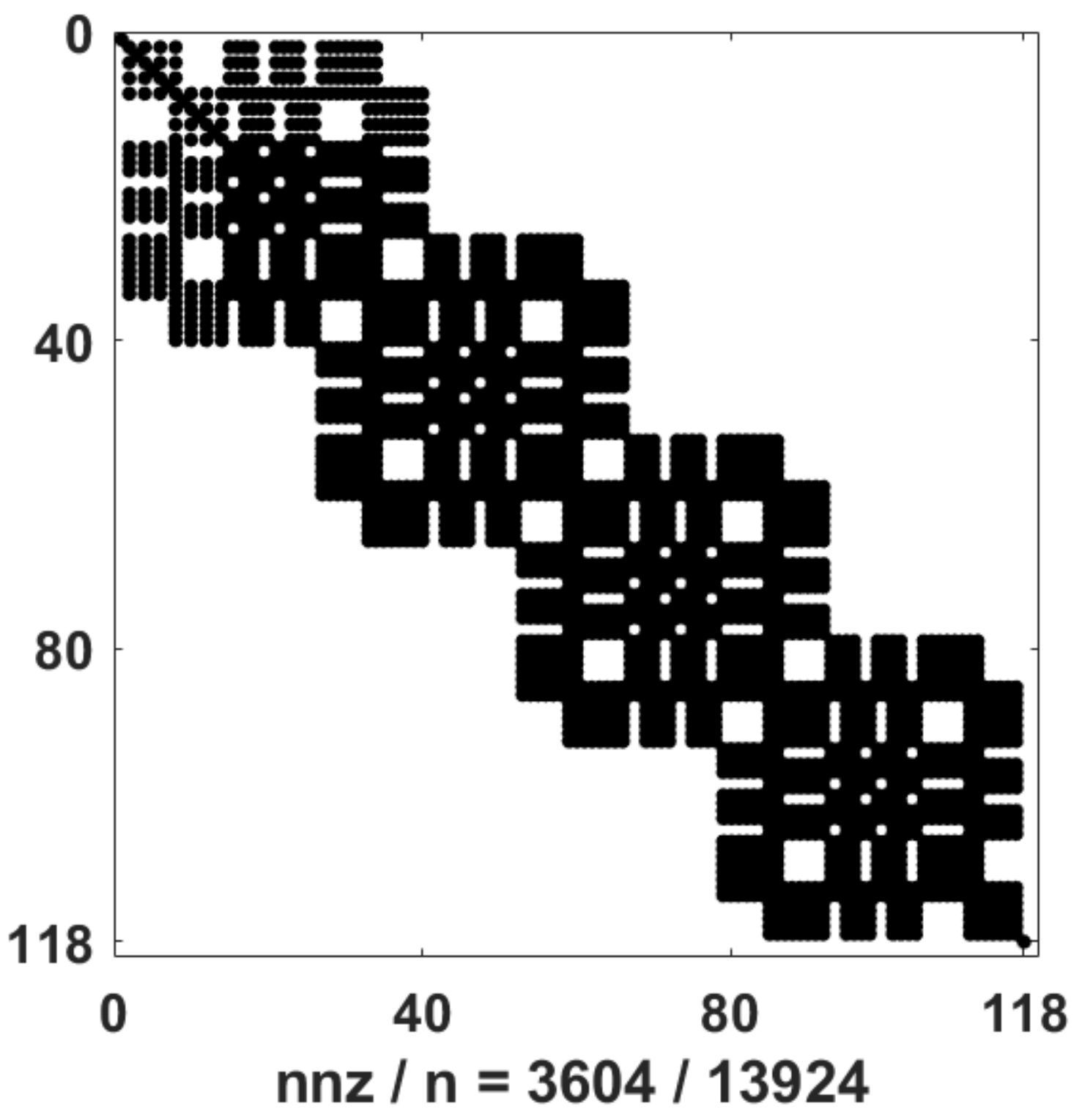}}\quad
  \subfigure[Substructuring]{\includegraphics[width=0.3\textwidth]{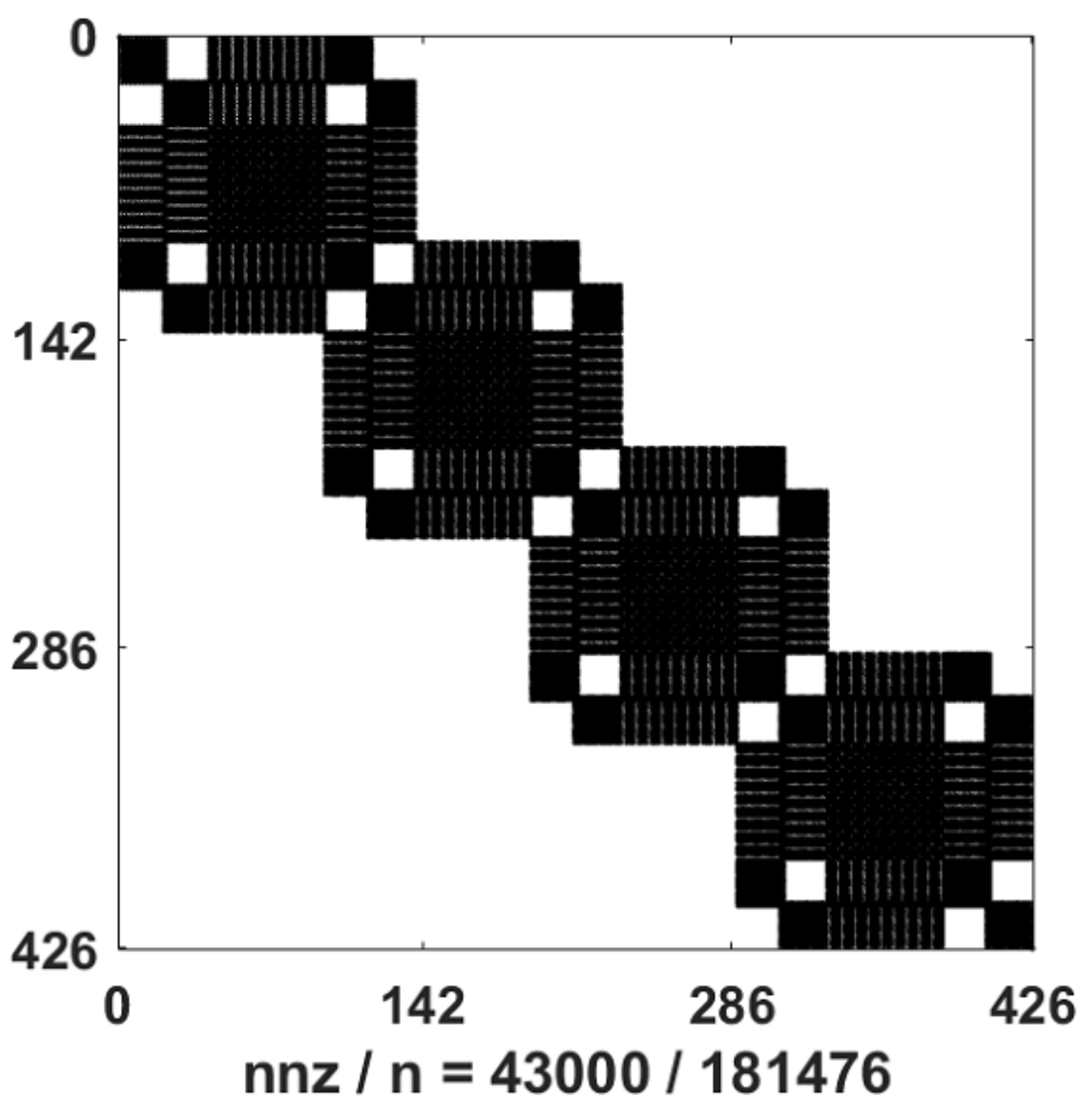}} \\
  \subfigure[Benchmark]{\includegraphics[width=0.58\textwidth]{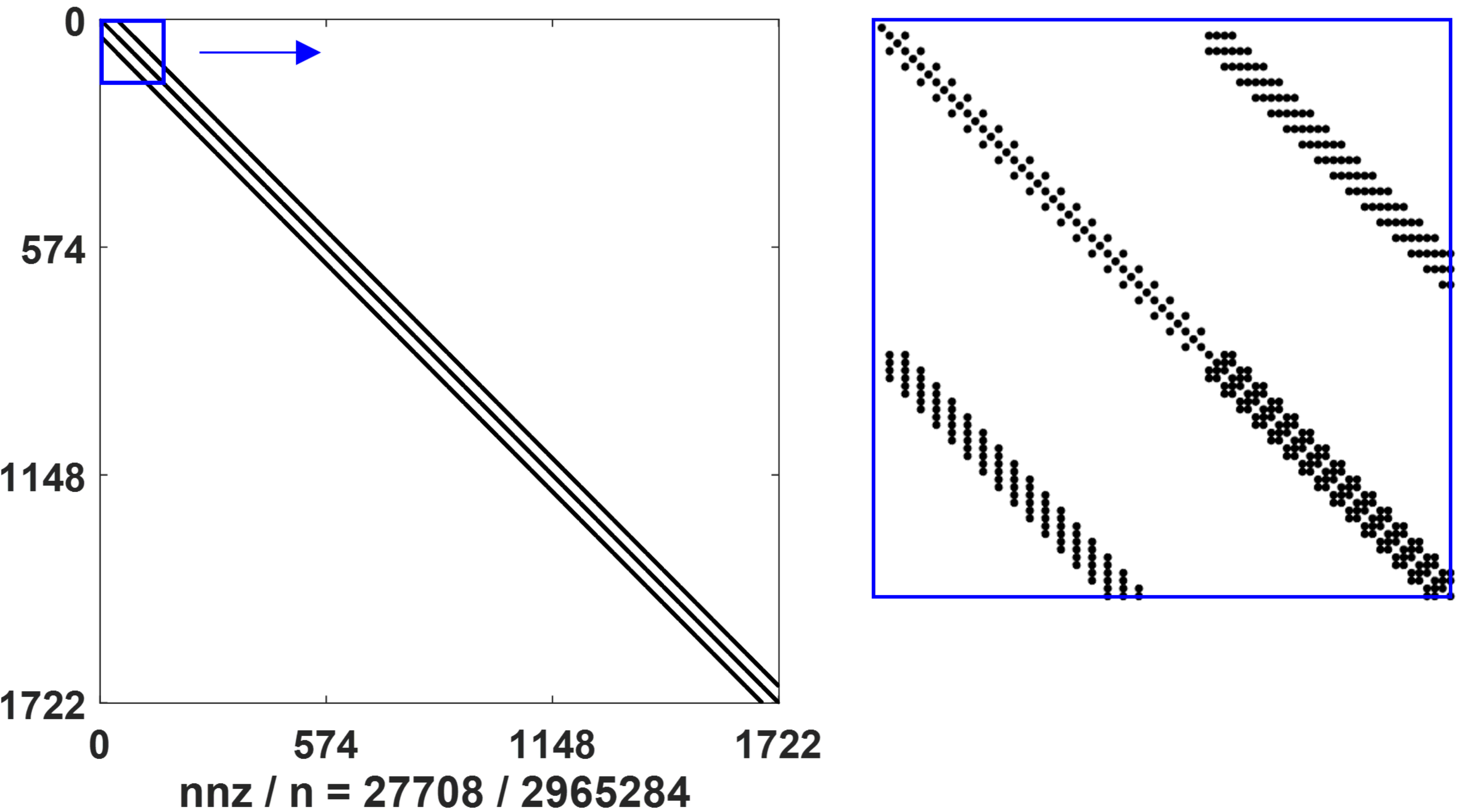}}
  \caption{Plots of the sparsities of the global stiffness matrices for different heterogenous structure analysis approaches: Homogenization, CMCM, FE$^{2}$, Substructuring and Our-CBN. Here, $nnz$ denotes the number of non-zero values, and $n$ is the size of the matrix. Note that $nnz$ of substructuring is almost $1.5$ times as that of the benchmark.}
  \label{fig:sparsity}
\end{figure}

\subsubsection{Complexity analysis of global displacement computation}
The complexity is determined by the DOFs in each coarse element. Homogenization, CMCM or FE$^2$ analysis only involves the corner nodes and has the node number of $8$ in 2D or $24$ in 3D; substructuring has the DOFs of $2b$ in 2D or $3b$ in 3D; our-CBN has the number of $6r$ in 2D or $q$ in 3D. The complexity analysis results are also summarized in Table~\ref{tab:complexity}.

For a more intuitive perspective, Fig.~\ref{fig:sparsity} further plots the sparsity of the global stiffness matrix of each method. Here, the coarse mesh has a dimension of $2\times 4$, and the local fine mesh has a dimension of $10\times 10$. Note that even for this simple example, substructuring has a dense stiffness matrix, and its number of nonzero elements is almost $1.5$ times that of the benchmark. Its direct use on large-scale industrial application problems may thus be impractical. Our CBN approach significantly reduces the number by introducing a B\'ezer interpolation matrix.
\section{Experiments}
The proposed approach of CBN heterogeneous structure analysis was implemented in MATLAB on an Intel Core i7, 3.7 GHz CPU and 64 GB RAM PC. Its performance was tested on various 2D and 3D examples. In all the examples, if not specifically stated, the matrices were assumed of Young's modulus $E = 1e^3$ and the inclusions of $E = 1$; both had a Poisson's ratio of $\mu=0.3$. Under these settings, the heterogeneous structures tended to present a large deformation that was more challenging to analyze with a high level of accuracy.

The substructuring approach is not further discussed as it constantly produces solutions of high accuracy with significantly high computational costs for large-size problems (see the complexity analysis in Section~\ref{sec:complexity}). We use \emph{Our-L} and \emph{Our-CBN} to denote our approach using linear interpolation or cubic B\'ezier interpolation; they have the same number of analysis DOFs for a fair comparison.

The analysis results on the global fine mesh were taken as the benchmark. In terms of the global energy or displacement, the analysis fidelity was measured via \emph{effectivity index} as the relative variation of the computed result with respect to the benchmark.
\begin{equation}
\begin{aligned}
    r_e = \frac{(e_1 - e_0)^2}{e_0^2}, \quad
    r_u = \frac{\int_\Omega (\bu_1 - \bu_0)^2}{\int_\Omega \bu_0^2},
  \end{aligned}
  \label{eq-effindex}
  \end{equation}
where $e_1,e_0$ are the energies of the computed and the benchmark, and $\bu_1,\bu_0$ are the displacements of the computed and the benchmark.

\begin{figure}[p]
  \centering
  \subfigure[Half MBB]{\includegraphics[width=0.4\textwidth]{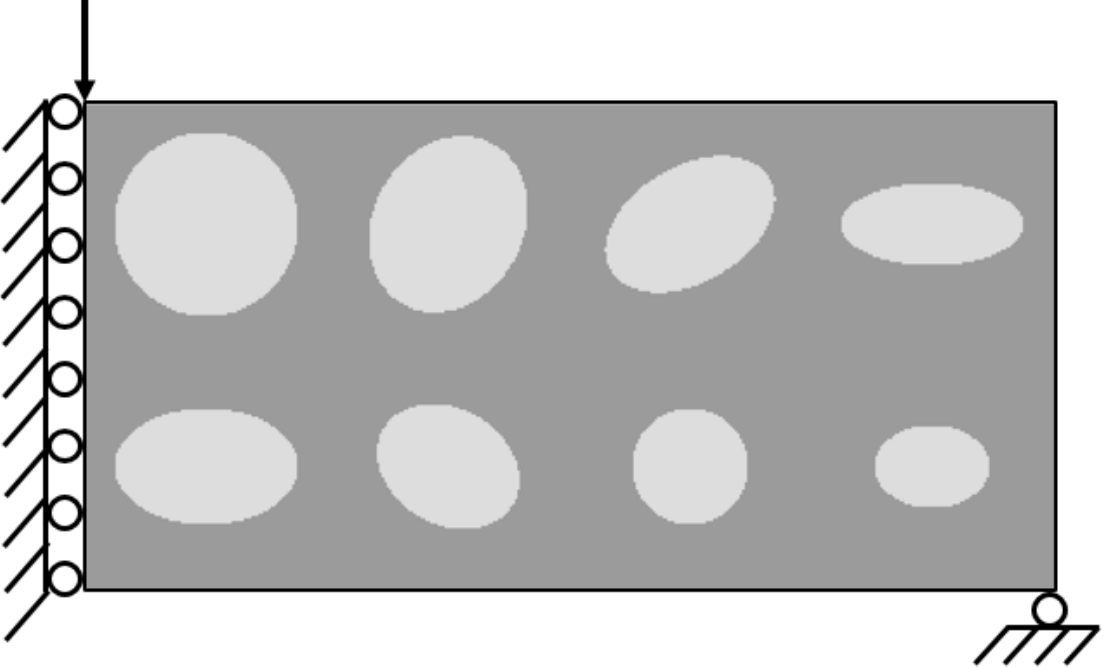}}
  \subfigure[Coarse mesh $\calM^H$ of $2 \times 4$]{\includegraphics[width=0.4\textwidth]{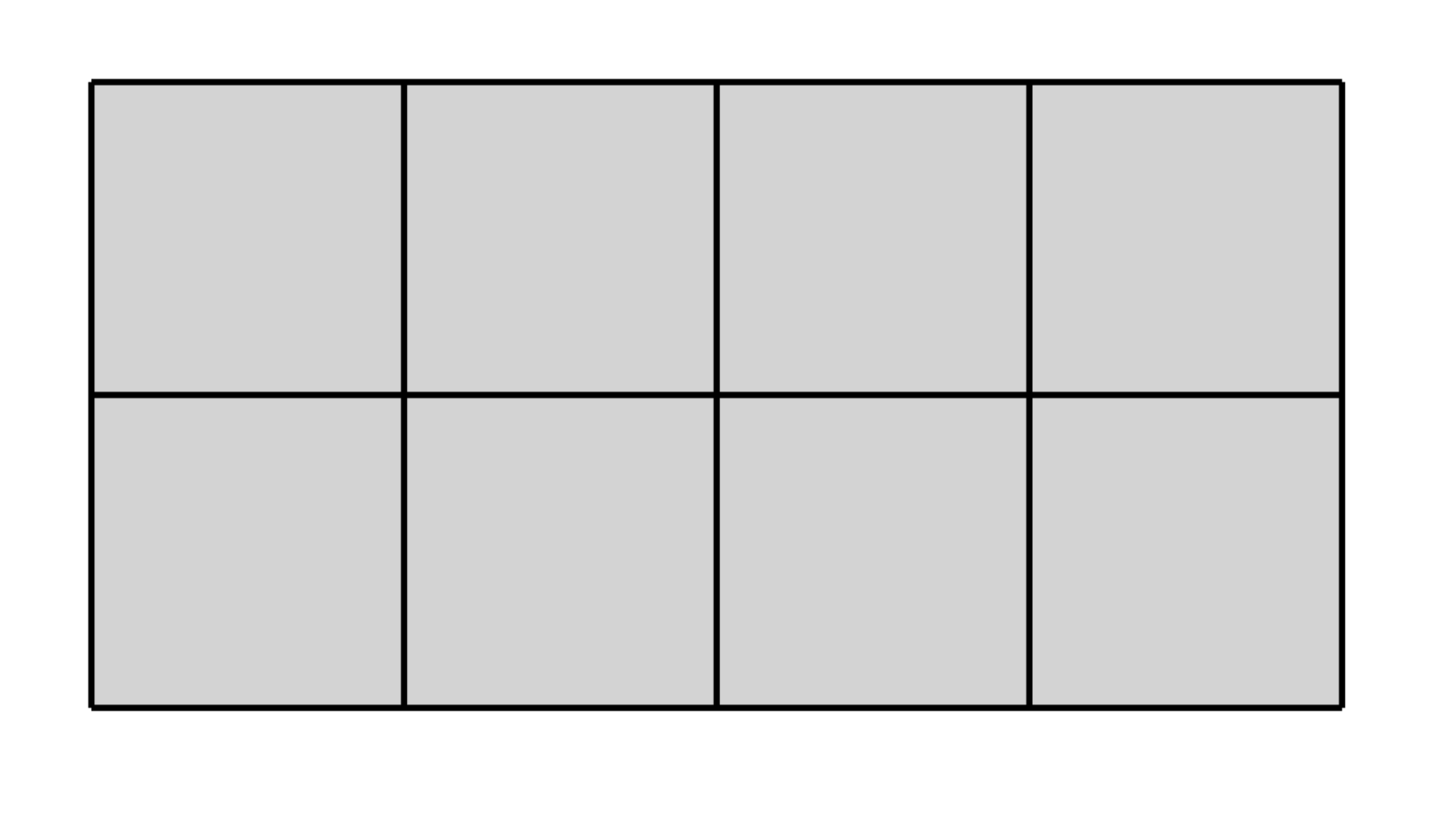}}
  \caption{The half MBB example has $2 \times 4$ coarse elements, each containing an elliptic inclusion.}
  \label{fig:MBB}
  \centering
  \subfigure[Benchmark]{\includegraphics[width=0.32\textwidth]{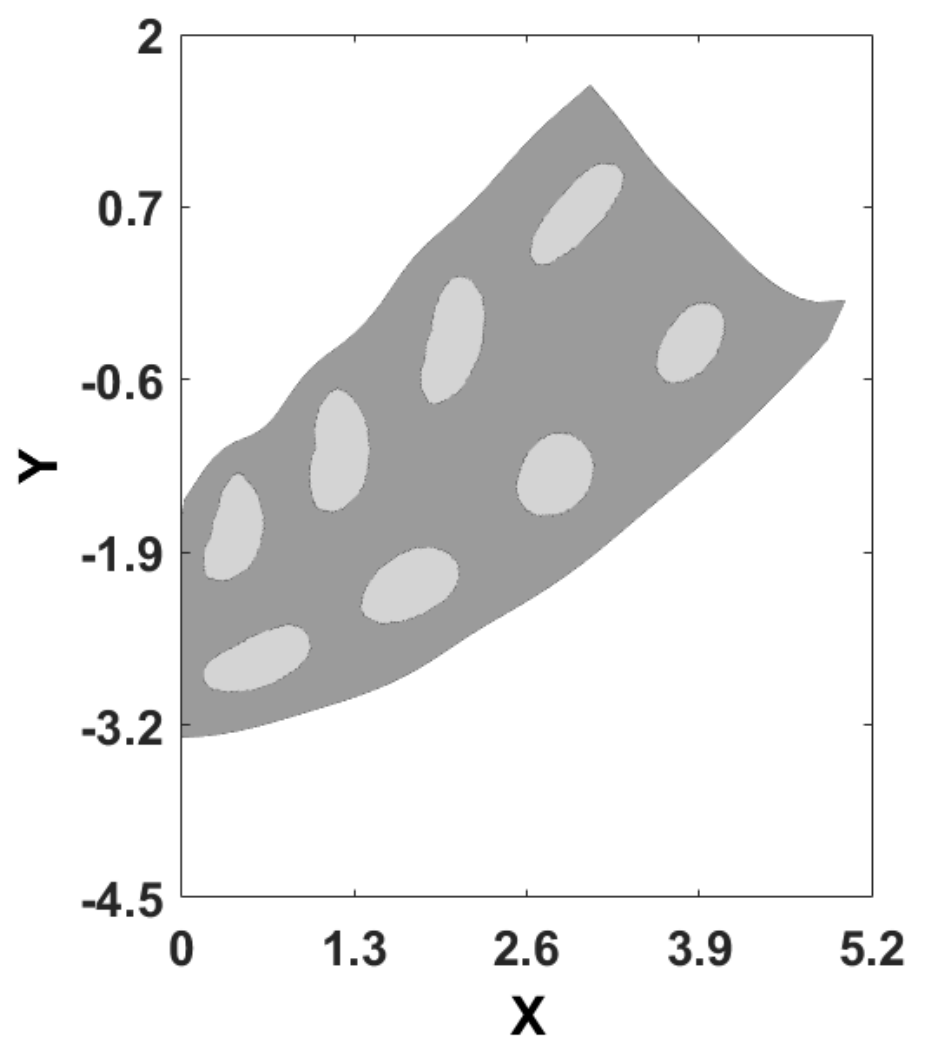}}
  \subfigure[Homogenization, $r_e = 0.47$, $r_u = 0.44$]{\includegraphics[width=0.32\textwidth]{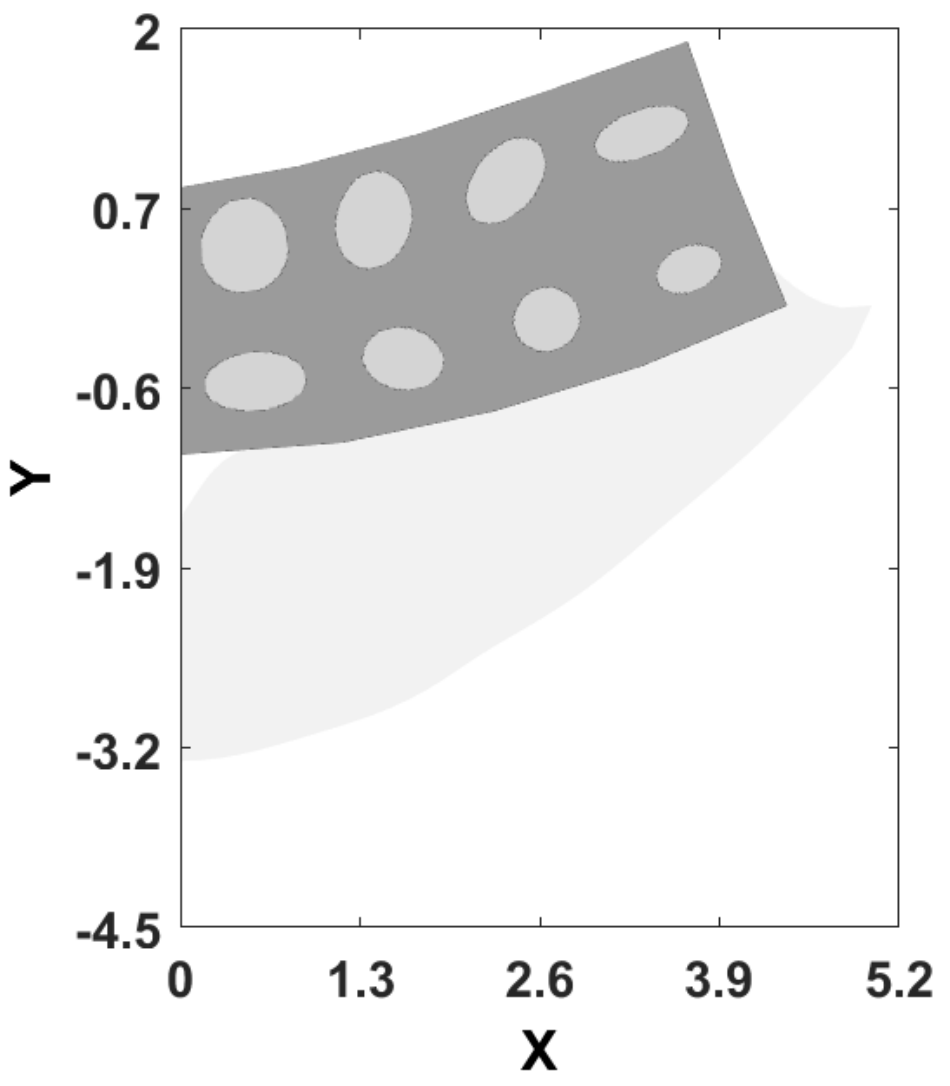}}
  \subfigure[FE$^2$, $r_e = 0.07, r_u = 0.06$]{\includegraphics[width=0.32\textwidth]{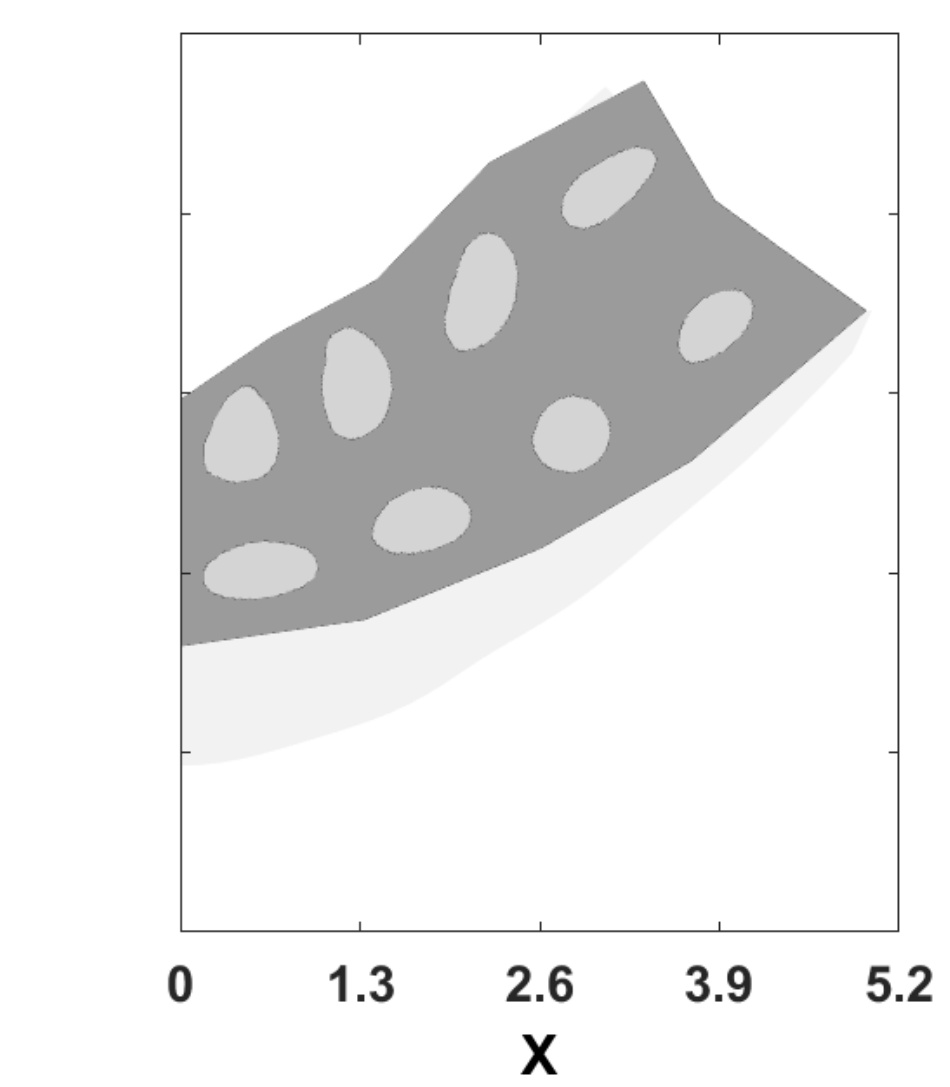}}
  \subfigure[CMCM, $r_e = 0.07, r_u = 0.08$]{\includegraphics[width=0.32\textwidth]{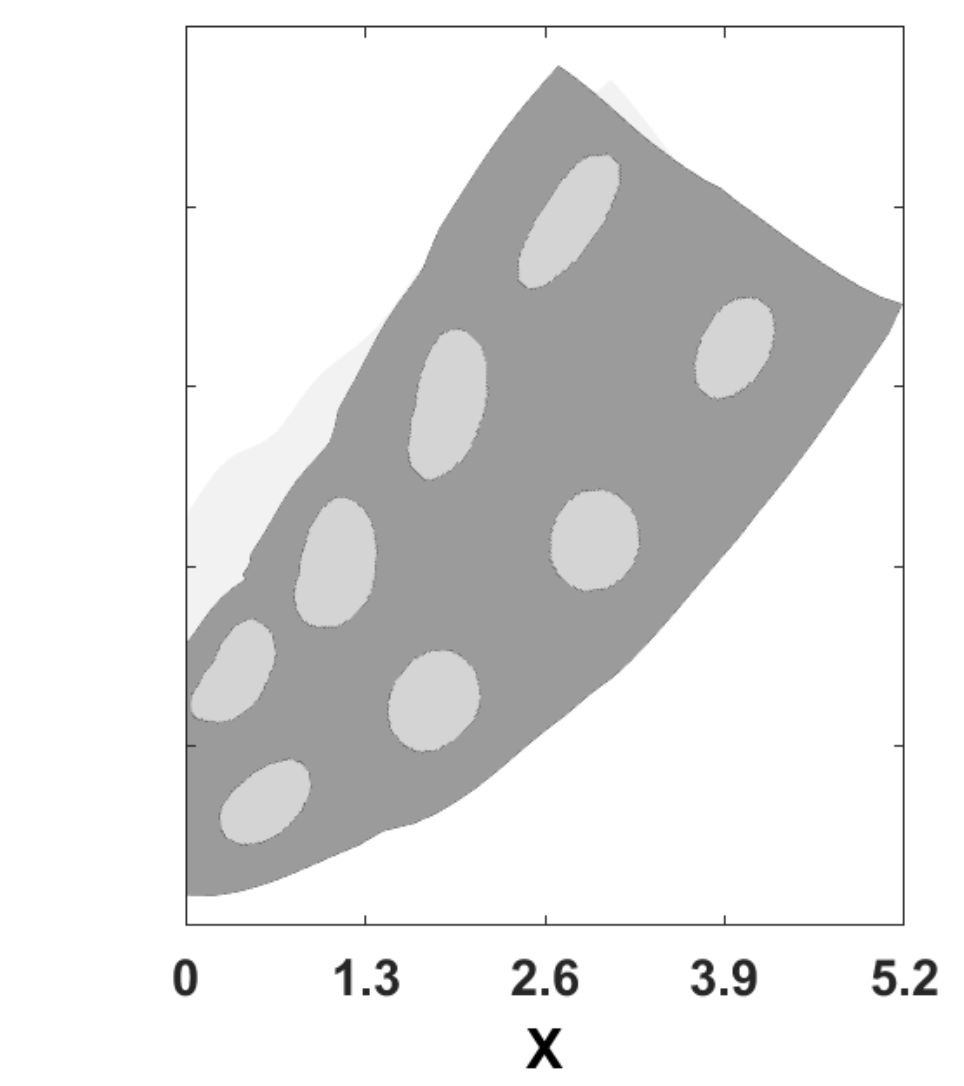}}
  \subfigure[Our-L, $r_e = 0.04, r_u = 0.04$]{\includegraphics[width=0.32\textwidth]{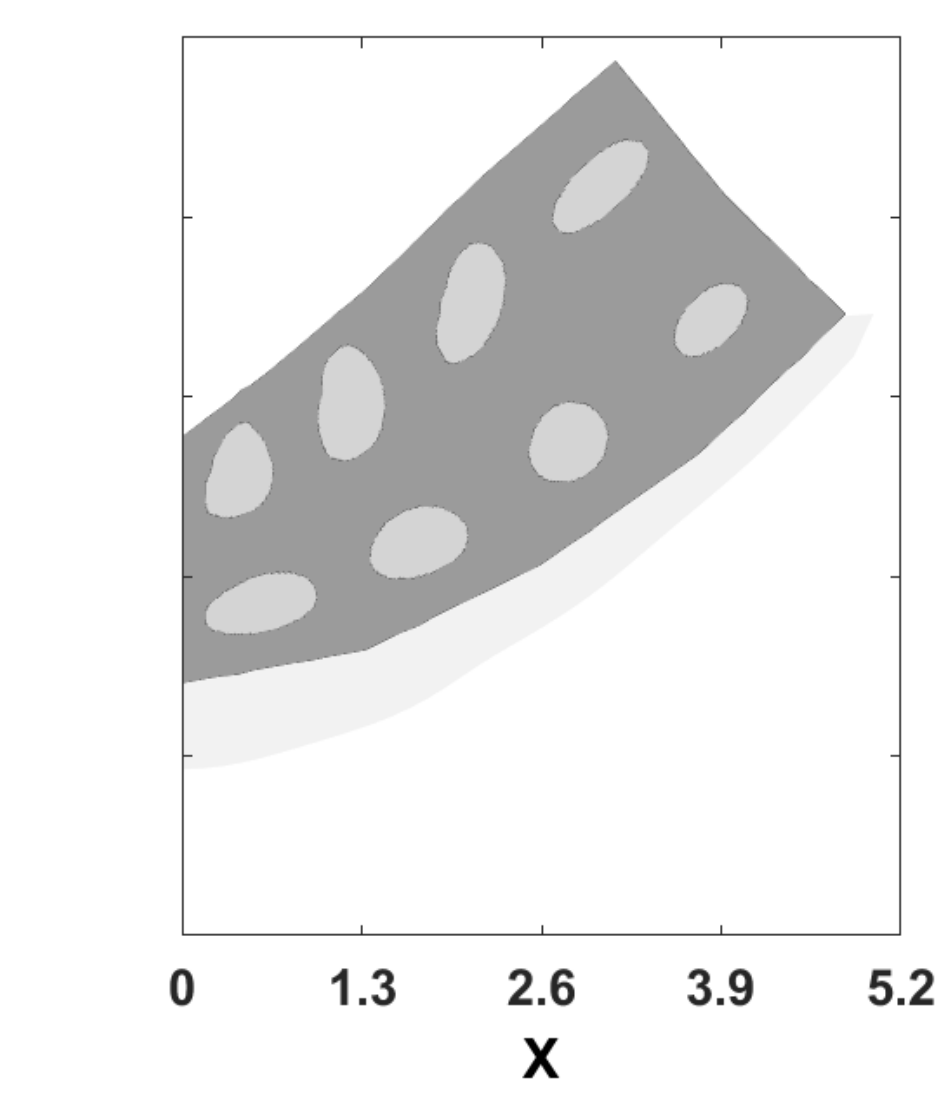}}
  \subfigure[Our-CBN, $r_e = 7.9e^{-4}, r_u = 9.1e^{-4}$]{\includegraphics[width=0.32\textwidth]{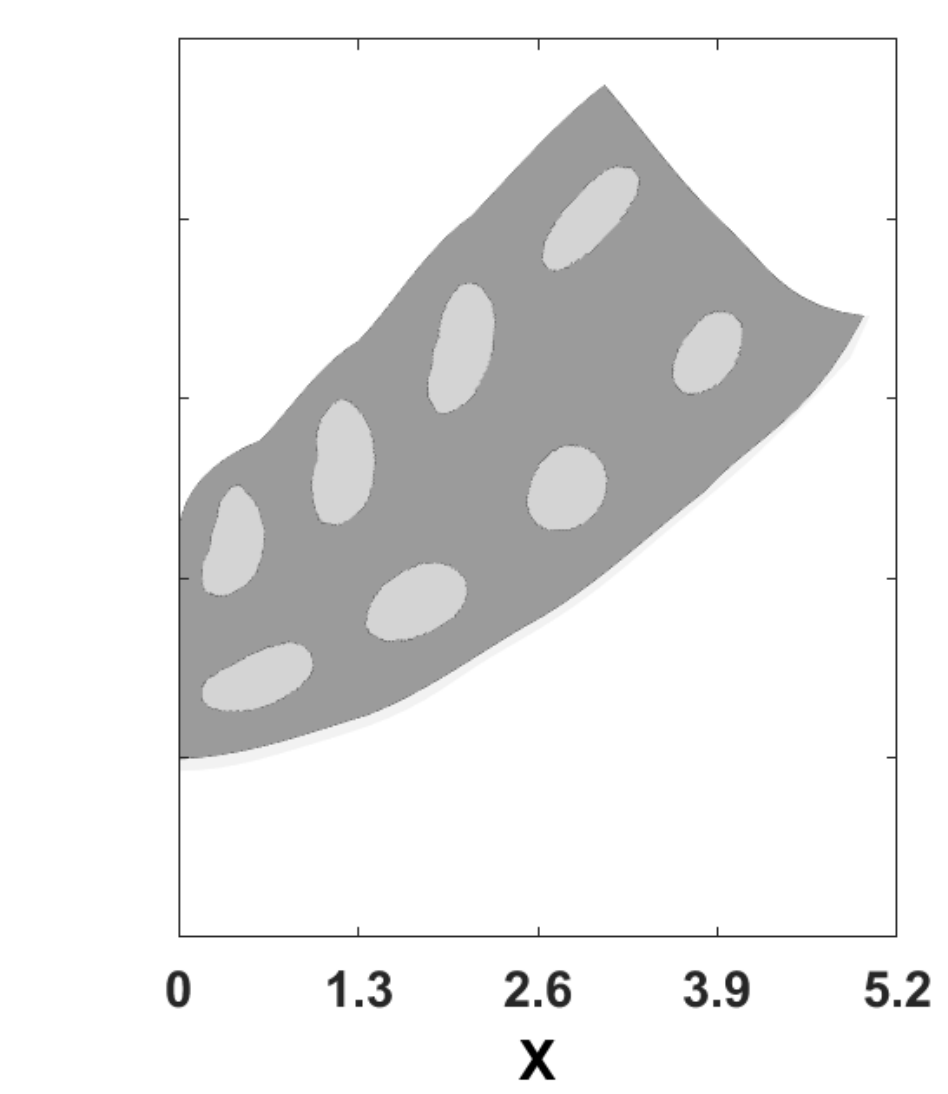}}\\
  \caption{Numerical results of the heterogeneous half MBB example in Fig.~\ref{fig:MBB} compared with benchmark and other related approaches: Homogenization, FE$^2$, CMCM and Our-L, where the shadow areas denote the benchmark deformations.}
  \label{fig:comparison}
\end{figure}

\subsection{Overall performance and comparisons with related approaches}\label{sec:overall}
We first show the overall performance and its comparison with related approaches using the half heterogeneous MBB example in Fig.~\ref{fig:MBB}. The coarse mesh is of size $2 \times 4$, and the local fine mesh is of size $10 \times 10$. The results are summarized in Fig.~\ref{fig:comparison} and Table~\ref{tab:time}.

The reconstructed deformation for each approach is shown in Fig.~\ref{fig:comparison} with the effective indices $r_e$ and $r_u$. Large deformation differences are clearly observed between results of the benchmark and homogenization, FE$^2$, CMCM, and Our-L. Instead, Our-CBN has a deformation almost identical to that of the benchmark, even at the local regions of large deformation. Their effectivity indices indicate similar phenomenon: homogenization has the largest index of 0.47, Our-CBN has the smallest of $7.9e^{-4}$, and FE$^2$, CMCM, Our-L have an index of approximately 0.07. A two-order improvement is observed using Our-CBN. Note that interestingly, homogenization, FE$^2$, and Our-L tend to over-stiffen the deformation while CMCM tends to soften it in this example.

The time costs are also summarized in Table.~\ref{tab:time}. As indicated, the benchmark takes the longest time and all the other approaches reduce it dramatically. In the local analysis for a coarse element, different approaches have similar timings, although Our-CBN and Our-L take slightly more time; the local computations can be conducted in parallel online except for FE$^2$. In online computation of the global coarse displacement, FE$^2$ takes much more time than the other approaches as it requires 13 iterations. Our-L and Our-CBN take more time than homogenization and CMCM, and this difference may further increase for super large-sized analysis problems. The numerical results are consistent with the algorithmic complexity analysis in Section ~\ref{sec:complexity}.

\begin{table}[t]
  \centering
  \caption{Approach timings (in s): Homogenization, FE$^2$, CMCM, Our-L, and Our-CBN for the half MBB example in Fig.~\ref{fig:MBB}.}
  \label{tab:time}
  \begin{threeparttable}
    \begin{tabular}{c|c|c}
      \hline
      \rowcolor[HTML]{C0C0C0}
      {\color[HTML]{333333} Approaches} & \begin{tabular}[c]{@{}c@{}} Solving linear equations systems  \\ for each coarse element\end{tabular} & \begin{tabular}[c]{@{}c@{}}Solving $\bK \bU=\bF$  \\ in online analysis\end{tabular}\\ \hline
      \rowcolor[HTML]{FFFFFF}
      Benchmark                      & -                                     & 0.21             \\
      \rowcolor[HTML]{EFEFEF}
      Homogenization                 & 3.2$e^{-4}$                           & 4$e^{-3}$        \\
      \rowcolor[HTML]{FFFFFF}
      FE$^2$                         & 3.2$e^{-4}$                           & 4$e^{-3} \times$ 13 \tnote{1}        \\
      \rowcolor[HTML]{EFEFEF}
      CMCM                           & 3.3$e^{-4}$                           & 4$e^{-3}$        \\
      \rowcolor[HTML]{FFFFFF}
      Our-L                     & 4.3$e^{-4}$                           & 5$e^{-3}$        \\
      \rowcolor[HTML]{EFEFEF}
      Our-CBN                            & 4.3$e^{-4}$                           & 5$e^{-3}$        \\ \hline
    \end{tabular}
    \begin{tablenotes}
      \footnotesize
      \item[1] FE$^2$ has 13 iterations in this example.
    \end{tablenotes}
  \end{threeparttable}
\end{table}

\subsection{Performance at different mesh settings}\label{sec:meshset}
\begin{figure}[tbh]
  \centering
\subfigure[Different numbers of coarse nodes]{\includegraphics[width=0.45\textwidth]{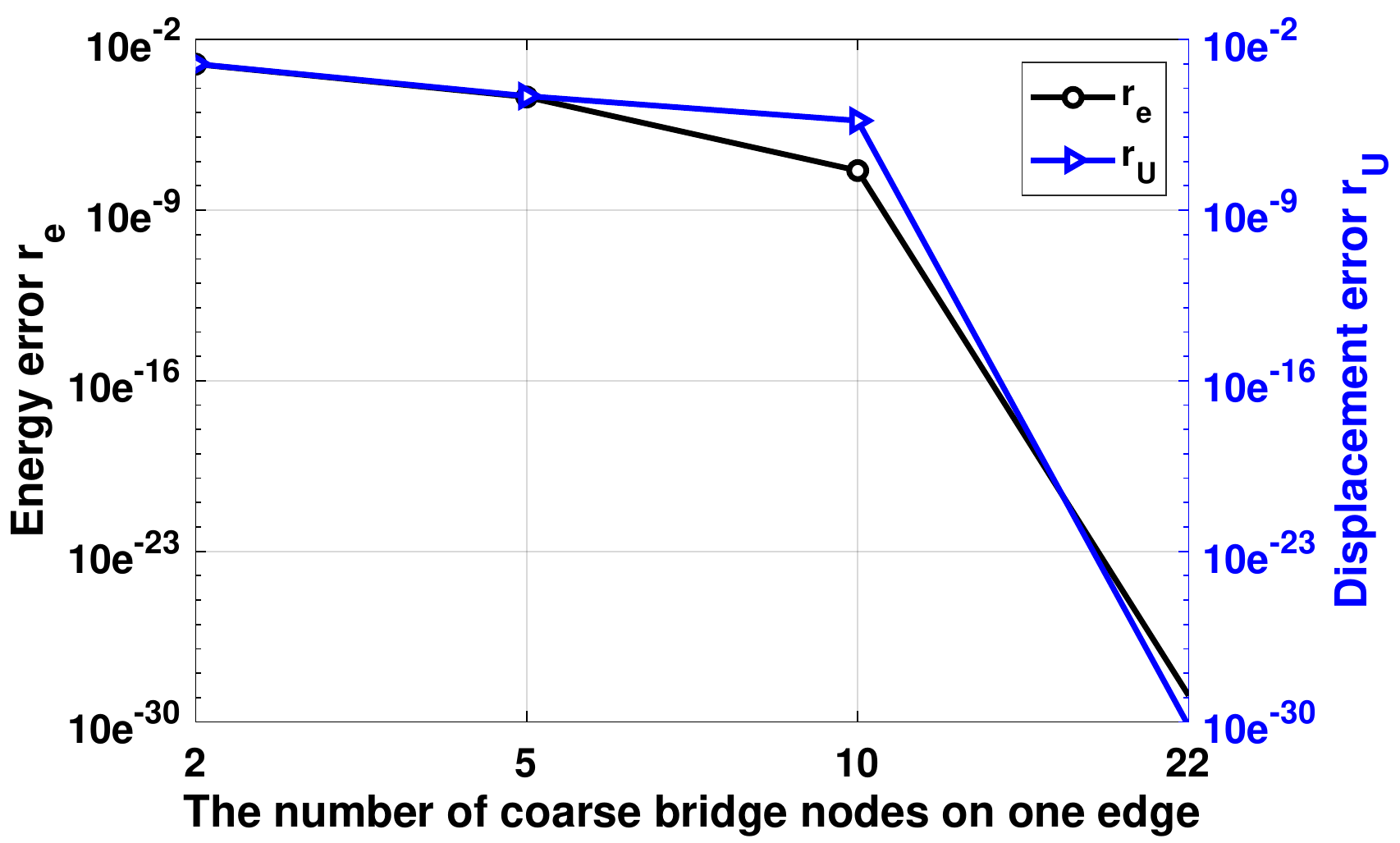}}  
\subfigure[Different sizes of coarse elements]{ \includegraphics[width=0.45\textwidth]{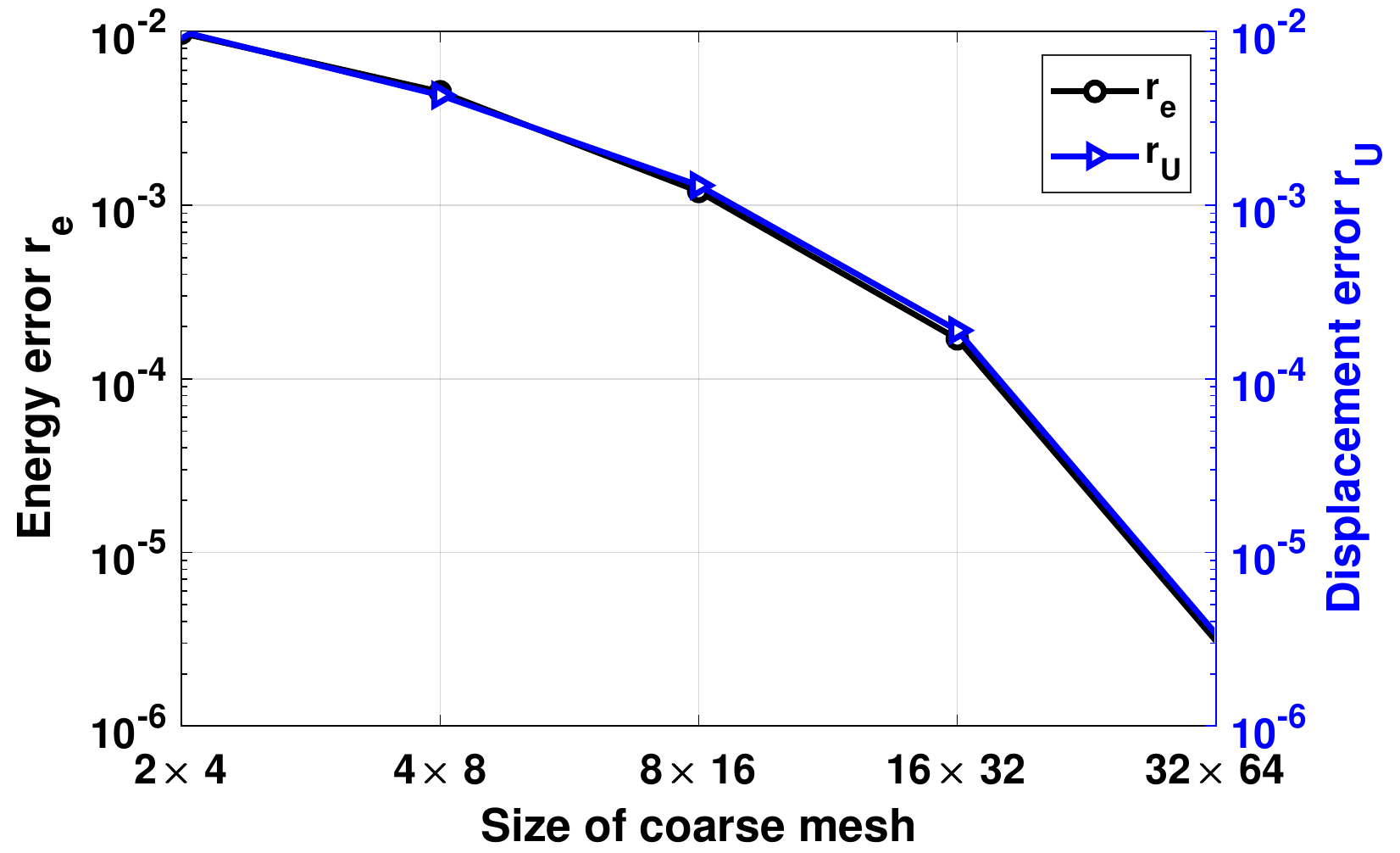}}  
\subfigure[Different material contrasts]{ \includegraphics[width=0.5\textwidth]{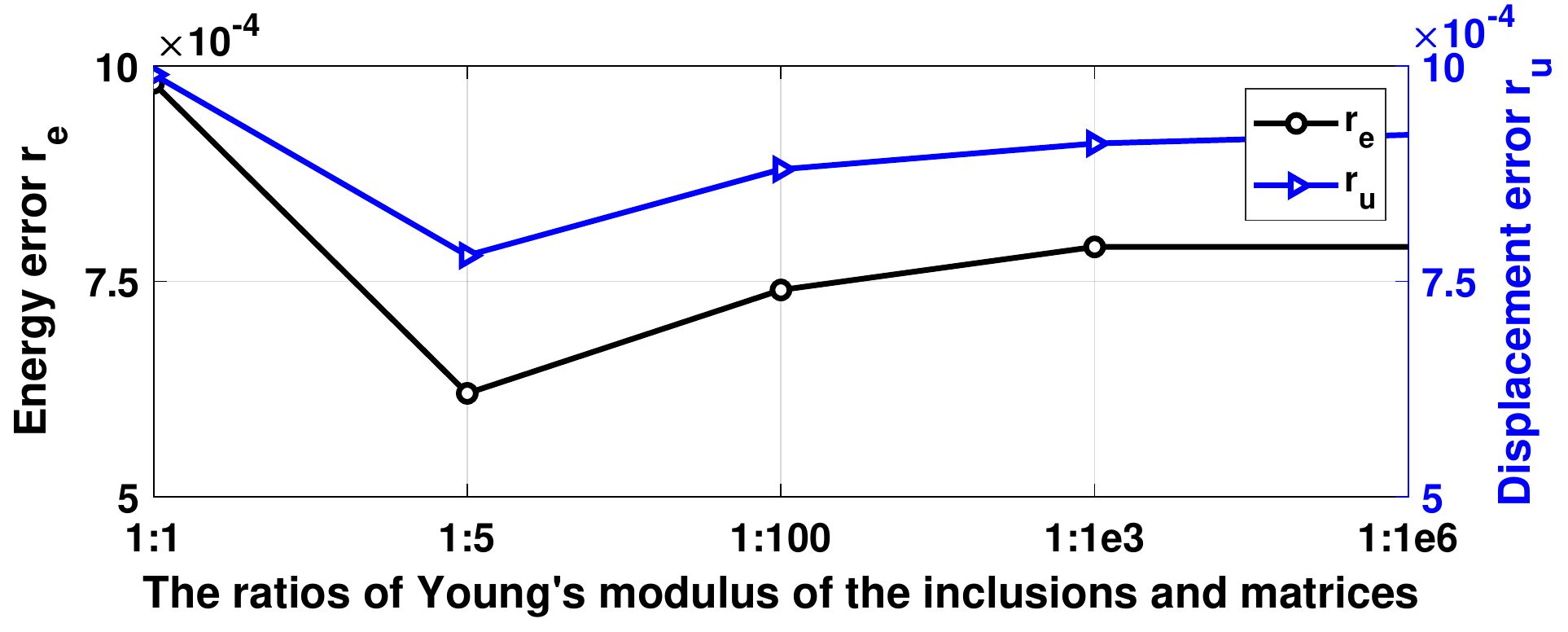}}   \label{fig:curve_mat}
  \caption{Variations of effectivity indices $r_e$, $r_u$ under different mesh settings.}
  \label{fig:parameters}
\end{figure}

The performance of Our-CBN is further tested at different analysis parameters: size of coarse mesh and number of bridge nodes or contrast of material stiffnesses. The half MBB in Fig.~\ref{fig:MBB} is used here, where the coarse mesh is of size $2\times 4$, and the local fine mesh is of size $64\times 64$. The results are presented in Fig.~\ref{fig:parameters}.

\paragraph{Different numbers of bridge nodes $\calV_r$} Our-CBN has the unique ability of choosing different numbers of bridge nodes. Its performance is tested at bridge node numbers of $2,5,10,22$ along a boundary, and Fig.~\ref{fig:parameters}(a) plots their effectivity indices $r_e$ and $r_u$. The indices decrease rapidly as the bridge node number increases, and have in particular a very high accuracy of $r_e=1.2e^{-28}$ and $r_u=6.7e^{-30}$ in the case of $22$ bridge nodes. In this case, the CBN number is equal to the number of all boundary nodes, and Our-CBN essentially conducts an identical analysis to the benchmark on the global fine mesh.

\begin{figure}[t]
  \centering
  \subfigure[$\calM^H$ of $4 \times 8$]{\includegraphics[width=0.23\textwidth]{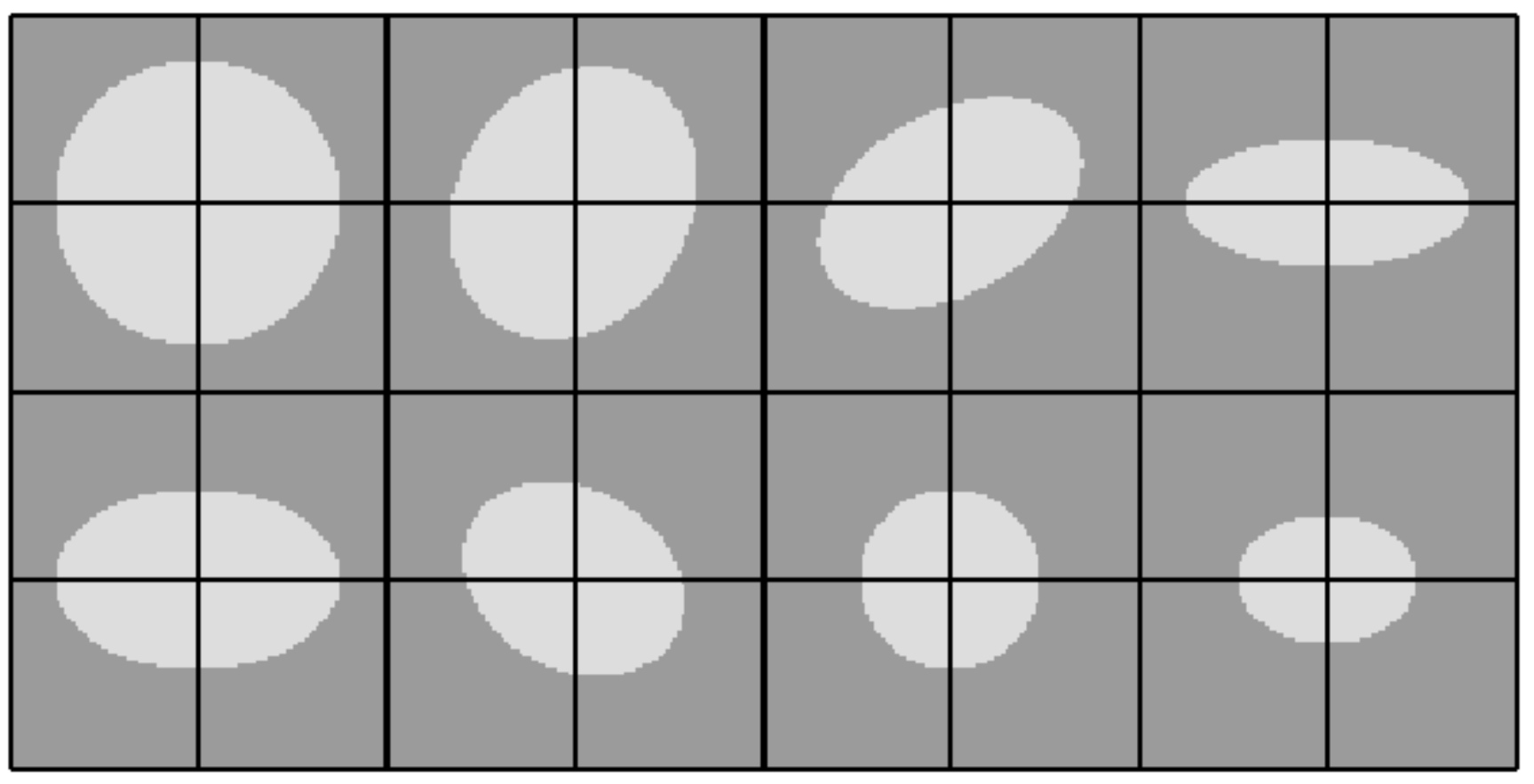}}\quad
  \subfigure[$\calM^H$ of $8 \times 16$]{\includegraphics[width=0.23\textwidth]{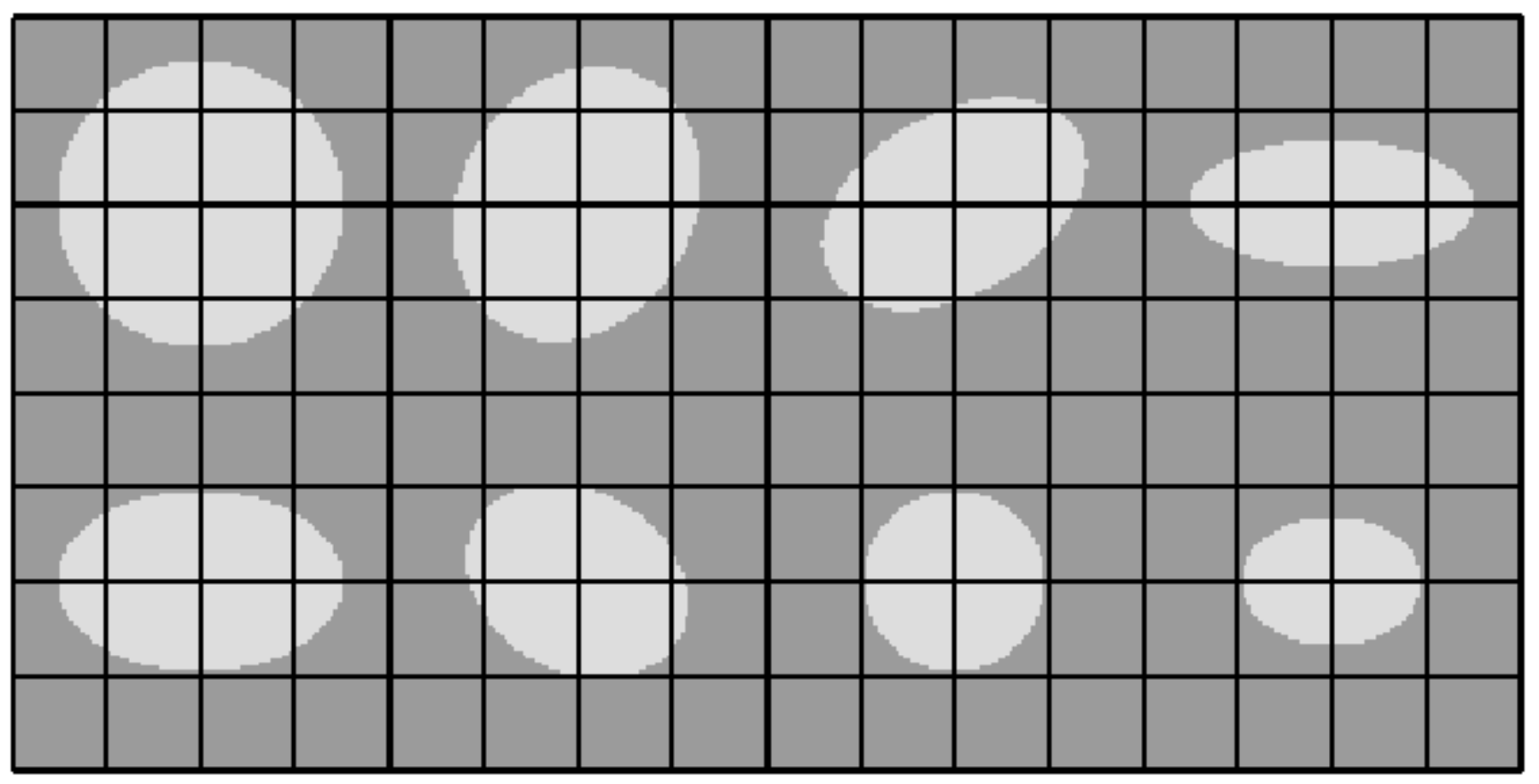}}\quad
  \subfigure[$\calM^H$ of $16 \times 32$]{\includegraphics[width=0.23\textwidth]{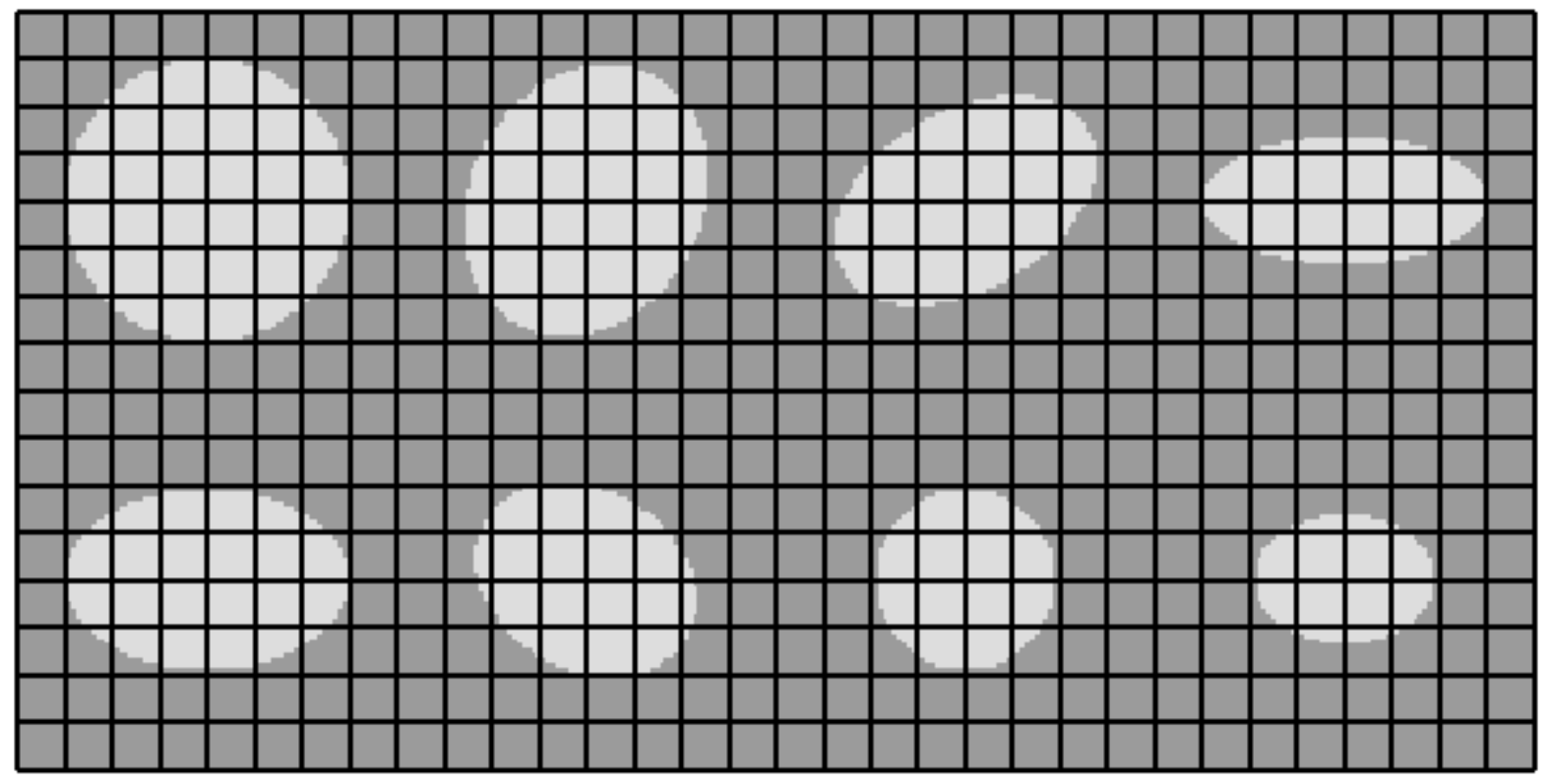}}\quad
  \subfigure[$\calM^H$ of $32 \times 64$]{\includegraphics[width=0.23\textwidth]{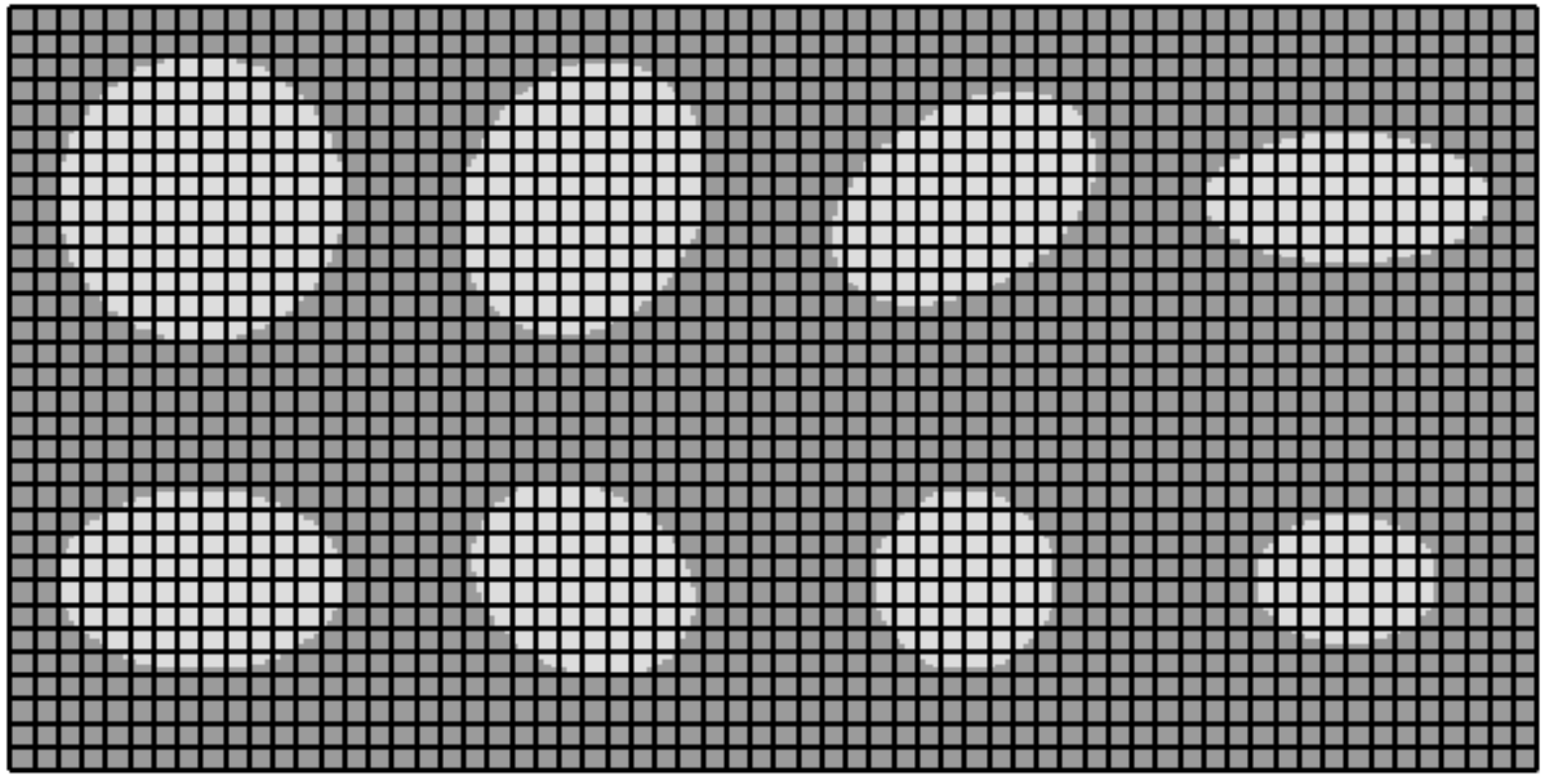}}
  \caption{Different coarse mesh sizes are set for the half MBB in Fig.~\ref{fig:MBB}.}
  \label{fig:MBB_size}
 \end{figure}

\paragraph{Different sizes of coarse meshes $\calM^H$} As indicated in Fig.~\ref{fig:MBB_size}, five different sizes of coarse mesh $\calM^H$ $2\times 4$, $4\times 8$, $8\times 16$, $16\times  32$, and $32\times 64$ are set, while the global fine mesh size is kept unchanged as $128\times 256$. Note that excluding the case of $2\times 4$, the coarse mesh boundaries will cross the different interface materials. The dramatic material variations along the boundaries pose a significant challenge in terms of maintaining analysis accuracy. Our approach still demonstrates its high performance as indicated by the effectivity indices in Fig.~\ref{fig:parameters}(b).

\paragraph{Different contrasts of material stiffness}
Different contrasts of Young's modulus are respectively set for the inclusions and matrices: $(1,1),~ (1,5),~ (1,100),~(1,1e^3),~(1,1e^6)$. Generally, the larger the ratios, the more difficult it is to achieve a reasonable result, considering the fact that the relatively softer inclusion tensors tend to show a large local deformation. Still, Our-CBN is able to maintain a high analysis accuracy, of all effectivity indices below $1e^{-3}$, even for the extreme ratio of $1:1e^6$ (Fig.~\ref{fig:parameters}(c)).

\subsection{Shape functions in terms of material distributions}\label{sec:sf}

\begin{figure}[t]
  \centering
  \includegraphics[width=1\textwidth]{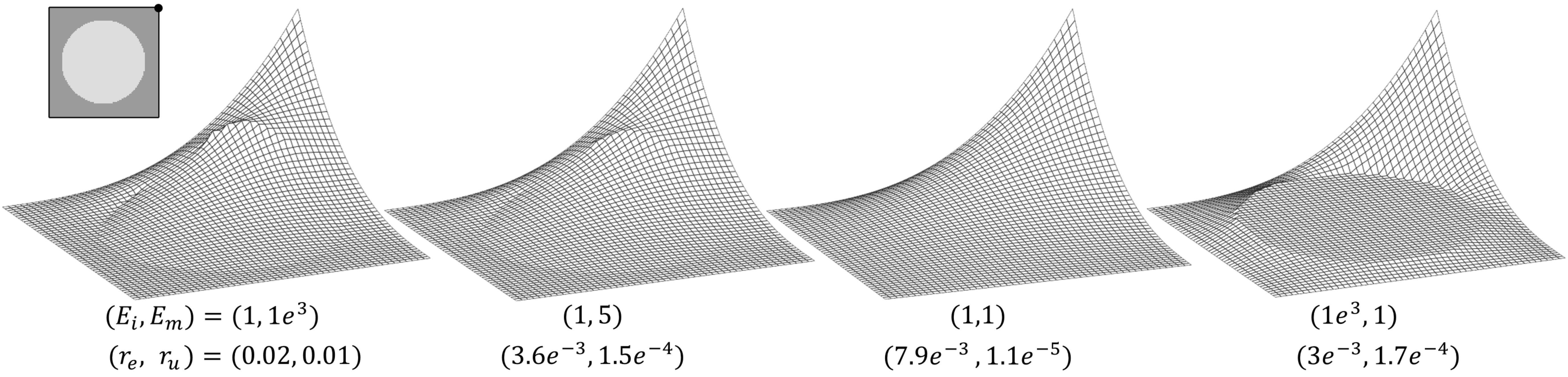}
  \caption{Surfaces of the first shape function component $\bN_{11}(\bx)$ and the effectivity indices for different material contrasts.}\label{fig:NH-E}
 \centering
  \includegraphics[width=1\textwidth]{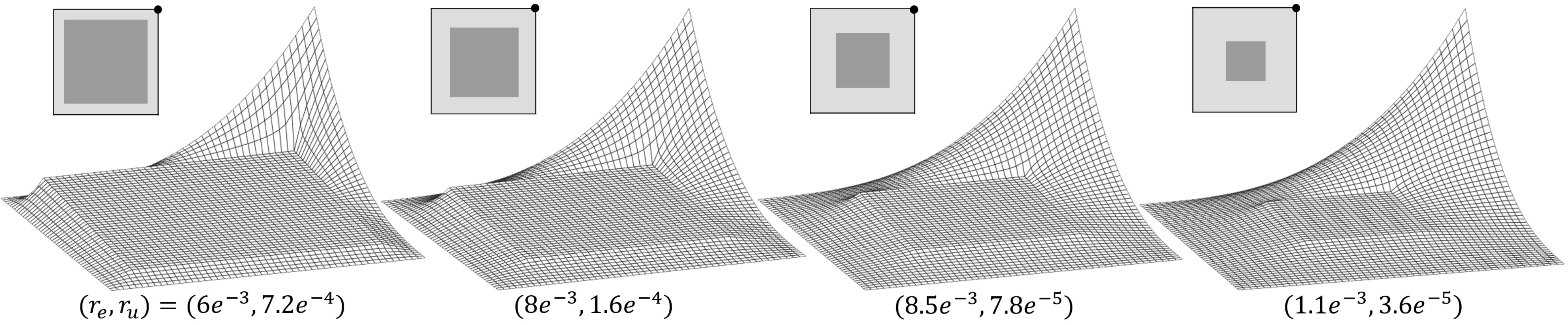}
  \caption{Surfaces of the first shape function components $\bN_{11}(\bx)$ and the effectivity indices for different-sized squared inclusions.}
  \label{fig:NH-shape}
\end{figure}

The constructed CBN shape functions are expected to closely reflect the interior material distributions for high accuracy analysis, irrespective of the imposed boundary conditions. We test this expectation by plotting in Fig.~\ref{fig:NH-E} the surfaces of the CBN shape functions on cases of different material distributions. The local fine mesh has a size of $50\times 50$, and we take the corner nodes as bridge nodes. Their associated effectivity indices on a $2\times2$ coarse mesh are also shown below each shape function.

\paragraph{Shape functions at different material contrasts} We plot in Fig.~\ref{fig:NH-E} surfaces of the first shape function component $\bN_{11}(\bx)$ at different contrasts of Young's modulus $(1,1e^3)$, $(1,5)$, $(1,1)$ and $(1e^3, 1)$ for inclusions and the matrices. The surface presents a smooth variation in the case of $(1,1)$. In contrast, the surface drops rapidly over the softer inclusion in case (a) while it remains almost unchanged over the stiffer inclusions in (d). The results demonstrate our CBN shape functions' level of adaptation to the variations of the material stiffness.

\paragraph{Shape functions at different sizes of inclusions}
We further plot in Fig.~\ref{fig:NH-shape} surfaces of the first component $\bN_{11}(\bx)$ on the case of different-sized squared inclusions, where the darker and lighter regions respectively have a Young's modulus of $E = 1e^3$ and $E = 1$. The surfaces of the shape functions present clearly flatter variations right above the squared stiffer area, which is consistent with our expectation.

\begin{figure}[htb]
  \centering
  \subfigure[Structure of a bending problem]{\includegraphics[width=1\textwidth]{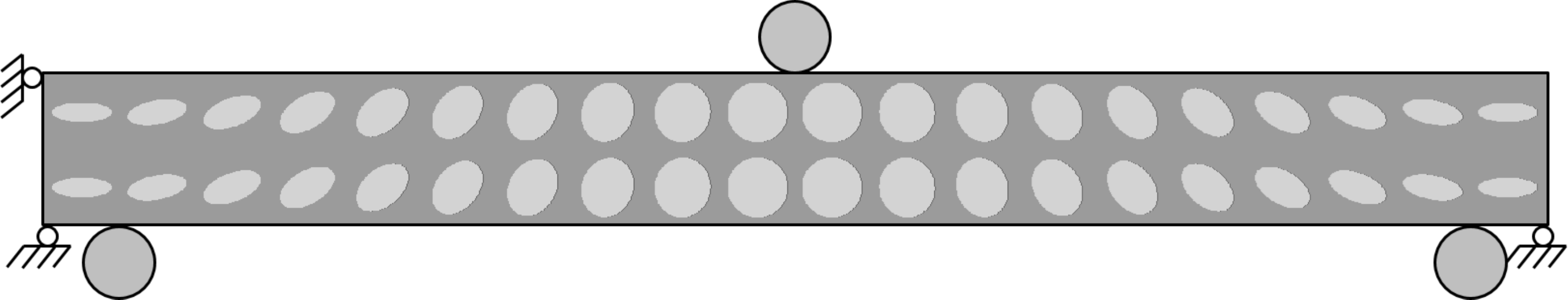}}
  \subfigure[$\calM^H$ of $1\times 10$]{\includegraphics[width=1\textwidth]{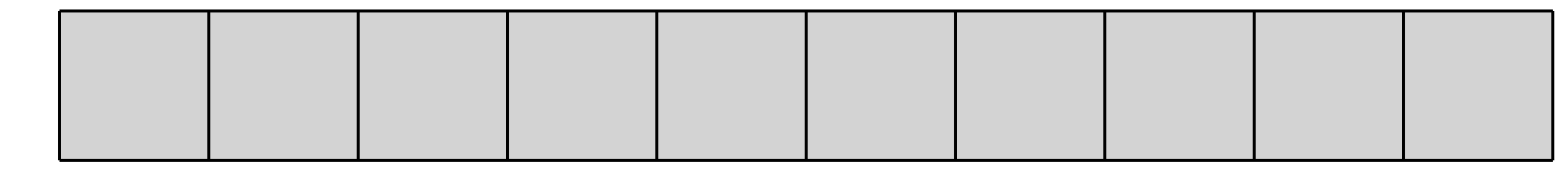}}
  \subfigure[$\calM^H$ of $2\times 20$]{\includegraphics[width=1\textwidth]{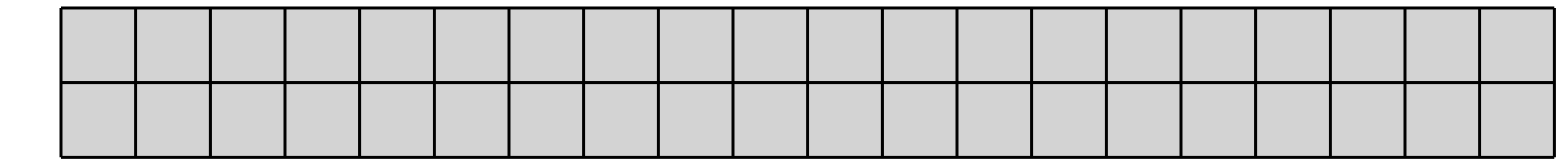}}
  \subfigure[$\calM^H$ of $4\times 40$]{\includegraphics[width=1\textwidth]{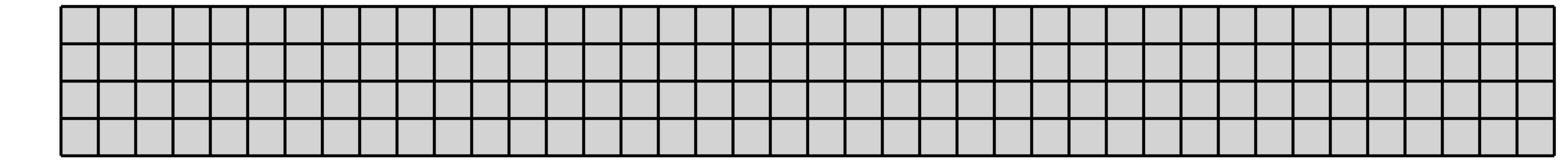}}
  \subfigure[$\calM^H$ of $8\times 80$]{\includegraphics[width=1\textwidth]{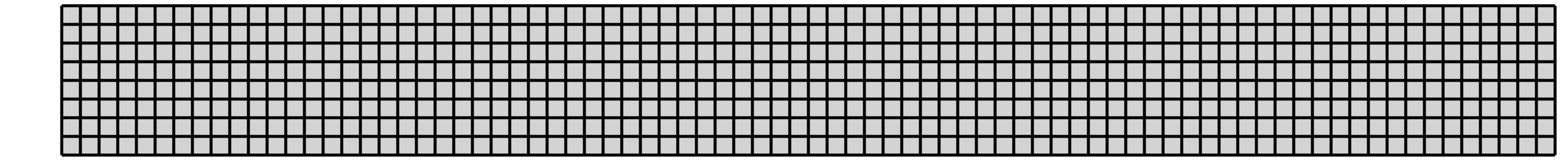}}
  \caption{Structure of a heterogeneous bending problem, and the different-sized coarse meshes $\calM^H$.}
  \label{fig:domain_bending}
\end{figure}

\subsection{A 2D heterogeneous bending beam}\label{sec:bending}
\begin{figure}[tb]
  \centering
  \includegraphics[width=0.6\textwidth]{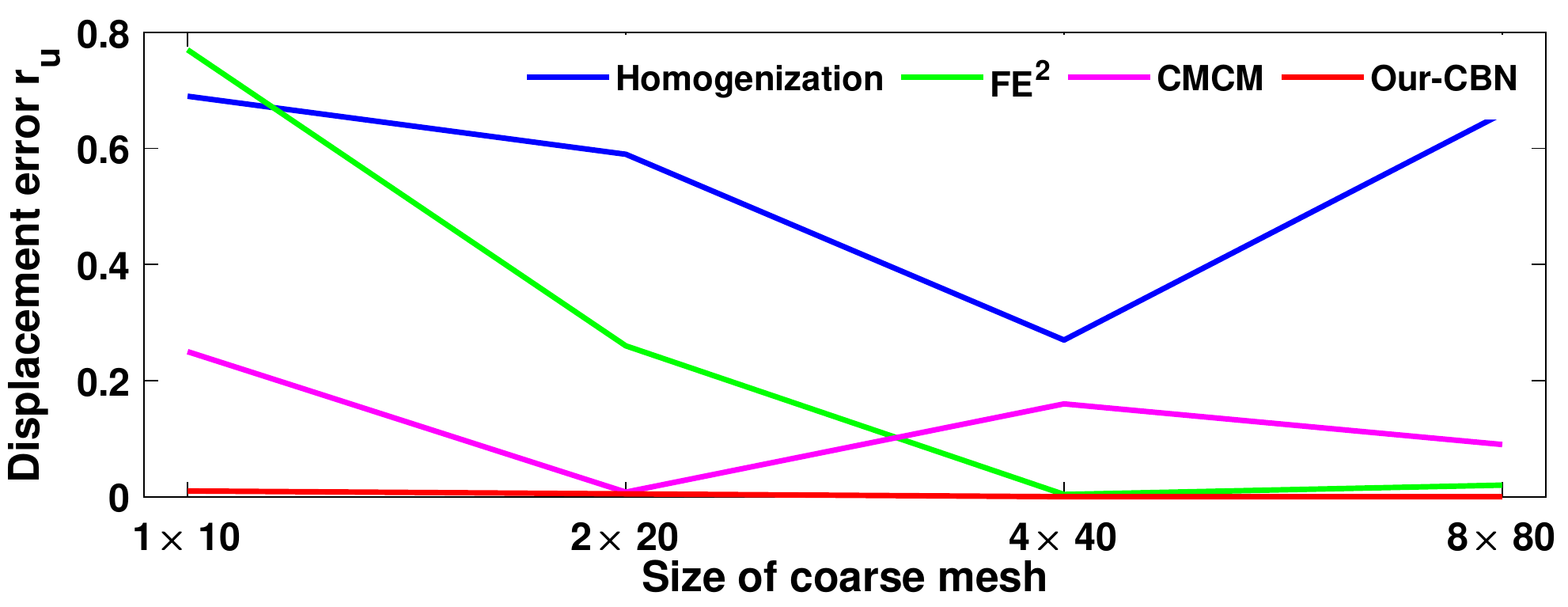}
  \caption{Variations of effectivity indices $r_u$ under different-sized coarse meshes: $1\times 10$, $2\times 20$, $4\times 40$, and $8\times 80$.}
  \label{fig:curve_bending}
\end{figure}

\begin{table}[t]
  \centering
  \caption{Effectivity indices $r_e$, $r_u$ at different-sized coarse meshes $\calM^H$, compared with homogenization, FE$^2$, CMCM, and Our-CBN.}
  \label{tab:macro_mesh}
  \begin{threeparttable}
    \begin{tabular}{c|c|c|c|c}
      \hline
      \rowcolor[HTML]{C0C0C0}
      \cellcolor[HTML]{C0C0C0}Size of $\calM^H$     & \begin{tabular}[c]{@{}c@{}}$1 \times 10$ \\ ($r_e$ / $r_u$)\end{tabular} & $2 \times 20$         & \begin{tabular}[c]{@{}c@{}}$4 \times 40$\\ (cutting interfaces)\end{tabular} & \multicolumn{1}{c|}{\cellcolor[HTML]{C0C0C0}\begin{tabular}[c]{@{}c@{}} $8 \times 80$\\ (cutting interfaces)\end{tabular}} \\ \hline
      \rowcolor[HTML]{EFEFEF}
      \cellcolor[HTML]{EFEFEF}Homogenization & 0.70 / 0.69                & 0.61 / 0.59           & 0.31 / 0.27                & 0.72 / 0.66                                                             \\
      FE$^2$                                 & 1.00 / 0.77 \tnote{1}      & 0.29 / 0.26           & 0.01 / 4e$^{-3}$           & 0.05 / 0.04                                                             \\
      \rowcolor[HTML]{EFEFEF}
      \cellcolor[HTML]{EFEFEF}CMCM           & 0.23 / 0.25                & 0.01 / 8e$^{-3}$      & 0.11 / 0.16                & 0.09 / 0.09                                                             \\
      \rowcolor[HTML]{EFEFEF}
      \cellcolor[HTML]{EFEFEF}Our-CBN            & 0.02 / 0.01                & 6e$^{-3}$ / 5e$^{-3}$ & 1e$^{-4}$ / 6e$^{-5}$      & 5e$^{-6}$ / 2e$^{-6}$                                                   \\ \hline
    \end{tabular}
    \begin{tablenotes}
      \footnotesize
      \item[1] FE$^2$ fails to converge in this context.
    \end{tablenotes}
  \end{threeparttable}
\end{table}

\begin{figure}[p]
  \centering
  \subfigure[Benchmark, $\varepsilon_{11}$]{\includegraphics[width=0.75\textwidth]{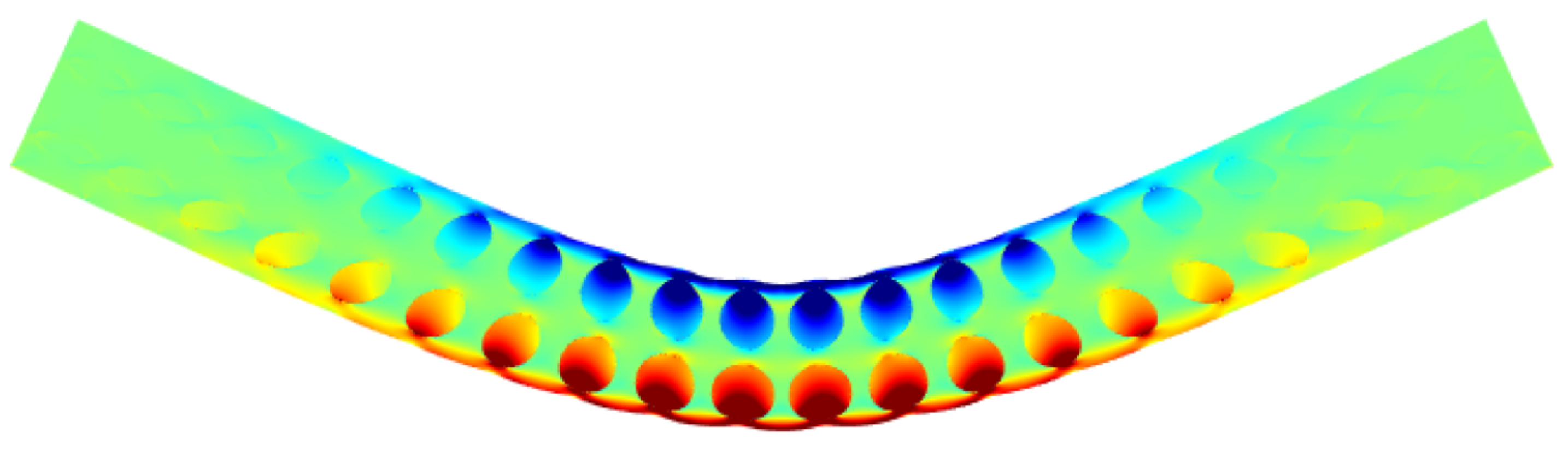}}
  \subfigure[Our-CBN, $\varepsilon_{11}$]{\includegraphics[width=0.75\textwidth]{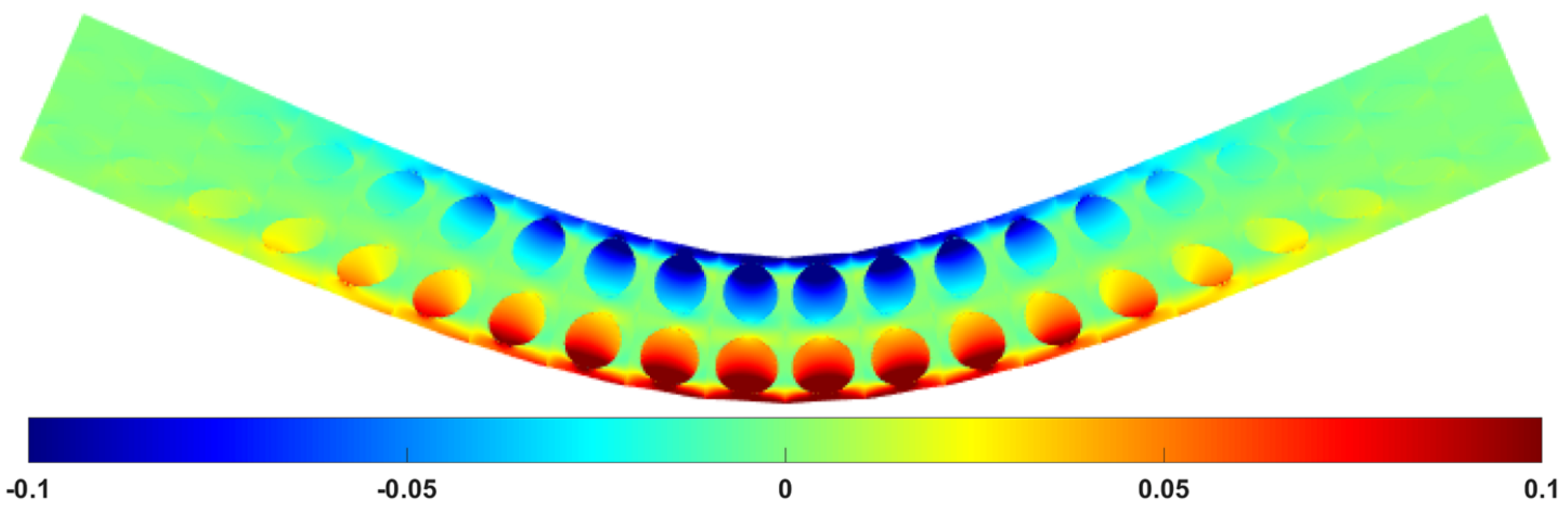}}
  \subfigure[Benchmark, $\varepsilon_{22}$]{\includegraphics[width=0.75\textwidth]{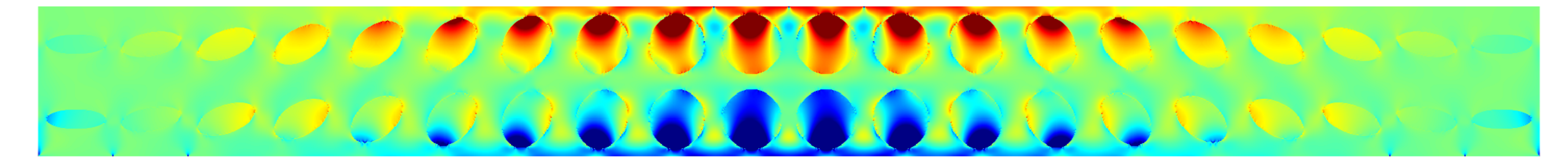}}
  \subfigure[Our-CBN, $\varepsilon_{22}$]{\includegraphics[width=0.75\textwidth]{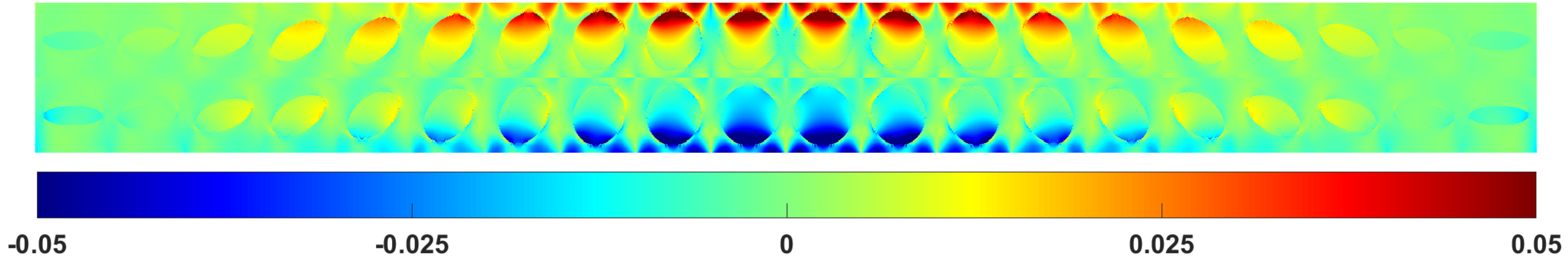}}
  \subfigure[Benchmark, $\varepsilon_{12}$]{\includegraphics[width=0.75\textwidth]{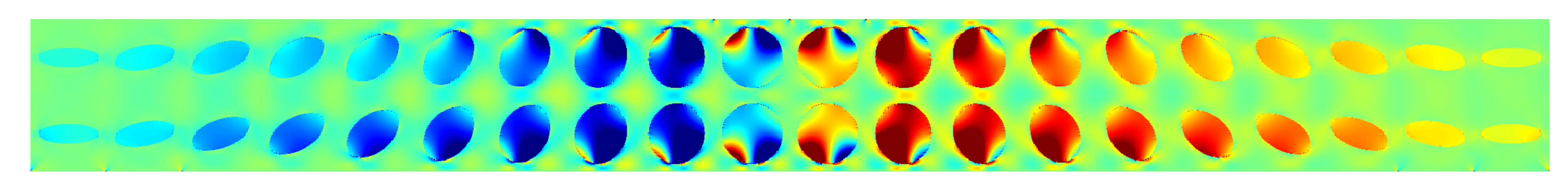}}
  \subfigure[Our-CBN, $\varepsilon_{12}$]{\includegraphics[width=0.75\textwidth]{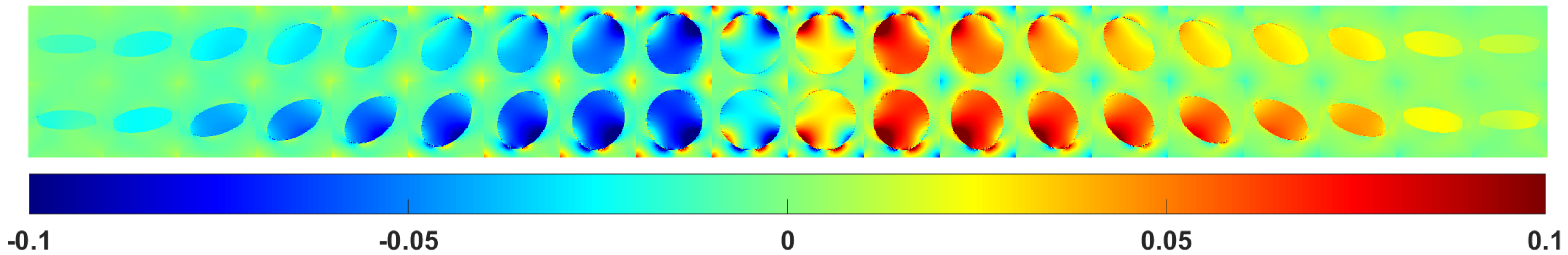}}
  \caption{Numerical results of Our-CBN approach in comparison with the benchmarks.}
  \label{fig:bending_strain}
\end{figure}

A more complex 2D heterogeneous bending beam example in Fig.~\ref{fig:domain_bending}(a) is also taken to test Our-CBN's performance in case of different coarse mesh sizes. The long beading beam is centered at $(0,0)$ in its left-bottom corner and has a size of $2\times 20$. At three different locations $(x^c_i,y^c_i)$ of coordinates of $(1, 0),(10, 2), (19, 0)$ from left to right, the body is imposed by a corresponding loading field in its vicinity, which described as
\eb
p(x) = p_i(1-(x - x^c_i)^2).
\ee
where $p_1 = p_2 = 1$ and $p_3 = 2$. The beam is fixed at the y-displacements on locations (0, 0) and (20, 0), and at the x-displacement on location (0, 2).

This example was modified from~\cite{le2020coarse} with two main changes. First, it contains elliptic holes of varied shapes, instead of homogeneous circular holes of the same shape and size. Second, the inclusion is softer, while~\cite{le2020coarse} has stiffer fibers; the former tends to produce a large deformation. In such a case, it is more challenging to produce a highly accurate result.

The global fine mesh was set at a fixed size of $200\times 2000$, and four coarse meshes of different sizes were set: $1\times 10$, $2\times 20$, $4\times 40$, and $8\times 80$. The effectivity indices $r_u,r_e$ are summarized in Table~\ref{tab:macro_mesh}, and variations of $r_u$ are also plotted in Fig.~\ref{fig:curve_bending} for a view.
The effectivity indices generally tend to decrease rapidly when the coarse element number increases (producing a small-sized local mesh).
This phenomenon can be explained via two observations. First, the shape functions tend to capture finer material variations for a local fine mesh of a smaller size.
Secondly, the global displacement tends to achieve a higher accuracy at the larger number of coarse DOFs.

Three exceptions are also observed. First, homogenization at $8\times 80$, where the coarse mesh displays a highly unordered distribution that strongly breaks the scale separation assumption. Second, FE$^2$ at $1\times 10$, which may come from the non-convergence of its nonlinear iteration. Third, CMCM at $4\times 40$, where a high material contrast along the coarse boundary poses a more challenging analysis task~\cite{le2020coarse}. Our-CBN handles all these situations well, and it has the smallest effectivity index of $r_e = 6e^{-3}$ and $r_u = 5e^{-3}$.

We also plot Fig.~\ref{fig:bending_strain} the beam's deformation and strain fields at a coarse mesh of size $2\times 20$, compared with the benchmark's. The computed results using Our-CBN are remarkably close to the benchmark's result, even in the middle large deformation area. Meanwhile, some local inconsistencies are also observed, particularly across the coarse elements' interfaces; similar phenomena were also observed in previous studies~\cite{le2020coarse}. Future research efforts are recommended to address the abovementioned issue.

\textbf{Comment.} The above results are built on the usage of the same size of coarse mesh, in which case we also notice that Our-CBN has more DOFs than the other approaches. For example, for a size $2\times 20$, Our-CBN has a DOF of $534$ (depending on the number of CBNs) while all the other approaches have a DOF of $126$ (depending on the corner node number).
To explore the topic deeply and more fairly, we further observe the other approaches' performance for the size of $8\times 80$ with $1458$ DOFs, around three times that of Our-CBN ($534$). However, using Our-CNB still approximately achieves an order of accuracy improvement. These observations indicate that Our-CBN has an intrinsic flexibility in closely capturing the coarse element's heterogeneity, which greatly improves its potentiality in the analysis of heterogeneous structures of non-separated scales. Its flexibility in choosing different DOFs further pronounces such potentialities.

\begin{figure}[p]
  \centering
  \subfigure[Structure of $2 \times 2 \times 8$ coarse mesh, each containing a $10\times 10 \times 10$ local fine mesh with $3 \times 3$ bridge nodes on one face]{\includegraphics[width=0.6\textwidth]{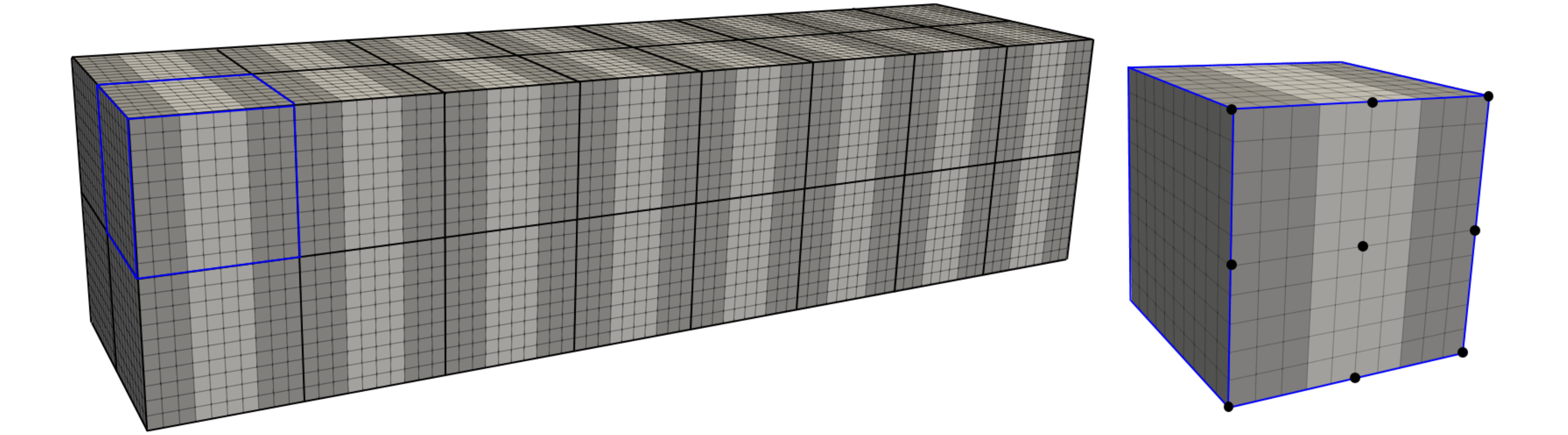}}\\
  \subfigure[Stretching]{\includegraphics[width=0.22\textwidth]{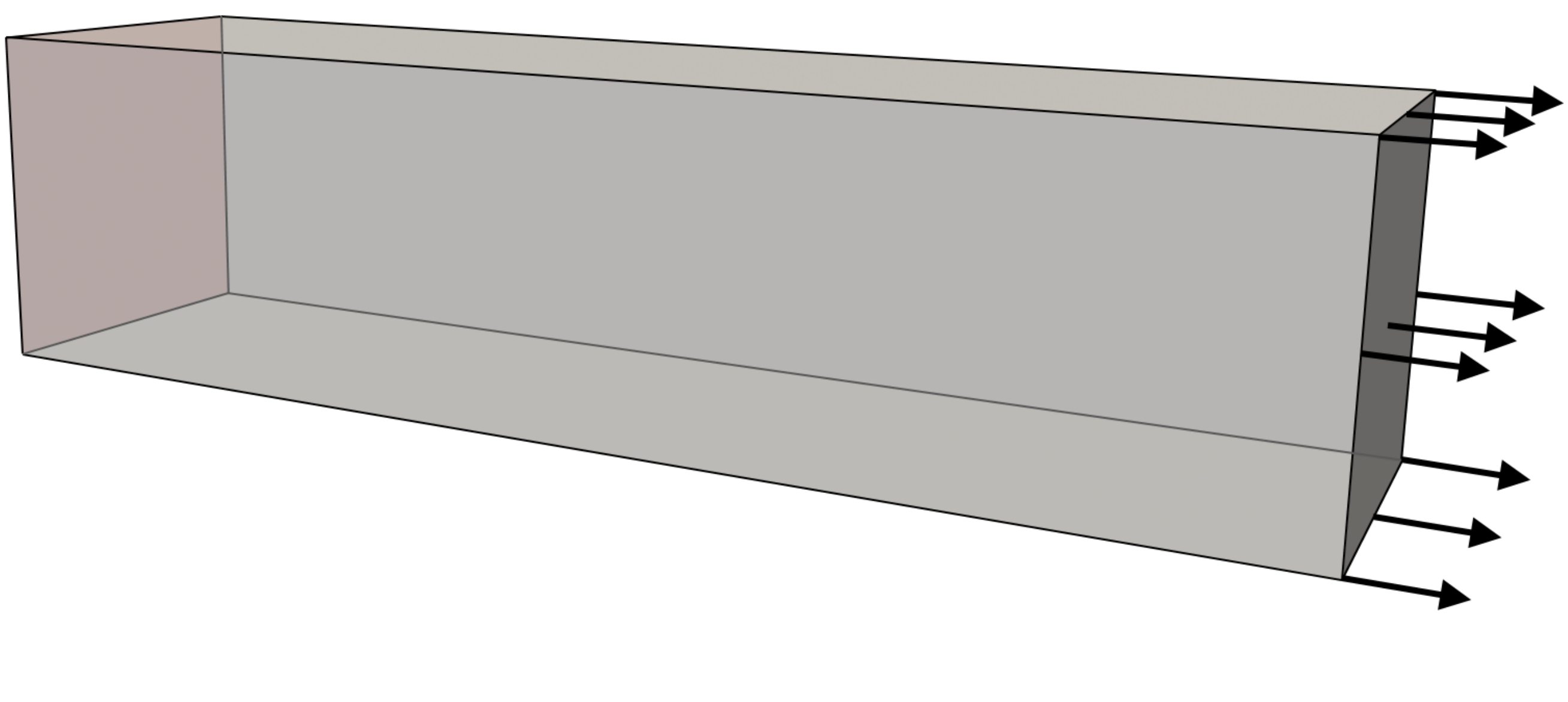}}
  \subfigure[Compressing]{\includegraphics[width=0.22\textwidth]{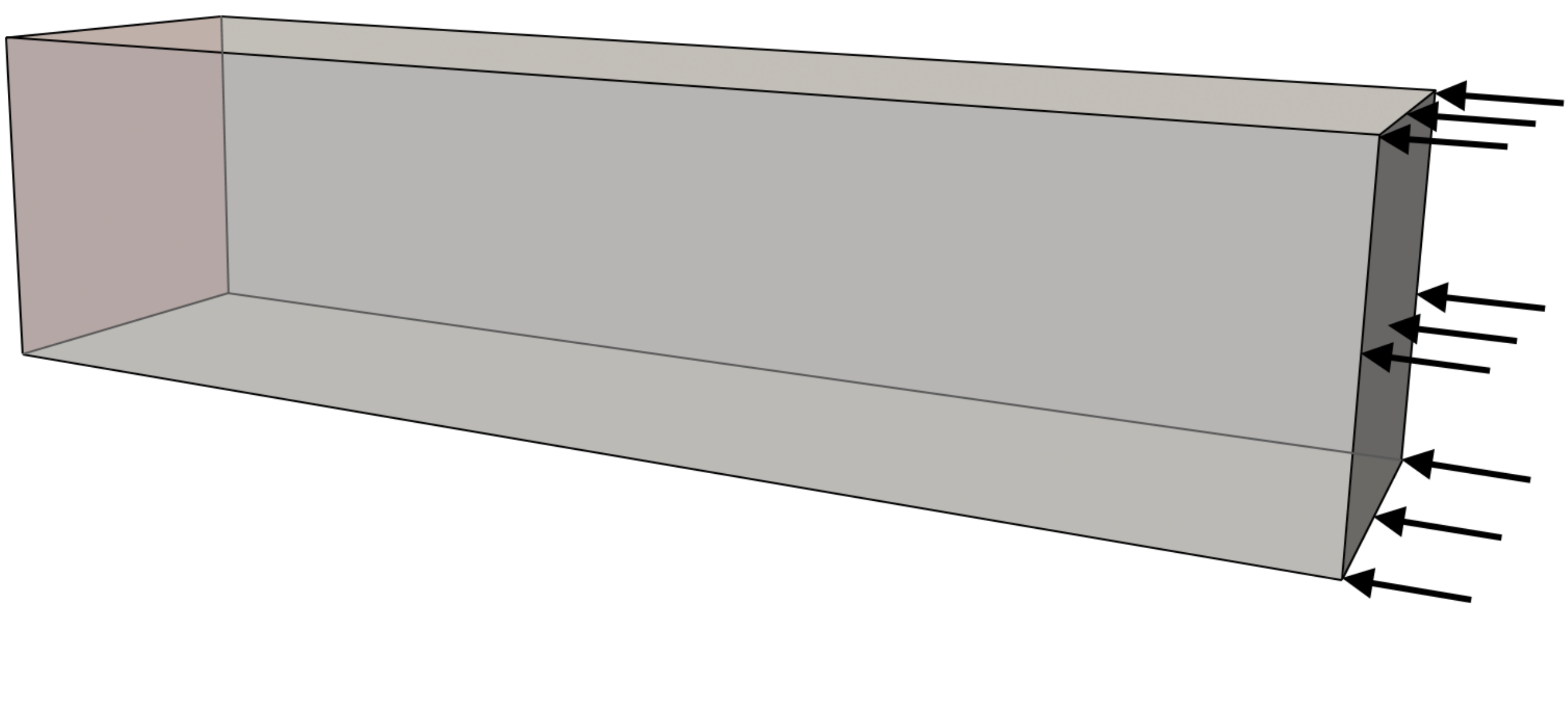}}
  \subfigure[Twisting]{\includegraphics[width=0.25\textwidth]{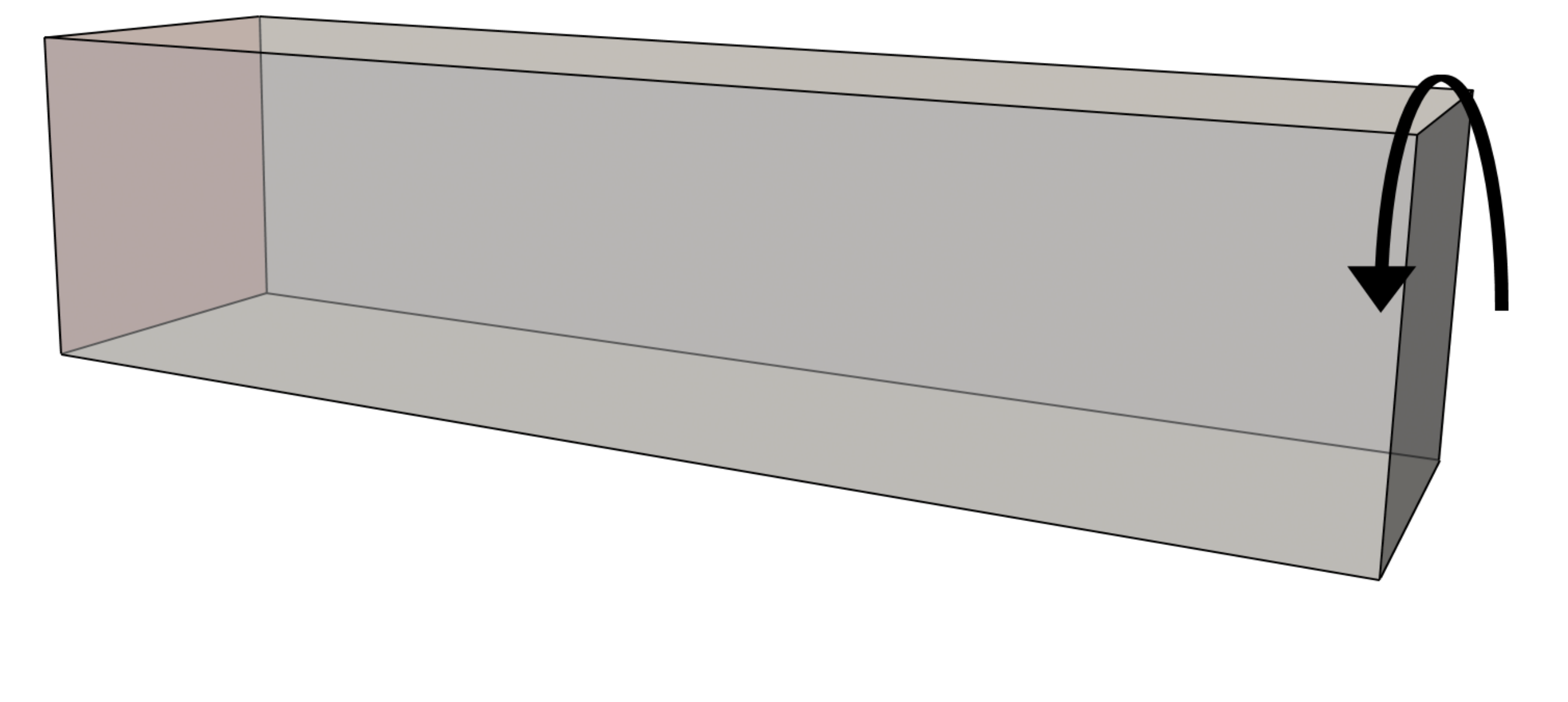}}
  \subfigure[Bending]{\includegraphics[width=0.25\textwidth]{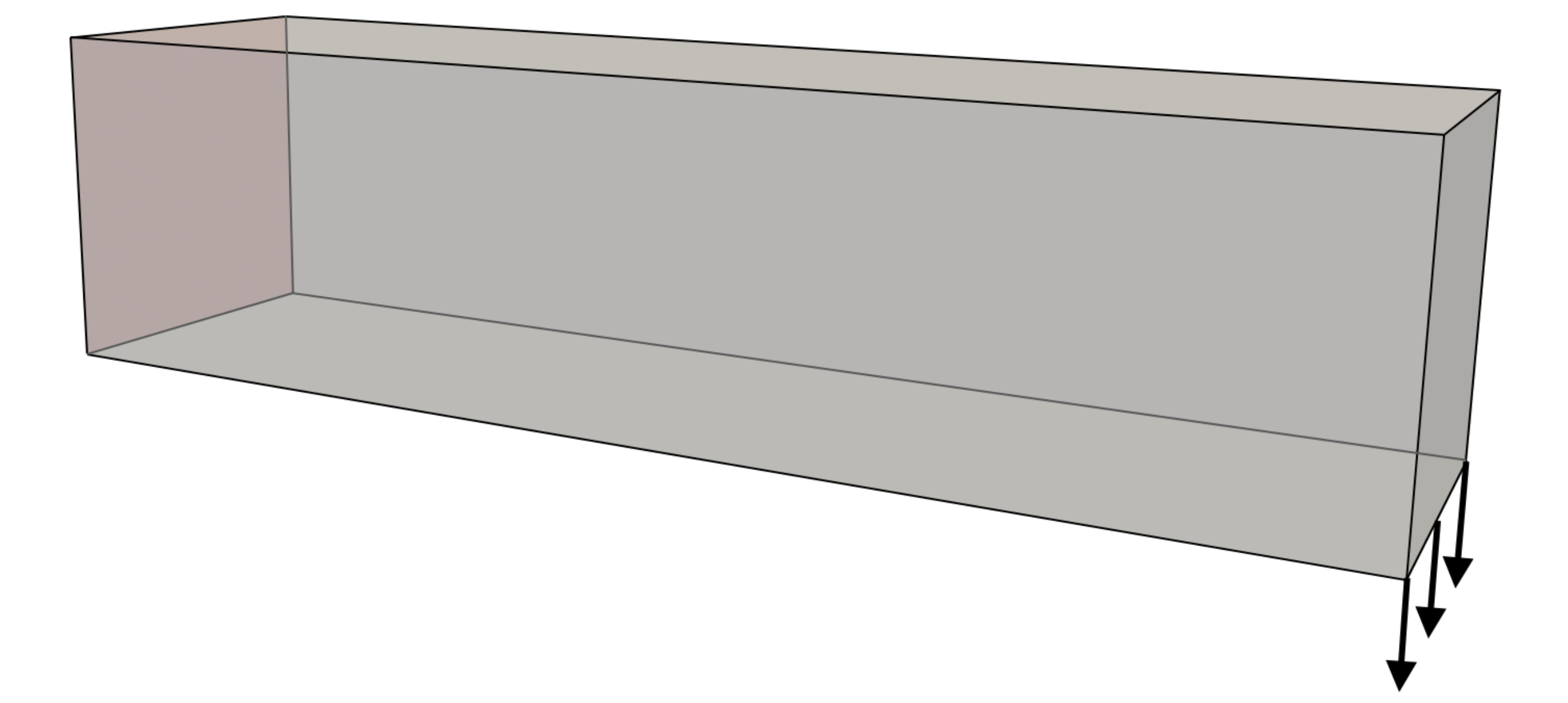}}
  \caption{A 3D example under four different loading conditions: stretching, compressing, twisting, and bending, where the red shadows denote the fixed area, and the black shadows and arrows denote the loading forces.}
  \label{fig:3d_domain}
%
  \centering
  \subfigure[Stretching, $r_e = 2.4e^{-4}, r_u = 8.9e^{-4}$]{\includegraphics[width=0.45\textwidth]{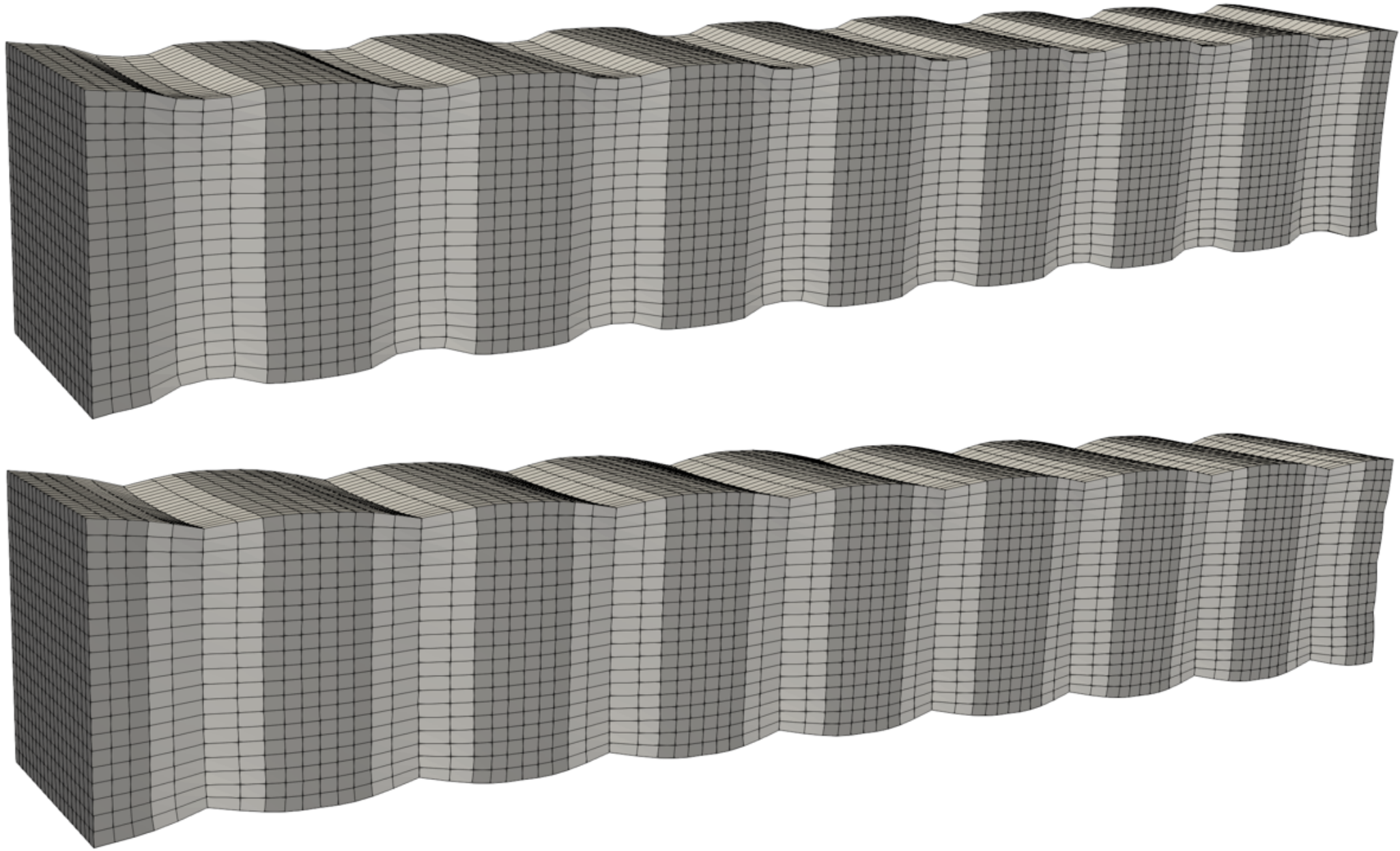}} \quad
  \subfigure[Compressing, $r_e = 6.5e^{-4}, r_u = 8.9e^{-4}$]{\includegraphics[width=0.3\textwidth]{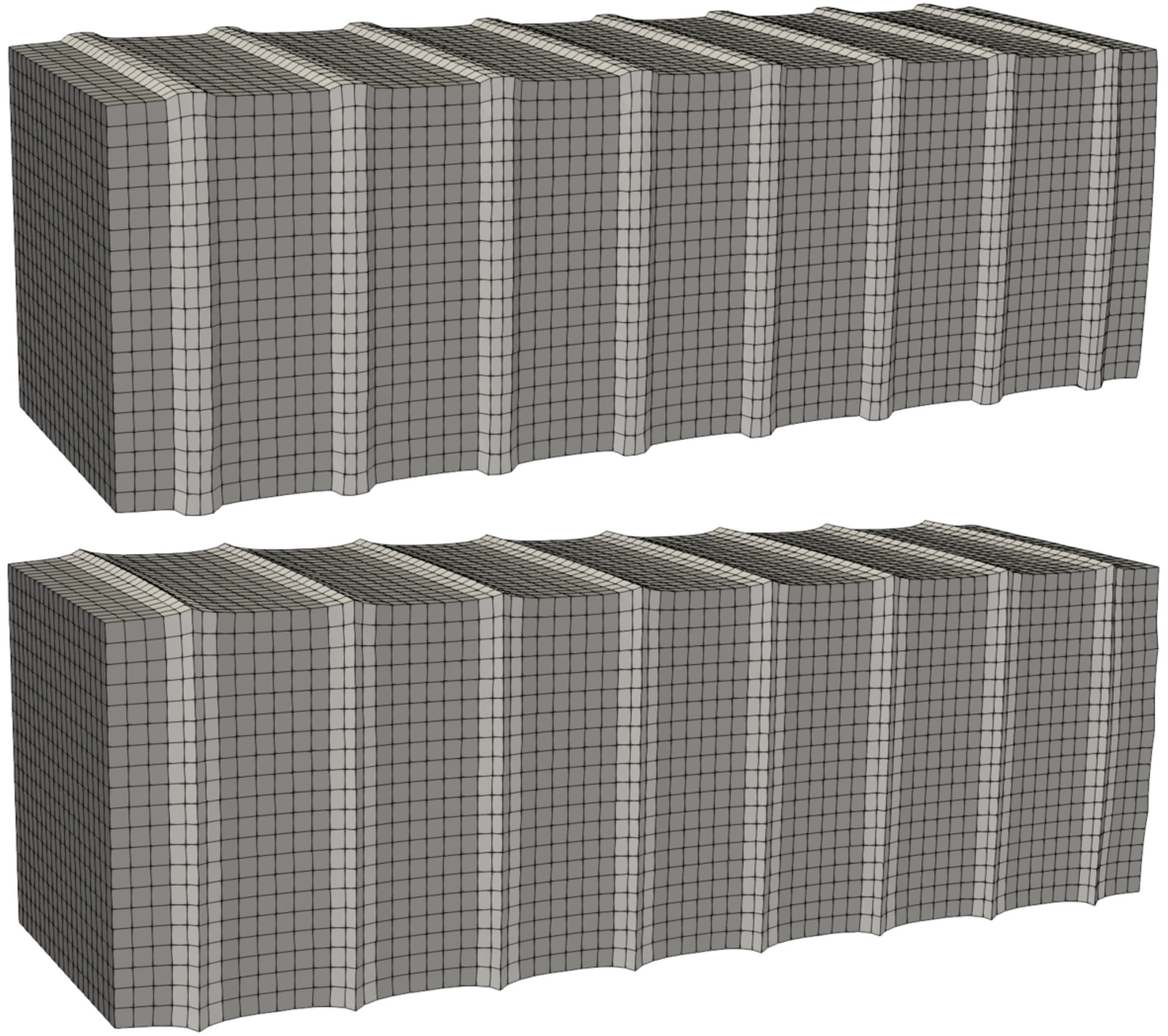}}\\
  \subfigure[Twisting, $r_e = 1.8e^{-4}, r_u = 3.1e^{-4}$]{\includegraphics[width=0.8\textwidth]{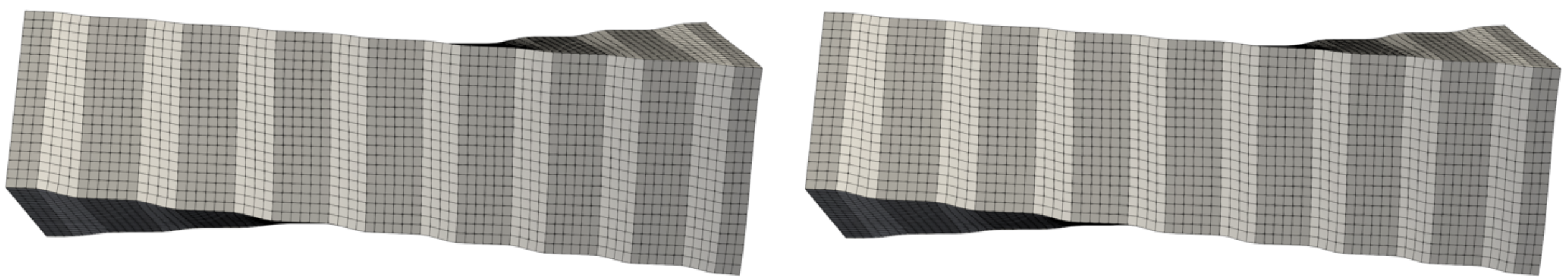}}\\
  \subfigure[Bending, $r_e = 1.4e^{-3}, r_u = 1.5e^{-3}$]{\includegraphics[width=0.8\textwidth]{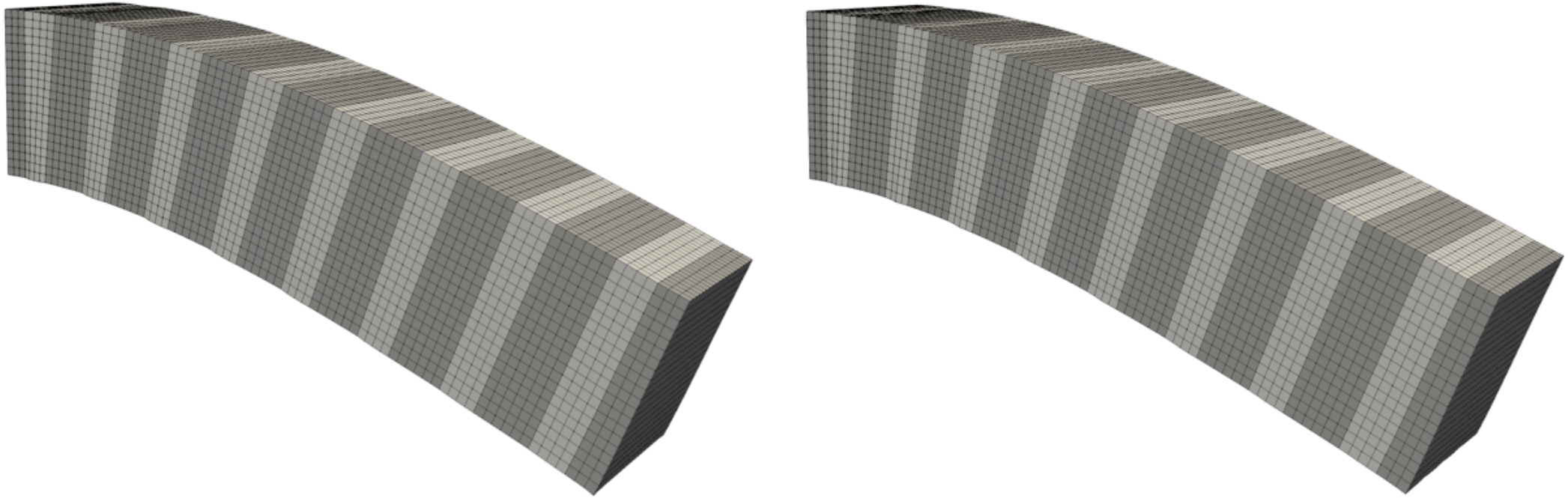}}
  \caption{Numerical results and the effectivity indices for the 3D example in Fig.~\ref{fig:3d_domain} compared with the benchmarks, where the benchmarks in (a),(b) are on the top and in (c),(d) are on the left.}
  \label{fig:3d_results}
\end{figure}

\subsection{A heterogeneous 3D example under different loading conditions}\label{sec:3dexample}
We also test Our-CBN's performance on the 3D example in Fig.~\ref{fig:3d_domain}, where the dark and light areas have Young's modulus of $E = 1e^4$ and $1e^3$. The model has a coarse mesh of size $2\times 2 \times 8$ and a local fine mesh of size $10\times 10 \times 10$, and each coarse mesh face has a set of $3\times 3$ coarse bridge nodes.

The example is tested under four classical loading cases: stretching, compressing, twisting, and bending, and the results are shown in Fig.~\ref{fig:3d_results} in comparison with the benchmark results. The largest effectivity indices in the tests have values of $r_e = 1.4e^{-3}$ and $r_u = 1.5e^{-3}$, demonstrating the approach's high approximation accuracy.
Note that in this analysis, all the coarse elements have the same heterogeneity distribution. Using Our-CBN, the shape functions can only be affordably computed offline once for a single coarse element, irrespective of these different loading conditions.

\begin{figure}[tbh]
    \centering
    \subfigure[Material distribution]{\includegraphics[width=0.48\textwidth]{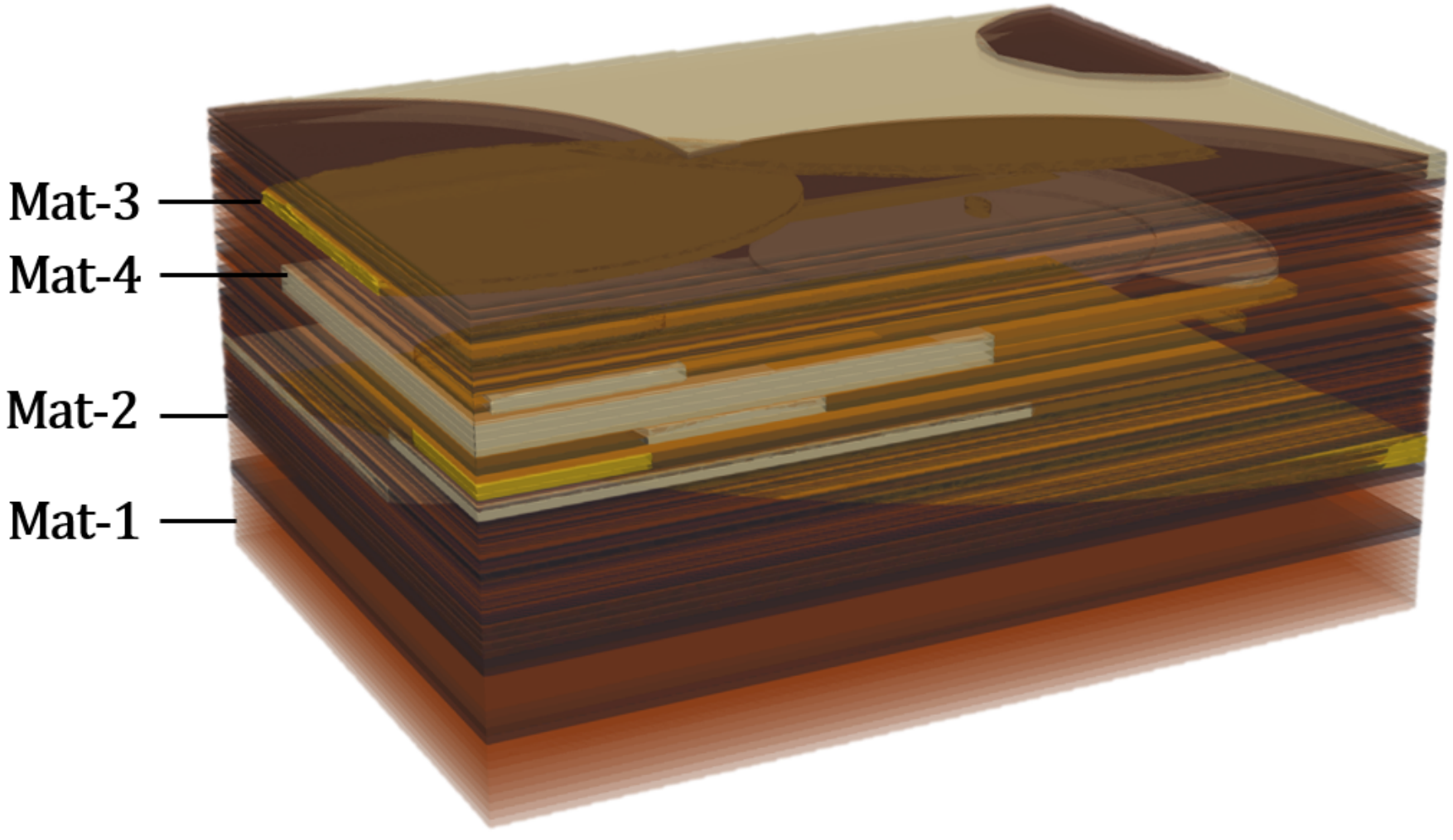}}\quad
    \subfigure[A bending example]{\includegraphics[width=0.4\textwidth]{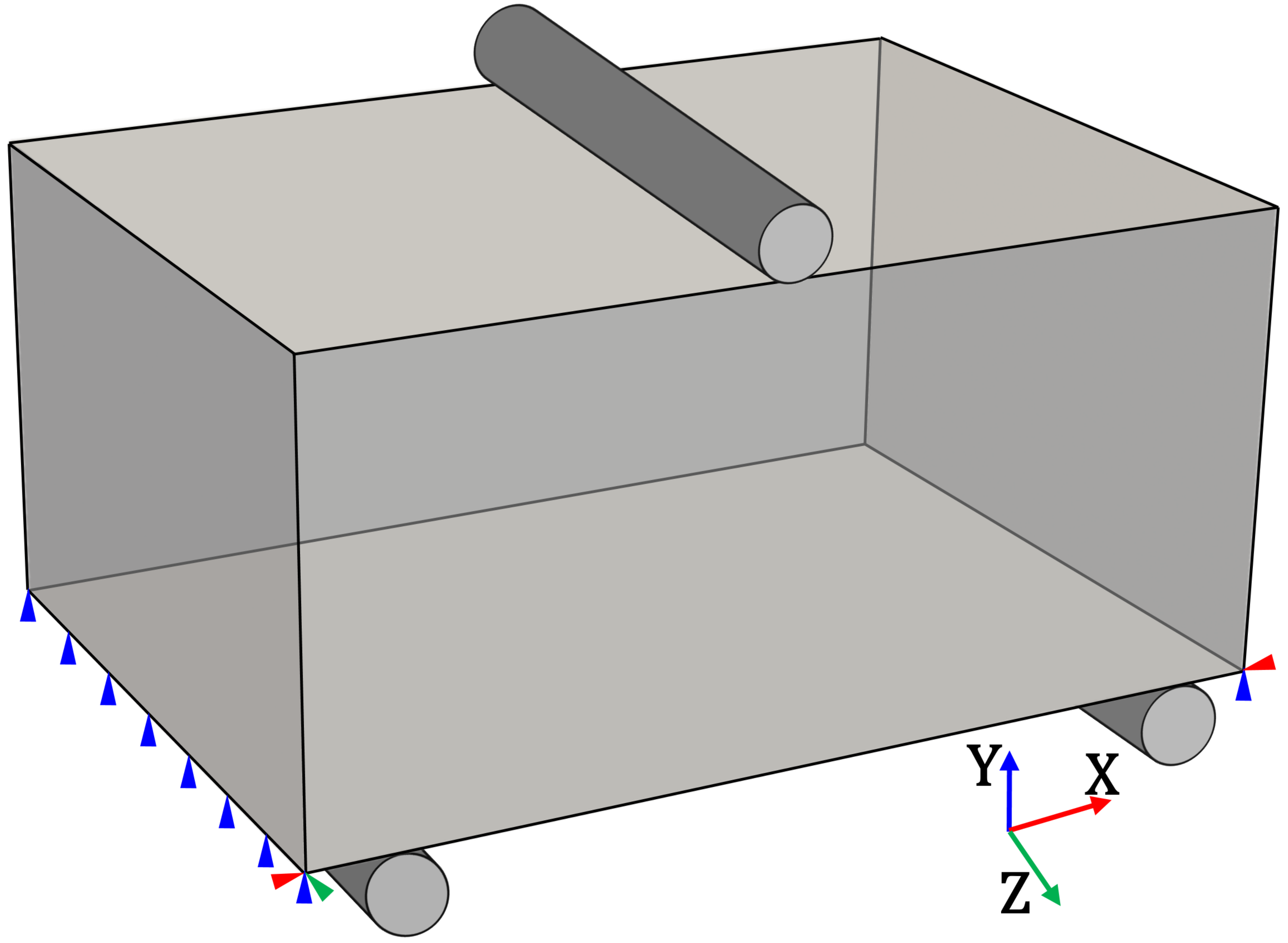}}
    \caption{An industrial 3D geologic model of 2.05 billion DOFs to test our approach's performance in heterogeneous structure analysis.}
    \label{fig:geologic_domain}
 \end{figure}

\begin{figure}[p]
 \centering
  \subfigure[Global deformation of $\calM^H$ in different views]{\includegraphics[width=0.9\textwidth]{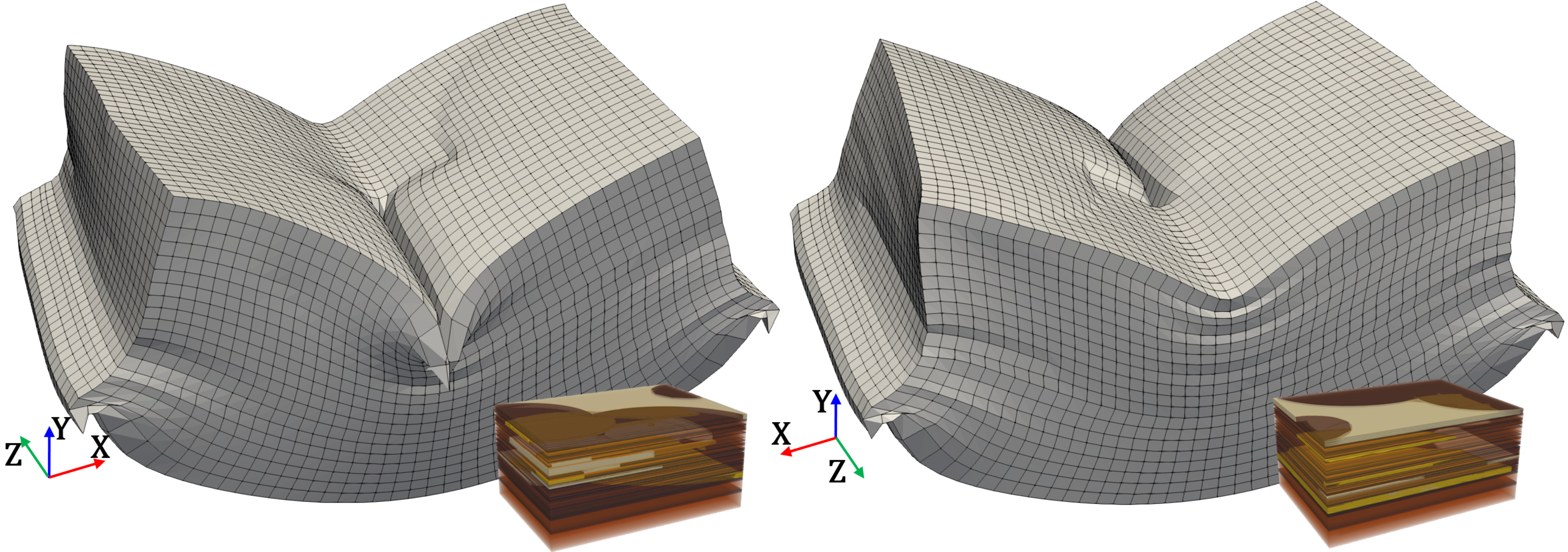}}
  \subfigure[Deformation of top layers in overall view]{\includegraphics[width=0.9\textwidth]{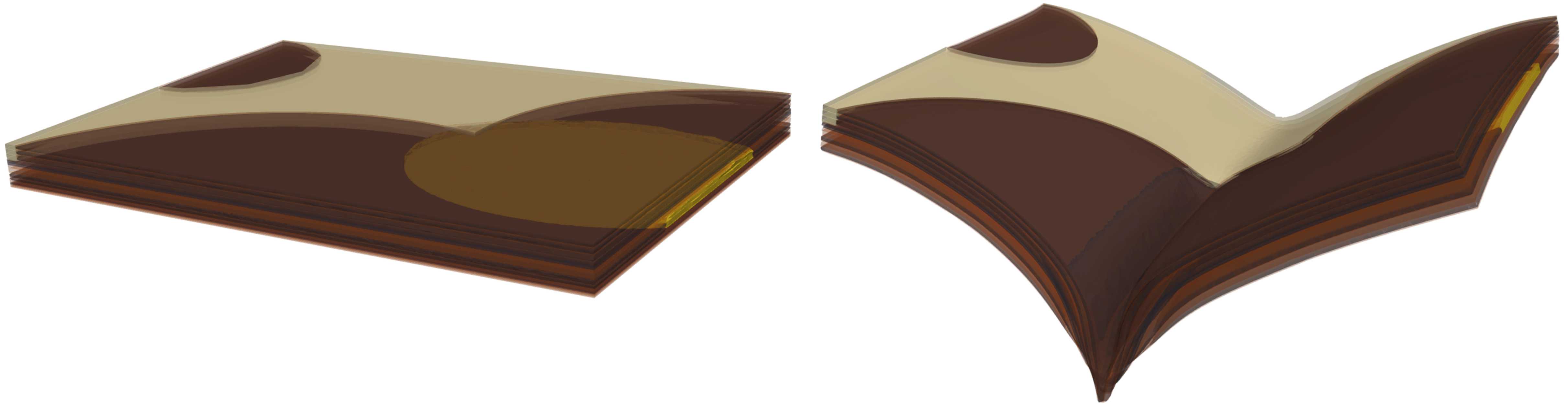}}\quad
  \subfigure[Deformations of top layers in isolation]{\includegraphics[width=0.9\textwidth]{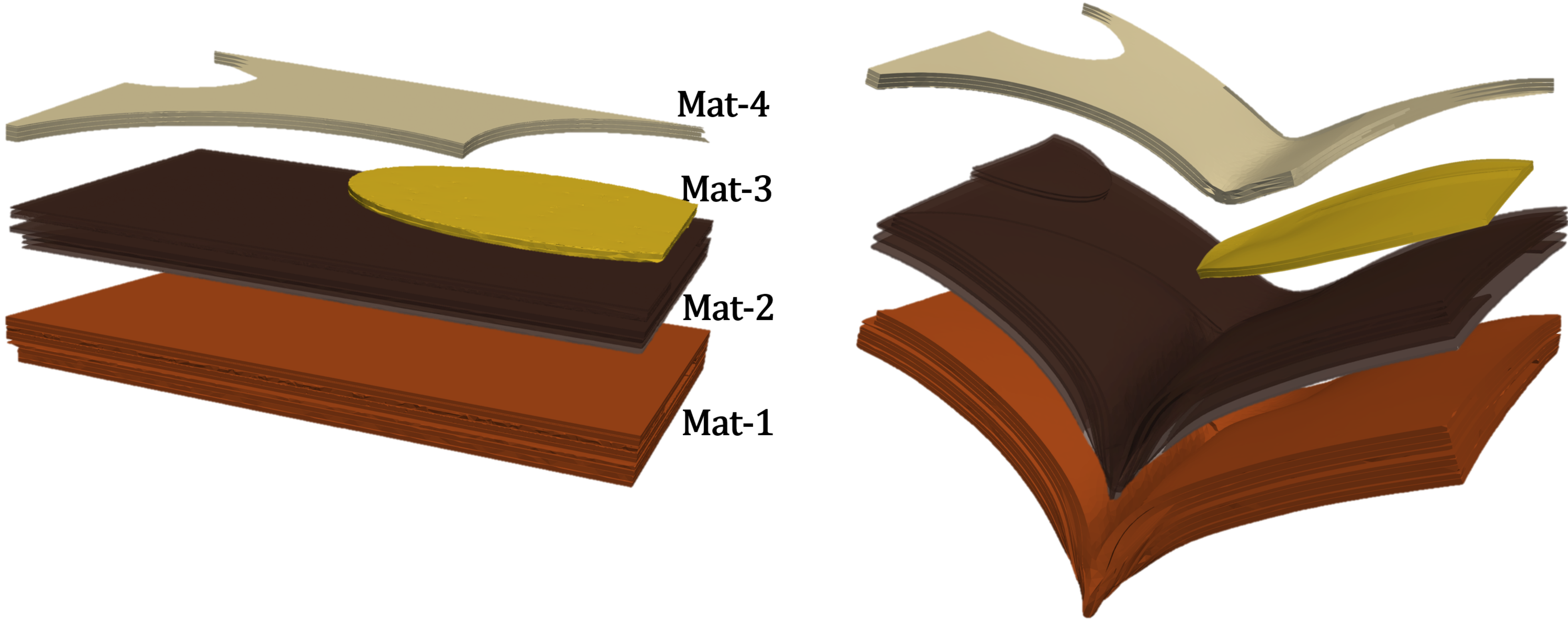}}
  \subfigure[Local deformation of four different $\calM^{\alpha,h}$]{\includegraphics[width=0.6\textwidth]{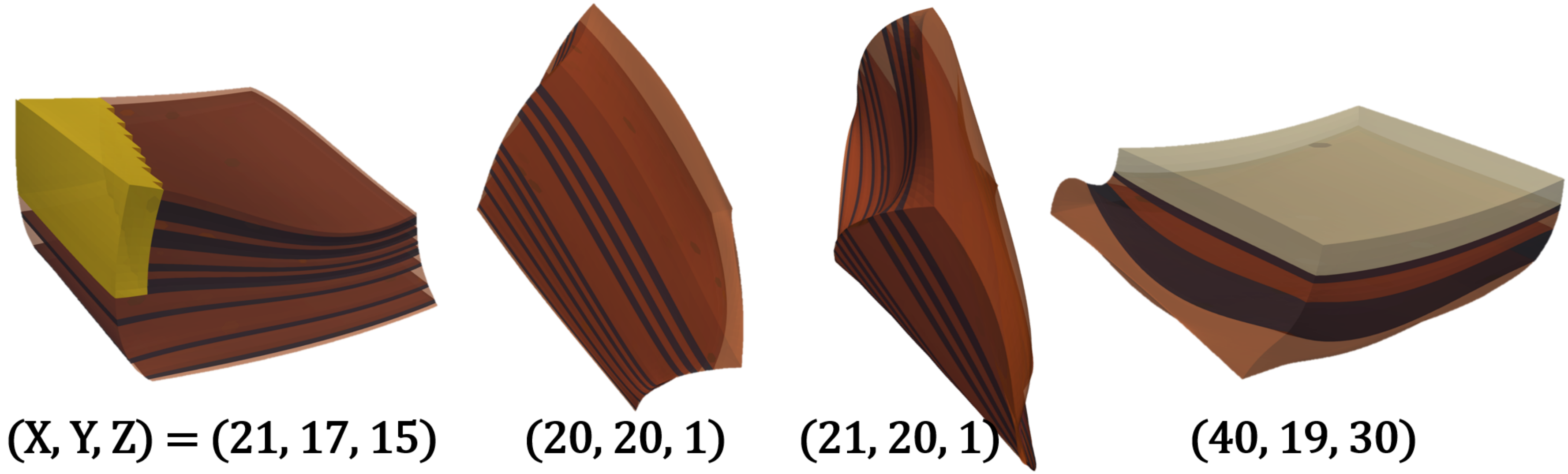}}
  \caption{The numerical results on the geologic model in Fig.~\ref{fig:geologic_domain}, where constructing the shape functions takes $6.6$ s per coarse element, and computing the global deformation takes $1.2$ h.}
  \label{fig:recover_dizhi}
\end{figure}

\subsection{A practical 3D large-scale geologic model with $2.05$ billion DOFs}
\label{sec:geologic}

A modified industrial 3D large-scale complex geologic model in Fig.~\ref{fig:geologic_domain}(a) is also analyzed using Our-CBN. The model contains four types of materials in different colors: miscellaneous fill (brown), silty clay (coffee), strong weathering rock (yellow), and middle weathering rock (buff); the property parameters are listed in Table.~\ref{tab:soil}. The model is fixed on its left and right sides, and imposed by three different pressure fields with $P=1e^5$ induced by the contact cylinders on its top or bottom, as shown in Fig.~\ref{fig:geologic_domain}(b).

In the analysis, we have a coarse mesh of size $20\times 30 \times 42$ and a local fine mesh of size $30\times 30 \times 30$, which combined give a global fine mesh of size $600\times 900 \times 1260$, involving approximately $2.05$ billion DOFs. The corner nodes are set as the bridge nodes. This turns out an offline local analysis problem of $89$-thousand DOFs, and an online global analysis problem of $1.5$ million DOFs. The local computation takes $6.6$ s per coarse element, and the global computation takes $1.2$ h using the conjugate gradient (CG) method ended with a relative residual $9.5e^{-7}$ in $8658$ iterations.

Fig.~\ref{fig:recover_dizhi} plots the overall deformations in (a), and the deformations of some top layers ($1/6$ of the whole height) in (b) and (c). A close-up of the deformation in some heterogeneous coarse elements are also shown in (d).
Different deformation behaviors are observed in these different regions---the softer region demonstrates relatively large deformations (in brown and coffee) while the stiffer region shows small deformations (in buff and yellow). Such phenomena are even observed in different areas in a single coarse element, as shown in (d), demonstrating the approach's ability in describing finely detailed local deformations of a heterogenous structure.

\subsection{Extension to nonlinear elastic model}
\label{sec:nonlinear_result}
\begin{figure}[htb]
  \centering
  \subfigure[Benchmark]{\includegraphics[width=0.4\textwidth]{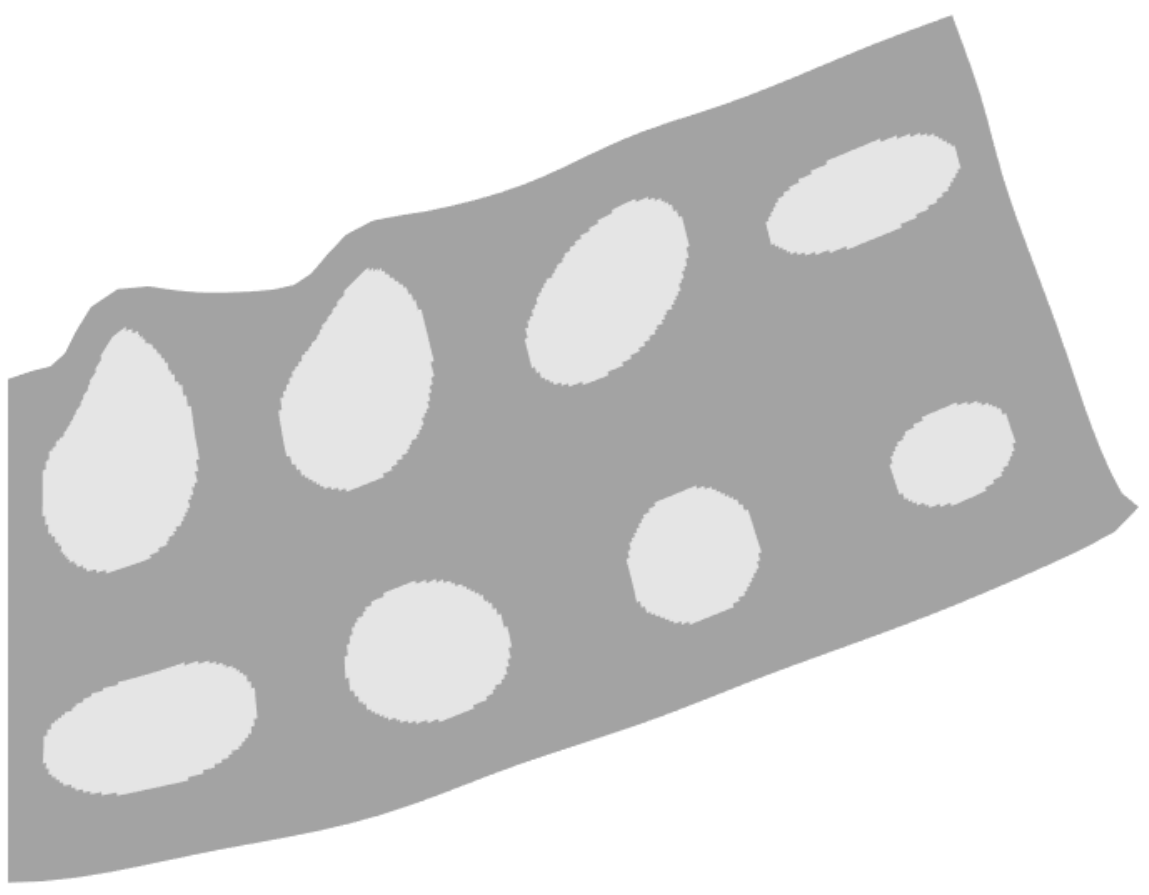}}\quad
  \subfigure[Our results, $r_e = 2.5e^{-6}$, $r_u = 0.07$]{\includegraphics[width=0.4\textwidth]{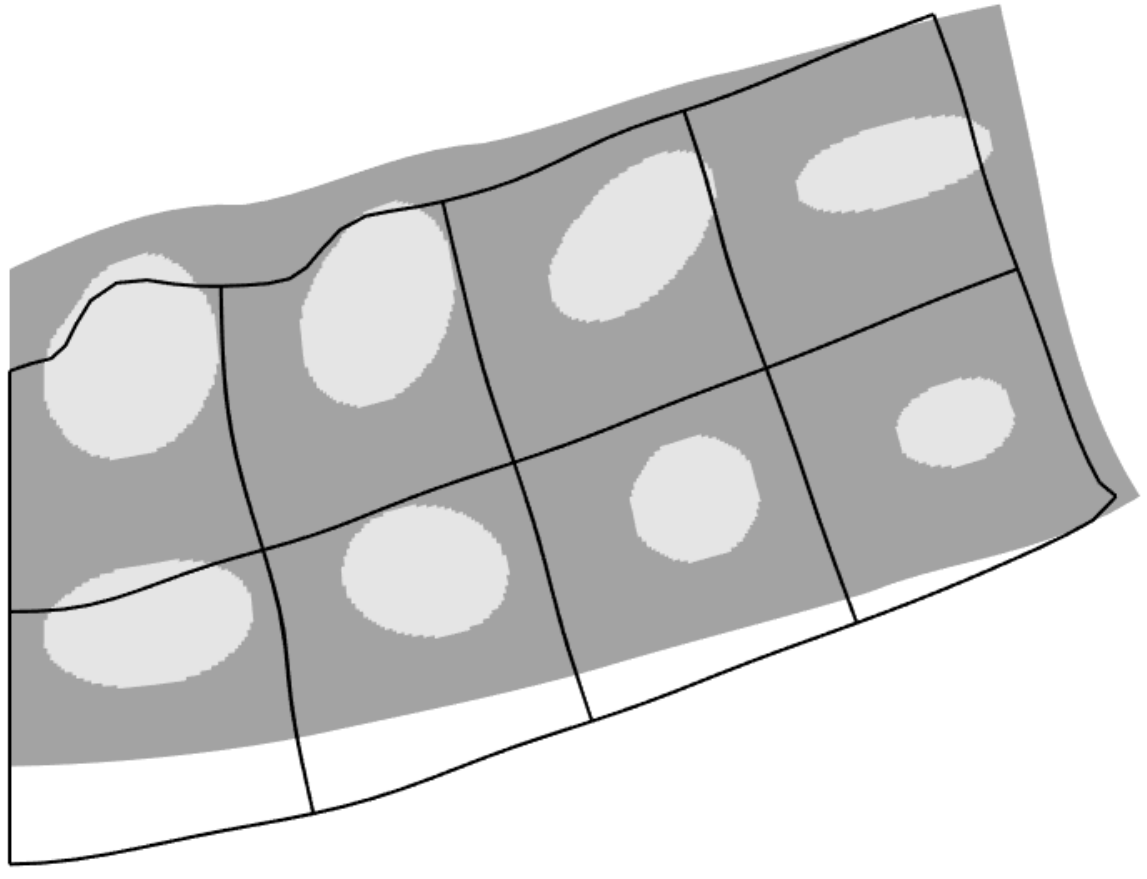}} \\
  \caption{Numerical results for half MBB in Fig.~\ref{fig:MBB} of nonlinear neo-Hookean material, where the black outlines in (b) denote the benchmark deformation.}
  \label{fig:nonlinear_2D}
 \end{figure}

Our CBN approach also works for the analysis of the nonlinear elastic model. Consider the half MBB in Fig.~\ref{fig:MBB} of neo-Hookean materials at a loading of $10$. The computed deformation is plotted in Fig.~\ref{fig:nonlinear_2D}, in comparison with the benchmark.

Unlike the linear case, $r_e$ and $r_u$ have very different values: $r_e = 2.5e^{-6}$ and $r_u = 0.07$. The small value of $r_e$ indicates that we have reached a global deformation energy, which is the same as that of the benchmark. However, pseudo-stiffness still exists as indicated by the large quantity $r_u = 0.07$. This is believed to be caused by the local linear elasticity analysis in constructing our CBN shape functions. Employing a nonlinear model to build more advanced shape function seems to be a reasonable choice for future research exploration.

\begin{table}[t]
  \centering
  \caption{Mechanical properties of different soil layers in the geologic model in Fig.~\ref{fig:geologic_domain}.}
  \label{tab:soil}
  \begin{tabular}{c|c|c|c}
  \hline
  \rowcolor[HTML]{C0C0C0}
  Index & Soil layer             & Bulk modulus K & Shear modulus G \\ \hline
  Mat-1 (brown) & Miscellaneous fill     & 7e$^6$       & 3.2e$^6$      \\
  \rowcolor[HTML]{EFEFEF}
  Mat-2 (coffee) & Silty clay             & 1.86e$^7$    & 9e$^6$        \\
  Mat-3 (yellow) & Strong weathering rock & 1.38e$^8$    & 5.96$e^7$     \\
  \rowcolor[HTML]{EFEFEF}
  Mat-4 (buff) & Middle weathering rock & 6.3e$^8$     & 3.86e$^8$     \\ \hline
  \end{tabular}
\end{table}

\section{Conclusions}
This study introduces the concept of curved bridge nodes (CBNs) and its associated CBN shape functions for the elasticity analysis of heterogeneous structures of non-separated scales. The shape functions are derived per coarse element as a product of a B\'ezier interpolation transformation and boundary-interior transformation and result in shape functions in an explicit matrix representation.

The B\'ezier interpolation transformation not only ensures the displacement smoothness between adjacent coarse elements but also provides additional variables in reducing the problem of inter-element stiffness. The boundary-interior transformation, derived from the local stiffness matrix to the local fine mesh, provides a prominent advantage to finely embed the intrinsic material heterogeneity into the shape functions. Finally, the derived shape functions have the properties of basic FE shape functions that avoid aphysical behavior. Extensive numerical examples indicate that Our-CBN has an intrinsic flexibility in closely capturing the coarse element's heterogeneity, and it may serve as a suitable method for the analysis and optimization of heterogeneous structures without scale separation  ~\cite{xia2014concurrent,wu2019topology,gao2019topology}.

Furthermore, the CBN shape functions directly work for nonlinear elasticity analysis problems but may encounter accuracy loss. Introducing nonlinear analysis in the shape function construction appears to be a very promising approach for the improvement of its analysis accuracy and warrants further research efforts. In addition, the shape functions can be computed in parallel, which boosts their applications in analysis of super-large problems, although their achievement still depends on the availability of appropriate computational facilities. Developing a surrogate model by using techniques on model reduction~\cite{el2013fe2,wu2019topology} or deep learning~\cite{2020Model,2021Iterative} is expected to help resolve the existing limitations and should be explored in future studies.

\section*{Acknowledgements}
The study described in this paper is partially supported by the National Key Research and Development Program of China (No. 2018YFB1700603) and the NSF of China (No. 61872320).

\bibliographystyle{elsarticle-num}
\bibliography{simulation}
\end{document}